\documentclass[a4paper,11pt]{memoir}

\usepackage[T1]{fontenc}
\usepackage[cp1250]{inputenc}

\usepackage[leqno]{amsmath}
\usepackage{amsfonts}
\usepackage{amssymb}
\usepackage{amsthm}
\usepackage{mathrsfs}

\usepackage{enumerate}
\usepackage{setspace}

\checkandfixthelayout

\makechapterstyle{MyPhdChapter}{%

}
\chapterstyle{MyPhdChapter}

\makepagestyle{titlepage} %define new pagestyle
\makeoddhead{titlepage}{}{{\Large\MakeUppercase{Institute of Mathematics Polish Academy of Sciences}}}{}
\makeoddfoot{titlepage}{}{\Large Łódź, 2013}{}
\setheadfoot{14.pt}{\baselineskip} %set header and footer height

\setsecnumdepth{subsection}

\maxtocdepth{subsection}

\newcommand{\al}{\alpha}
\newcommand{\be}{\beta}
\newcommand{\ga}{\gamma}
\newcommand{\del}{\delta}
\newcommand{\Del}{\Delta}
\newcommand{\eps}{\varepsilon}
\newcommand{\lam}{\lambda}
\newcommand{\si}{\sigma}

\newcommand{\Om}{\Omega}
\newcommand{\cl}[1]{\overline{#1}}
\newcommand{\inter}{\operatorname{int}}
\newcommand{\oline}[1]{\overline{#1}}
\newcommand{\uline}[1]{\underline{#1}}
\newcommand{\bsh}{\backslash}

\newcommand{\re}{\text{Re}}

\newcommand{\power}[1]{2^{#1}} %power set of given #1
\newcommand{\bbC}[1][]{\mathbb{C}^{#1}}
\newcommand{\bbR}[1][]{\mathbb{R}^{#1}}
\newcommand{\bbN}[1][]{\mathbb{N}^{#1}}
\newcommand{\bbNo}[1][]{\mathbb{N}^{#1}_0}
\newcommand{\bbZ}[1][]{\mathbb{Z}^{#1}}
\renewcommand{\H}{H} %main Hilbert space
\newcommand{\reH}{H^{(r)}} %real Hilbert space obtained from H
\newcommand{\X}{X} %Banach space
\newcommand{\weak}{\rightharpoonup}

\newcommand{\linspan}{\text{span}}
\newcommand{\conv}{\text{conv}}
\newcommand{\LR}[1][]{L^2(\mathbb{R}^{#1}\!,\bbR)}
\newcommand{\LC}[1][]{L^2(\mathbb{R}^{#1}\!,\bbC)}

\newcommand{\LL}{L} %main linear operator
\newcommand{\N}{N} %main nonlinear operator
\newcommand{\R}[1]{R_{#1}} %resolvent operator of #1
\newcommand{\dom}{\operatorname{Dom}}
\newcommand{\ran}{\operatorname{Ran}}
\renewcommand{\ker}{\operatorname{Ker}}
\newcommand{\graph}{\operatorname{Gra}}
\newcommand{\ind}{\text{ind}} %index
\newcommand{\slim}[1]{\underset{#1}{s\text{-}\!\lim}\,}
\newcommand{\lin}{\mathcal{L}} %general linear operators
\newcommand{\clos}{\mathcal{C}} %closed operators
\newcommand{\bound}{\mathcal{B}} %bounded operators
\newcommand{\proj}{\mathcal{P}} %orthoprojections
\renewcommand{\sp}{\si}
\newcommand{\psp}{\si_{p}} %point spectrum
\newcommand{\csp}{\si_{c}} %continuous spectrum
\newcommand{\rsp}{\si_{r}} %residual spctrum
\newcommand{\essp}{\si_{e}} %essential spectrum
\newcommand{\dsp}{\si_d} %discrete spectrum
\newcommand{\rez}{\rho} %residual set
\newcommand{\lH}{\H_1}
\newcommand{\rH}{\H_2}
\newcommand{\lP}{P_1}
\newcommand{\rP}{P_2}
\newcommand{\lL}{\LL_1}
\newcommand{\rL}{\LL_2}
\newcommand{\lN}{\N_1}
\newcommand{\rN}{\N_2}
\newtheorem{theorem}{Theorem}[chapter]
\newtheorem*{twierdzenie}{Twierdzenie}
\newtheorem{lemma}[theorem]{Lemma}

\theoremstyle{definition}
\newtheorem{definition}[theorem]{Definition}
\theoremstyle{remark}
\newtheorem{remark}[theorem]{Remark}
\let\oldproof\proof
\let\oldendproof\endproof
\newenvironment{myproof}[1][\proofname]{%
  \oldproof[\bfseries{#1}]%
}{\oldendproof}

\begin{document}

%\vspace*{1cm}

\begin{center}
{\Large\MakeUppercase{Doctor of Philosophy Dissertation}}
\end{center}
\vspace{2.5cm}

\begin{center}
{\Large \MakeUppercase{Przemysław Zieliński}}
\end{center}
\vspace{0.5cm}
\begin{center}
\textbf{\huge\DoubleSpacing\MakeUppercase{Spectral and topological}}\smallskip
\textbf{\huge\DoubleSpacing\MakeUppercase{methods in the study of}}\smallskip \textbf{\huge\DoubleSpacing\MakeUppercase{solvability of semilinear}}\smallskip
\textbf{\huge\DoubleSpacing\MakeUppercase{equations in Hilbert spaces}}
\end{center}\medskip

\begin{center}
{\Large\OnehalfSpacing\MakeUppercase{Metody spektralne i topologiczne w badaniu rozwiązalności równań semi-liniowych w przestrzeniach Hilberta}}
\end{center}
\vspace{0.5cm}
%\begin{center}
%{\LARGE \MakeUppercase{(Rough Copy)}}
%\end{center}
\vspace{3cm}
\noindent{\Large Supervised by\\ Professor Bogdan Przeradzki\\ Institute of Mathematics, Lodz University of Technology}

\thispagestyle{titlepage}
\clearpage

\section*{Podziękowania}

Serdeczne podziękowania składam Panu prof. dr hab. Bogdanowi Przeradzkiemu za poświęcony czas, dzielenie się wiedzą oraz cenne wskazówki, które ułatwiły mi pisanie tej pracy.\medskip

\noindent Dziękuję również moim najbliższym za pomoc, wparcie i wiarę we mnie.

\thispagestyle{empty}

\cleardoublepage

\frontmatter

%\begin{center}
%\textbf{\Large Opis pracy}
%\end{center}
%\vspace{1cm}
\chapter{Opis pracy}

\noindent\textit{Słowa kluczowe:} równania semi-liniowe, widmo istotne operatora liniowego, stopień topologiczny, operatory maksymalnie monotoniczne, operatory typu monotonicznego
\vspace{0.5cm}

Celem tej dysertacji jest podanie warunków dostatecznych rozwiązalności w przestrzeni Hilberta $\H$ równań semi-liniowych postaci
\begin{equation}\label{eq:glowne}\tag{$\star$}
\LL u+\N(u)=h,
\end{equation}
gdzie $\LL$ jest operatorem liniowym samosprzężonym (zwykle nieograniczonym), $\N$ jest nieliniowy oraz $h\in\H$. Podstawową nowością w realizowanym tu podejściu jest osłabienie założeń nakładanych na postać widma części liniowej $\LL$.\medskip

W większości zastosowań operatory liniowe fizyki matematycznej pochodzące od zagadnień brzegowych na zbiorach ograniczonych, posiadają zwartą rezolwentę. Oznacza to, po pierwsze, że widmo $\sp(\LL)$ operatora $\LL$ jest dyskretne, czyli składa się z izolowanych wartości własnych o skończonych krotnościach geometrycznych. Po drugie, dla dowolnego $\lambda\notin\sp(\LL)$ rezolwenta
\[
\R{\LL}(\lam)=(\LL-\lam I)^{-1},
\]
gdzie $I$ jest identycznością na $\H$, jest operatorem zwartym określonym na całym $\H$, czyli przeprowadza zbiory ograniczone w relatywnie zwarte. Zatem jeżeli $\LL$ jest operatorem ze zwartą rezolwentą, to w równaniu (\ref{eq:glowne}) możemy wyróżnić dokładnie dwa przypadki. Jeśli $0\notin\sp(\LL)$, to wyjściowe równanie możemy zamienić na zagadnienie poszukiwania punktów stałych
\[
u=\R{\LL}(0)(h-N(u))
\]
i ze względu na zwartość rezolwenty stosować teorię Leray-Schaudera. W drugim przypadku mamy $0<\dim\ker\LL<+\infty$. Mówimy wtedy, że równanie (\ref{eq:glowne}) jest \emph{rezonansowe}. Wówczas częściowe odwrócenia $\LL$ na zakresie $\ran(\LL)$ są zwarte, więc możemy stosować teorię stopnia koincydencji Mawhina do poszukiwania rozwiązań (\ref{eq:glowne}) (\cite{KarakostasTsamatos2001}). Abstrakcyjne podejście jest przedstawione w \cite{BrezisNirenberg1978}. Warunki gwarantujące istnienie rozwiązania równania (\ref{eq:glowne}) są wtedy związane z zachowaniem się części nieliniowej $\N$ na $\ker\LL$ i noszą nazwę \emph{warunków typu Landesmana-Lazera} (\cite{LandesmanLazer1970},\cite{Hess1974}).\medskip

W przypadku, gdy rozważamy zagadnienia semi-liniowe na zbiorach nieograniczonych, to związane z nimi operatory liniowe mogą posiadać niepuste widmo istotne $\essp(\LL)$. Jest tak na przykład dla operatora Laplace'a $-\Delta$ rozważanego na całym $\bbR[n]$ (wtedy jego widmo jest w całości istotne i pokrywa się z nieujemną półosią rzeczywistą \cite[Ch. V.5.2, p. 299]{Kato1995}) lub w większości zagadnień związanych z operatorami Schr\"odingera $S=-\Delta+V$ (\cite{ReedSimon1978}). Praca ta ma na celu zbadanie rozwiązalności równania (\ref{eq:glowne}), gdy \emph{$0$ jest elementem $\essp(\LL)$ i punktem brzegowym luki w widmie $\LL$}. Podstawowe konsekwencje tego założenia, które odróżniają tę sytuację od przedstawianych w poprzednim paragrafie są następujące:
\begin{itemize}
\item po pierwsze jądro $\ker\LL$ może być trywialne, a co za tym idzie nie możemy w ogólności opierać się na badaniu części nieliniowej $\N$ na tej podprzestrzeni;
\item po drugie rezolwenta $\LL$ a także częściowe odwrócenia tego operatora tracą w tym przypadku zwartość, więc nie możemy stosować teorii stopnia koincydencji.
\end{itemize}
Chociaż $0$ nie musi już być punktem izolowanym widma operatora $\LL$, to jednak z powyższego założenia wynika, że istnieje prawo- lub lewostronne sąsiedztwo $0$ które jest w całości poza widmem $\LL$. Dzięki temu możemy dokonać rozkładu przestrzeni
\begin{equation}\tag{$\star\star$}\label{eq:rozklad}
\H=\lH\oplus\rH,
\end{equation}
zgodnego z podziałem widma na część poniżej i powyżej $0$. Wówczas zakładając dla ustalenia uwagi, że $0$ ma lewostronne sąsiedztwo w całości poza widmem, mamy następujące własności tego podziału. Część operatora $\LL$ działająca w $\lH$ ma widmo zawarte w ujemnej półosi rzeczywistej oraz, ze względu na istnienie luki, odseparowane od zera (czyli $\LL$ jest na tej podprzestrzeni odwracalny). Część $\LL$ działająca w $\rH$ ma widmo zawarte w nieujemnej półosi rzeczywistej. W szczególności jest ona operatorem nieujemnym.\medskip

Na poniższą dysertację składają się trzy główne części. W pierwszym rozdziale wprowadzamy potrzebne pojęcia oraz przedstawiamy podstawowe rezultaty wykorzystywane dalej. Na początku przywołujemy definicje i fakty dotyczące operatorów liniowych (w ogólności nieograniczonych) w przestrzeniach Hilberta i ich widm. Ta część kończy się podaniem twierdzenia spektralnego dla operatorów samosprzężonych (twierdzenie \ref{thm:spectral}), które jest jednym z podstawowych narzędzi wykorzystywanych w tej pracy. Stwierdza ono, że dla każdego operatora liniowego i samosprzężonego $\LL$ istnieje jednoznacznie wyznaczona, prawostronnie ciągła, jednoparametrowa rodzina samosprzężonych projekcji (ortoprojekcji) w $\H$ $\{E_\mu:\ \mu\in\bbR{}\}$, dla której
\[
\LL=\int_{-\infty}^{\infty}\mu\,dE_\mu,
\]
gdzie całka po prawej stronie (tzw. całka spektralna) jest rozumiana jako granica w $\H$ odpowiednich sum całkowych. Ponadto możemy wówczas dla dowolnej funkcji ciągłej $f\colon\bbR\to\bbC$ zdefiniować
\[
f(\LL)=\int_{-\infty}^{\infty}f(\mu)\,dE_\mu.
\] 
W ten sposób możemy określić tzw. rachunek funkcyjny dla operatorów samosprzężonych w przestrzeni Hilberta $\H$. W dalszej części pracy, przy użyciu wprowadzonych wcześniej pojęć, formułujemy precyzyjnie problem badany w tej dysertacji oraz opisujemy motywacje za nim stojące. Następnie, przy wykorzystaniu twierdzenia spektralnego, definiujemy podprzestrzenie $\lH$ oraz $\rH$ i dowodzimy prawdziwości rozkładu (\ref{eq:rozklad}). W sekcji \ref{sec:monot_oper} zajmujemy się operatorami nieliniowymi monotonicznymi i maksymalnie monotonicznymi. Podajemy podstawowe definicje oraz fakty dotyczące tych odwzorowań, które wykorzystywane są potem w rozdziale \ref{ch:trivial_ker}. Ponadto przywołujemy definicje operatorów quasi-monotonicznych oraz klasy $(S_+)$. Są one przykładami tak zwanych odwzorowań typu monotonicznego, wykorzystywanych do badania rozwiązalności równań nieliniowych. W sekcji \ref{sec:top_deg} wprowadzamy na podstawie pracy Kartsatosa i Skrypnika \cite{KartsatosSkrypnik1999} pewne rozszerzenie klasy $(S_+)$ obejmujące operatory, które są tylko gęsto określone w przestrzeni Hilberta. Na koniec rozdziału \ref{ch:prelim}, w sekcji \ref{sec:top_deg} podajemy twierdzenie o stopniu topologicznym dla tej klasy odwzorowań które stosujemy w rozdziale \ref{ch:inf_essential}.\medskip

W następnych dwóch rozdziałach koncentrujemy się na dowodzeniu twierdzeń o rozwiązalności równania (\ref{eq:glowne}). Odbywa się to przy udziale dodatkowych założeń narzuconych na część nieliniową $\N$. Podstawową metodą stosowaną w tej pracy jest metoda perturbacyjna. Najpierw pokazujemy że istnieją rozwiązania zaburzeń wyjściowego równania postaci
\begin{equation}\tag{$\star\star\star$}\label{eq:zaburzone}
\eps P_2u+\LL u+\N(u)=h,
\end{equation}
gdzie $P_2$ jest ortoprojekcją na podprzestrzeń $\rH$ a $\eps>0$ jest parametrem perturbacyjnym. Rozwiązania równania (\ref{eq:glowne}) szukamy następnie wśród granic
\[
\lim_{\eps\to0}u_\eps,
\]
gdzie $u_\eps\in\H$ dla $\eps>0$ są rozwiązaniami równań zaburzonych (\ref{eq:zaburzone}). Ten proces odbywa się w dwóch krokach. W pierwszym wykazujemy, że przy jednostajnej ograniczoności w normie rozwiązań $u_\eps$ względem ograniczonego parametru $0<\eps<C$ powyższe granice, rozumiane w sensie słabej zbieżności, generują rozwiązania (\ref{eq:glowne}). W drugim kroku podajemy dodatkowe warunki, które zapewniają tę ograniczoność. Ponadto wykorzystujemy stopień topologiczny dla operatorów typu monotonicznego oraz własności odwzorowań maksymalnie monotonicznych w przestrzeniach Hilberta. Podkreślmy także, że w przypadku gdy równanie (\ref{eq:glowne}) posiada rozwiązanie trywialne (tzn. $N(0)=h$), to metody stosowane w tej pracy nie pozwalają wykazać istnienia rozwiązań nietrywialnych. Wynika to z faktu, że nie mamy wystarczającej kontroli nad ciągami rozwiązań równań zaburzonych (\ref{eq:zaburzone}). W tym przypadku udowodnione tutaj twierdzenia o rozwiązalności się trywializują.  %oraz teorię punktów stałych w przestrzeniach lokalnie wypukłych.
\medskip

W drugim rozdziale zakładamy dodatkowo że $0$ jest kresem dolnym widma istotnego. Wówczas okazuje się, że operator $(\LL+\eps I)P_2$ działający w $\rH$ przynależy do klasy odwzorowań gęsto określonych spełniających warunek $(S_+)$, która jest zdefiniowana w sekcji \ref{sec:top_deg} rozdziału \ref{ch:prelim}. Jeżeli spojrzymy na lewą stronę równania (\ref{eq:zaburzone}) jak na nieliniowe zaburzenie $\LL+\eps P_2$, to właściwym założeniem o operatorze $\N$, które pozwala pozostać w obrębie wspomnianej klasy, jest jego quasi-monotoniczność, tzn.
\[
\limsup_{k\to\infty}\langle\N(u_k), u_k-u_0\rangle_r\geqslant 0,
\]
o ile $u_k\weak u_0$. Tutaj $\weak$ oznacza słabą zbieżność w $\H$ a $\langle\cdot\,,\cdot\rangle_r:= \re\langle\cdot\,,\cdot\rangle$ jest iloczynem skalarnym w urzeczywistnieniu przestrzeni Hilberta $\H$. Operatory zwarte lub monotoniczne są też quasi-monotoniczne, jednak warunek ten jest istotnie od nich słabszy. Stosując odpowiedni stopień topologiczny dowodzimy w tym rozdziale twierdzenia kontynuacyjnego dla równań zaburzonych (theorem \ref{thm:continuation}) i przy jego pomocy oraz założenia o subliniowym wzroście $\N$, tzn.
\[
\lim_{k\to\infty}\frac{\|\N(u_k)\|}{\|u_k\|}=0\quad \text{o ile}\ \|u_k\|\to\infty,
\]
wykazujemy istnienie rozwiązań równań (\ref{eq:zaburzone}). Rozwiązalność (\ref{eq:glowne}) otrzymujemy badając tzw. funkcjonał recesji $J_n\colon\H\to[-\infty,+\infty]$, który jest określony formułą
\[
J_{\N}(u)=\inf\left\{\liminf_{k\to+\infty}\langle N(t_kv_k),v_k\rangle_r:\ t_k\to+\infty,\ \{v_k\}_{k\in\bbN}\subset\H,\ v_k\weak u\right\}.
\]
Funkcja $J_N$ jest związana z zachowaniem się w nieskończoności funkcjonału postaci $u\mapsto\langle\N(u),u\rangle_r$. Głównym wynikiem tej części pracy jest następujące
\begin{twierdzenie} %[\ref{thm:main_ex}]
Niech $\LL$ będzie operatorem liniowym samosprzężonym i ograniczonym z dołu w zespolonej i ośrodkowej przestrzeni Hilberta $\H$. Załóżmy ponadto, że $0$ jest kresem dolnym widma istotnego i $\LL$ ma tylko skończenie wiele ujemnych wartości własnych. Niech ponadto $h\in\H$. Jeżeli $\N\colon\H\to\H$ jest odwzorowaniem ograniczonym, demiciągłym, quasi-monotonicznym oraz zachodzi
\begin{enumerate}[(i)]
\item $\displaystyle\lim_{k\to+\infty}\frac{\|\N(u_k)\|^2}{\|u_k\|}=0$ dla każdego ciągu $\{u_k\}_{k\in\bbN}\subset\H$ takiego że $\|u_k\|\to+\infty$,
\item $\displaystyle \limsup_{k\to\infty}\frac{\langle\N(u_k),u_k\rangle_r}{\|u_k\|}>0$ dla wszystkich $\{u_k\}_{k\in\bbN}\subset\H$ takich że $\|u_k\|\to+\infty$.
\item $J_{\N}(u)>\langle h,u\rangle_r$ dla każdego $u\in\ker\LL$, $\|u\|=1$,
\end{enumerate}
to równanie (\ref{eq:glowne}) posiada rozwiązanie.
\end{twierdzenie}
Punkt \textit{(ii)} w powyższym twierdzeniu jest warunkiem znakowym niezbędnym przede wszystkim w sytuacji gdy $\ker\LL$ jest trywialne, czego nie można wykluczyć gdy wiemy tylko że $0\in\essp(\LL)$. Na końcu tego rozdziału podajemy przykład zastosowania powyższego twierdzenia w przypadku gdy
\[
\N(u)=\phi(\|u\|)P_B(u),\quad u\in\H,
\]
gdzie $\phi\colon[0,+\infty)\to\bbC$ a $P_B$ jest projekcją metryczną na domkniętą kulę jednostkową w $\H$.\medskip

W rozdziale trzecim opuszczamy założenie, że $0$ jest kresem dolnym widma istotnego. W tym przypadku, w przeciwieństwie do poprzedniego rozdziału, podprzestrzeń $\lH$ jest nieskończenie wymiarowa. Podstawowym spostrzeżeniem w tej części pracy jest fakt, że operator $\LL P_2$ działający z $\H$ do $\rH$ jest operatorem maksymalnie monotonicznym. Aby odwzorowanie $\LL P_2+\N$ pozostało w tej klasie zakładamy, że $\N$ spełnia pewien warunek monotoniczności. Jest to założenie silniejsze niż quasi-monotoniczność z drugiego rozdziału. Następnie dokonujemy wielowartościowego odwrócenia równań zaburzonych (\ref{eq:zaburzone}) i pokazujemy, że są one równoważne inkluzji
\[
\LL^{-1}_\eps h\in A(v)+\LL^{-1}_\eps v,\quad v\in\dom(A)
\]
gdzie $\LL_{\eps}=\LL P_1+\eps P_2$ oraz $A=(\LL P_2+\N)^{-1}$. Tutaj $P_1$ jest ortoprojekcją na $\lH$. Wykorzystując własności operatorów maksymalnie monotonicznych (w szczególności twierdzenia dotyczące ich suriektywności) wykazujemy rozwiązalność powyższej inkluzji. Istnienie rozwiązań równania (\ref{eq:glowne}) dowodzimy przy wykorzystaniu funkcji recesji oraz dodatkowego warunku znakowego dla operatora nieliniowego $\N$. Głównym wynikiem tego rozdziału jest
\begin{twierdzenie}
Niech $\LL$ będzie operatorem liniowym samosprzężonym i ograniczonym z dołu przez $-\ga <0$ w zespolonej i ośrodkowej przestrzeni Hilberta $\H$. Załóżmy ponadto, że $0$ należy do widma istotnego $\LL$ oraz istnieje $\del>0$ taka że przedział $(-\del,0)$ jest w całości poza widmem. Niech ponadto $h\in\H$. Jeżeli $\N\colon\H\to\H$ jest odwzorowaniem ograniczonym, demiciągłym, $\N(0)=0$ oraz
\begin{enumerate}[(i)]
\item istnieje $\al>\ga/\del^2$ taka że dla dowolnych $u,u'\in\H$
\[
 \langle\N(u)-N(u'),u-u'\rangle_r\geqslant\al\|\N(u)-\N(u')\|^2,
\]
\item $\displaystyle \limsup_{k\to\infty}\frac{\langle\N(u_k),u_k\rangle_r}{\|u_k\|}>\frac{\ga\|h\|}{\del^2\al-\ga}$ dla wszystkich $\{u_k\}_{k\in\bbN}\subset\H$ takich że $\|u_k\|\to+\infty$.
\item $J_{\N}(u)>\langle h,u\rangle_r$ dla każdego $u\in\ker\LL$, $\|u\|=1$,
\end{enumerate}
to równanie (\ref{eq:glowne}) posiada rozwiązanie.
\end{twierdzenie}
Zauważmy, że punkt \textit{(i)} implikuje w szczególności że operator $\N$ jest monotoniczny, tzn. spełnia
\[
\langle\N(u)-N(u'),u-u'\rangle_r\geqslant 0
\]
dla dowolnych $u,u'\in\H$. Ponadto jeżeli przyjmiemy $u'=0$ to otrzymujemy
\[
\|\N(u)\|^2\leqslant\frac{1}{\al}\langle\N(u),u\rangle_r,
\]
dla każdego $u\in\H$. Zatem można go też interpretować jako pewien warunek wzrostu dla $\N$. W szczególności pokazujemy, że w przypadku gdy $\N$ jest operatorem Niemyckiego w $\LR[n]$, tzn.
\[
\N(u)(x)=f(x,u(x)),\ x\in\bbR[n],
\]
dla dowolnego $u\in\LR[n]$, gdzie $f\colon\bbR[n]\times\bbR\to\bbR$, to $\N$ spełnia \textit{(i)} jeżeli $f$ jest ze względu na drugą współrzędną monotoniczna oraz jej wzrost jest jednostajnie ograniczony. Jako przykład zastosowania powyższego twierdzenia wykazujemy rozwiązalność równania Schr\"{o}dingera
\[
-\Del u+V(x)u+f(x,u)=h(x)
\]
w $\LR[3]$.

%W czwartym rozdziale opuszczamy założenie o ograniczoności $\LL$ z dołu. Jest to najbardziej ogólny zestaw założeń o $\LL$, który rozważamy i zarazem wydaje się najtrudniejszy z punktu widzenia zastosowań. Ze względu na wspomniany wcześniej brak zwartości zakładamy w tej części, że operator $\N$ jest ciągowo ciągły w słabych topologiach na $\H$. Wyklucza to możliwość rozpatrywania jako $\N$ operatora superpozycji (za wyjątkiem przestrzeni z miarą dyskretną, czyli tak naprawdę równań różnicowych) i kieruje nas w stronę równań z nielokalną częścią nieliniową. Przy wykorzystaniu twierdzenia o punkcie stałym w $\H$ ze słabą topologią oraz subliniowego wzrostu $\N$ wykazujemy rozwiązalność równań zaburzonych. Istnienie rozwiązań (\ref{eq:glowne}) otrzymujemy przy dodatkowym założeniu, będącym warunkiem znakowym między wartościami $\N$ na elementach z zewnętrza kuli o dostatecznie dużym promieniu, a rzutami tych elementów na $\rH$.

\cleardoublepage

\tableofcontents

\chapter{Introduction}

Keywords: semilinear equations, essential spectrum of linear operator, topological degree, maximal monotone operators, monotony type operators
\bigskip

The main goal of this dissertation is to find conditions which will guarantee the existence of solutions in the Hilbert space $\H$ of semilinear equation
\begin{equation}\label{eq:main_intro}\tag{$\star$}
\LL u+\N(u)=h
\end{equation} 
where $\LL$ is a linear and self-adjoint operator, $\N$ a non-linear mapping and $h\in\H$. In this project we concentrate on the case when $0$ belongs to the essential spectrum of operator $\LL$ which was not previously studied in this general setting.\medskip

The essential spectrum emerges in particular when $\LL$ comes from a differential operator defined on an unbounded domain in $\bbR[n]$. For example this is the case for the Laplace operator $-\Del$ on $\bbR[n]$ \cite[Ch. V.5.2, p. 299]{Kato1995} or for a wide class of Schr\"{o}dinger operators $S=-\Del+V$ \cite{ReedSimon1978}. The goal of this work is to study the solvability of equation (\ref{eq:main_intro}) when \emph{$0$ is an element of essential spectrum of operator $\LL$ and the boundary point of its spectral gap}. The basic feature of this situation, which distinguish it from the case of operators with purely discrete spectrum, is that the kernel of $\LL$ can be trivial and the resolvent of $\LL$ is not a compact operator. To compensate for this lack of compactness we explore various monotonicity conditions. We also rely on the perturbation technique in which we firstly consider equations with parameter $\eps$ such that for $\eps=0$ we get (\ref{eq:main_intro}). We prove the solvability of these supplementary equations and then search for the solutions of (\ref{eq:main_intro}) as limits of sequences of perturbed solutions as $\eps$ goes to zero.
\medskip

This thesis is divided into three parts.\medskip

Chapter \ref{ch:prelim} is concerned with prerequisites. Firstly, we recall the basic definitions and facts about linear (generally unbounded) operators in a Hilbert space and their spectra. This culminates with the formulation of the spectral theorem for self-adjoint operators (theorem \ref{thm:spectral}) which is one of our basic tools. Then we precisely describe the problem of this thesis, give the motivation for its study and make the decomposition of a Hilbert space $\H$ according to spectral properties of $\LL$. Next in section \ref{sec:monot_oper} we move to non-linear mappings and give a brief survey of the theory of maximal monotone operators which is used in chapter \ref{ch:trivial_ker}. Moreover we formulate definitions and give some remarks concerning quasi-monotone and monotone type operators of class $(S_+)$ applied in chapter \ref{ch:inf_essential}. Finally in section \ref{sec:top_deg} we present the topological degree for mappings of class $(S_+)$.\medskip

In chapter \ref{ch:inf_essential} we additionally assume that $0$ is the infimum of the essential spectrum of $\LL$. Hence, we can only have a discrete set of eigenvalues of finite multiplicity below zero. We apply the degree theory for mappings of class $(S_+)$ to the operator given by the left hand side of equation (\ref{eq:main_intro}). We assume that non-linear part $\N$ is quasi-monotone and satisfies sublinear growth condition. The former assumption is crucial to ensure that the perturbation of the left hand side of (\ref{eq:main_intro}) is of class $(S_+)$ on certain subspace of $\H$ and the latter allows us to control the perturbed solutions. Moreover, since $0$ can have non-trivial eigenspace, we make use of the so called \emph{recession functional} connected with $\LL$ and $\N$. It allows us to control the behaviour of non-linear part on the kernel of $\LL$. \medskip

In chapter \ref{ch:trivial_ker} we allow the essential spectrum of $\LL$ to lay below zero. Our method is based on the observation that certain perturbation of operator $\LL$ is maximal monotone on the subspace of $\H$ corresponding to non-negative part of the spectrum of $\LL$. In order to extend this property to the perturbation of the left hand side of (\ref{eq:main_intro}) we assume that operator $\N$ satisfies as well certain monotonicity condition. Next we explore the surjectivity properties of maximal monotone operators to show the existence of solutions to the class of perturbed equations. Finally with the help of recession functional and a growth conditions on $\N$ we make the limiting step and prove the solvability of (\ref{eq:main_intro}).

\mainmatter

\chapter{Preliminaries}\label{ch:prelim}

\section{Formulation of the problem and decomposition of space}\label{ssec:problem_descr}

Let $\H$ be a real or complex Hilbert space. The symbol $\langle\cdot\,,\cdot\rangle$ will indicate the scalar product on $\H$ (in complex case we assume the linearity on the second coordinate) and $\|\cdot\|$ will denote the norm given by this product. Symbol $\bbN$ stands for the set of natural numbers $\{1,2,3,\ldots\}$ and $\bbNo=\bbN\cup\{0\}$. If $\{u_k\}_{k\in\bbNo}\subset\H$ then the formula $u_k\to u_0$ will represent convergence in the norm of the space $\H$, that is
\[
\lim_{k\to\infty}\|u_k-u_0\|=0,
\]
and by $u_k\weak u_0$ we will mean the convergence of this sequence in the weak topology $\si(H,H^*)$ of $\H$, equivalent to the following condition:
\[
\lim_{k\to\infty}\langle u_k-u_0,v\rangle=0
\]
for every $v\in\H$.

A subset $S\subset\H$ will be called a \emph{subspace} of space $\H$ if it is a linear space which is closed in the norm of $\H$. In case of lack of closeness we will call $S$ a \emph{(linear) submanifold} of space $\H$. If $r>0$ and $u\in\H$ the ball of radius $r$ centred at $u$ will be denoted $B(u,r)$ and corresponding sphere as $S(u,r):=\partial B(u,r)$.\medskip

\subsection{Linear operators in Hilbert space}\label{sssec:lin_oper}

Let $\lin(\H)$ denote the family of all linear operators $L\colon \dom(L)\subset\H\to\H$ where $\dom(L)$ is a dense linear submanifold of space $\H$. The \emph{adjoint operator} $\LL^*\in\lin(H)$ of $\LL\in\lin(\H)$ is given by the formula
\begin{equation*}
\left\{\begin{array}{l}\medskip
\dom(\LL^*)=\{v\in\H:\ (\exists w\in\H)\ (\forall u\in \dom(\LL))\ \langle v,\LL u\rangle=\langle w,u\rangle\},\\ 

L^*v=w,\ \text{for all}\ v\in \dom(\LL^*).
\end{array}\right.
\end{equation*}
We say that operator $\LL\in\lin(\H)$ is \emph{symmetric} if $\LL\subset\LL^*$ (that is $\dom(\LL)\subset \dom(\LL^*)$ and $\LL u=\LL^* u$ for all $u\in \dom(\LL)$), \emph{normal} if $\LL\LL^*=\LL^*\LL$ and \emph{self-adjoint} if $\LL=\LL^*$. Let us mark out a subfamily of symmetric operators. A symmetric operator $\LL\in\lin(\H)$ is called \emph{bounded from below} if there exists a number $\ga\in\bbR$ such that
\[
\langle u,\LL u\rangle\geqslant\ga\langle u,u\rangle
\]
for each $u\in\dom(\LL)$. When it is true we will simply write $\LL\geqslant\ga$ and if $\LL\geqslant 0$ we say that $\LL$ is \emph{non-negative}. Moreover a symmetric operator $\LL\in\lin(\H)$ is \emph{bounded from above} if $-\LL$ is bounded from below.

We say that operator $\LL\in\lin(\H)$ is \emph{closed} if its graph 
\[
\graph(\LL):=\{(u,\LL u):\ u\in \dom(\LL)\}
\]
is a closed set (subspace) of direct sum $\H\oplus\H$ considered with scalar product
\[
\langle u,v\rangle +\langle u',v'\rangle
\]
for all $(u,u'),(v,v')\in\H\times\H$. The family of all closed operators in Hilbert space $\H$ is denoted by $\clos(\H)$.

An operator $\LL\in\lin(\H)$ is \emph{bounded} if
\[
\sup\{\|\LL u\|:\ u\in \dom(\LL),\ \|u\|\leqslant 1\}<\infty,
\]
and $\bound(\H)$ will stand for the family of all bounded operators in $\H$. If $\LL\in\bound(\H)$ then, using the density of $\dom(\LL)$ in $\H$, we can uniquely extend it to the whole space $\H$ without increasing this supremum \cite[Thm. I.7, p. 9]{ReedSimon1980}. According to this fact we will henceforward always assume that if $\LL\in\bound(\H)$ then $\dom(\LL)=\H$ and we will call the \emph{norm} of $\LL$ the number
\[
\|\LL\|:=\sup\{\|\LL u\|:\ u\in\H,\ \|u\|\leqslant 1\}.
\]
Let $\{\LL_k\}_{k\in\bbNo}\subset\bound(\H)$. Since $\bound(\H)$ with the above norm becomes a Banach space we can consider the convergence in this norm which will be denoted by 
\[
\lim_{k\to\infty}\LL_k=\LL_0.
\]
We will also use \emph{strong (pointwise) convergence} defined by condition 
\[
\lim_{k\to\infty}\|(\LL_k-\LL_0)u\|=0,
\]
for all $u\in\H$, and indicated as
\[
\slim{k\to\infty}\LL_k=\LL_0.
\]
In the space of bounded and symmetric operators we can introduce partial ordering $"\leqslant"$. To this end let $\LL_1,\LL_2\in\bound(\H)$ be symmetric. Then we write $\LL_1\leqslant\LL_2$ if for all $u\in\H$ we have
\[
\langle u,\LL_1u\rangle\leqslant\langle u,\LL_2u\rangle.
\]
If $\LL_1=\ga I$, where $I$ is an identity on $\H$, the foregoing definition coincides with the notion of boundedness from below defined earlier.

Finally let us define the family of \emph{orthoprojections} in $\H$, which we will denote $\proj(\H)$. It consists of all self-adjoint operators $\LL\in\bound(\H)$ satisfying the idempotent condition $\LL\circ\LL=\LL$. For all $\LL\in\proj(\H)$ we have 
\[
\ran(\LL)\perp\ker\LL,\ \H=\ran(\LL)\oplus\ker\LL,\ \text{and}\  0\leqslant\LL\leqslant I.
\]

To conclude this review let us write some rather trivial connections between various classes defined above.
\begin{enumerate}
\item $\proj(\H)\subset\bound(\H)\subset\clos(\H)\subset\lin(\H)$
\item If $\LL\in\lin(\H)$ is symmetric then $\LL\in\clos(\H)$.
\item If $\LL\in\clos(\H)$ and $\dom(\LL)=\H$ then $\LL\in\bound(\H)$ (closed graph thm.)
\end{enumerate}

\subsection{Spectrum and the spectral theorem}\label{sssec:spectrum}

Let $\H$ be a complex Hilbert space. If $\LL\in\clos(\H)$ we define the \emph{resolvent set} $\rez(\LL)\subset\bbC$, where $\bbC$ stands for the field of complex numbers, of an operator $\LL$ in the following way
\[
\rez(\LL):=\{\lam\in\bbC:\ \LL-\lam I\ \text{is injective and}\ (\LL-\lam I)^{-1}\in\bound(\H)\}.
\]
\medskip

\begin{remark}\label{rem:resolvent}

\begin{enumerate}[(a)]
\item Throughout this paper we will use a convention in which the symbol of identity operator $I$ will be omitted in formulas like the one above, i.e., we use the notation
\[
\LL-\lam:=\LL-\lam I
\]
for any number $\lam\in\bbC{}$.
\item If $\LL\in\clos(\H)$, according to closed graph theorem, we have the following characterisation of the resolvent set
\[
\rez(\LL)=\{\lam\in\bbC:\ \LL-\lam\ \text{is bijection from}\ \dom(\LL)\ \text{onto}\ \H\}.
\]
\item The assumption that $\LL\in\clos(\H)$ is important. The operator $\LL\in\lin(\H)$ is closed iff for every $\lam\in\bbC$ the operator $\LL-\lam$ is closed iff $(\LL-\lam)^{-1}$ is closed (if it exists). Hence if $\LL\in\lin(\H)$ is not closed then even if $\LL-\lam$ is bijection from $\dom(\LL)$ onto $\H$, for some $\lam\in\bbC$, operator $(\LL-\lam)^{-1}$ is not closed and in particular cannot belong to $\bound(\H)$. So in this case we always have $\rez(\LL)=\emptyset$.
\end{enumerate}
\end{remark}\medskip

Assuming that $\lam\in\rho(\LL)$ we can consider the operator 
\[
\R{\LL}(\lam):=(\LL-\lam)^{-1}
\]
which is called \emph{the resolvent} of $\LL$ at point $\lam$. Finally, \emph{the spectrum} of $\LL\in\clos(\H)$ is the subset of the complex plane defined as
\[
\sp(\LL):=\bbC{}\bsh\rez(\LL).
\]
Remark \ref{rem:resolvent}$(b)$ allows us to further divide $\sp(\LL)$ into disjoint parts as follows
\[
\sp(\LL)=\psp(\LL)\cup\csp(\LL)\cup\rsp(\LL),
\]
where
\begin{align*}
&\psp(\LL)=\{\lam\in\sp(\LL):\ \ker(\LL-\lam)\neq\{0\}\},\\
&\csp(\LL)=\{\lam\in\sp(\LL)\bsh\psp(\LL):\ \cl{\ran(\LL-\lam)}=\H \},\\
&\rsp(\LL)=\{\lam\in\sp(\LL)\bsh\psp(\LL):\ \cl{\ran(\LL-\lam)}\neq\H \},
\end{align*}
which we call respectively \emph{point, continuous and residual spectrum}. Every $\lam\in\psp(\LL)$ is called an \emph{eigenvalue} of an operator $\LL$, the set $\ker(\LL-\lam)$ is the \emph{eigenspace} associated with $\lam$ and its dimension $\dim\ker(\LL-\lam)\in\bbNo\cup\{\infty\}$ is the (geometrical) \emph{multiplicity} of $\lam$.\medskip

\begin{remark}\label{rem:spectrum}
\begin{enumerate}[(a)]
\item The spectrum $\sp(\LL)$ of $\LL\in\clos(\H)$ is always closed. If $\LL$ is self-adjoint and $\H$ is not trivial then it is also nonempty \cite[Prop. 5.4.1, p. 172]{BlaExnHav2008}.
\item The spectrum of self-adjoint operator $\LL\in\lin(\H)$ is real. The converse holds only if we know \textit{a priori} that $\LL$ is symmetric \cite[Ch. V.4, V.5]{Kato1995}. Moreover the residual spectrum of a self-adjoint operator is always empty (since $\ker\LL=\ran\LL^*$, see \cite[Prop. 4.1.1(b), p. 93]{BlaExnHav2008}).
\item For a self-adjoint operator $\LL\in\lin(\H)$ the assumption about its boundedness from below is equivalent to the fact that its spectrum $\sp(\LL)$ (contained in $\bbR$, see item $(b)$ above) is bounded form below. In this case the infimum of $\sp(\LL)$ is equal to supremum of lower bounds of operator $\LL$ \cite[Ch. V.10]{Kato1995}. The analogous remark applies to boundedness from above and in particular we see that a self-adjoint operator belongs to $\bound(\H)$ iff its spectrum is a bounded subset of $\bbR$. 
\end{enumerate}
\end{remark}\medskip

The salient result of the theory of self-adjoint operators in Hilbert space is the following
\begin{theorem}[spectral theorem for self-adjoint operators]\label{thm:spectral}
Let $\H$ be a complex and separable Hilbert space and let $\LL\in\lin(\H)$ be self-adjoint. There exists one-parameter family $\{E_\mu:\ \mu\in\bbR{}\}\subset\proj(\H)$ such that
\begin{enumerate}[(i)]
\item\label{thm:spec_monot} $E_\mu E_\nu=E_{\min\{\mu,\nu\}}$ for every $\mu,\nu\in\bbR$,
\item\label{thm:spec_right_cont} $E_{\mu+}:=\slim{k\to\infty}E_{\mu+1/k}=E_\mu$ for all $\mu\in\bbR{}$,
\item $E_{-\infty}:=\slim{k\to\infty}E_{-k}=0$, $E_\infty:=\slim{k\to\infty}E_k=1$.
\end{enumerate}
Moreover for every continuous function $f\colon\bbR\to\bbC$ the formulas 
\begin{align*}
& \dom(f(\LL))=\left\{u\in\H:\ \int_{-\infty}^{\infty}|f(\mu)|^2d\langle u,E_\mu u\rangle<\infty\right\},\\
& \langle v,f(\LL)u\rangle=\int_{-\infty}^{\infty}f(\mu)\ d\langle v,E_\mu u\rangle,\ v\in\H,\ u\in\dom(f(\LL)),
\end{align*}
where integrals are understood in the sense of Riemann-Stieltjes, define the normal operator $f(\LL)\in\lin(\H)$, which we denote
\[
f(\LL)=\int_{-\infty}^{\infty}f(\mu)dE_\mu,
\]
such that
\[
\|f(\LL)\|^2=\int_{-\infty}^{\infty}|f(\mu)|^2d\langle u,E_\mu u\rangle.
\]
\end{theorem}
For the proof see \cite[Ch. VI.2 and VI.3]{Kato1995}. For more general version with borel functions we refer to \cite[Ch. VIII.3]{ReedSimon1980} or \cite[Ch. 5.3 and 5.5]{BlaExnHav2008}.\medskip

\begin{remark}\label{rem:spectral}
\begin{enumerate}[(a)]
\item Condition $(i)$ from the theorem can be formulated equivalently -- for every $\mu,\nu\in\bbR$
\[
E_\mu\leqslant E_\nu,
\]
when $\mu\leqslant\nu$ (which is also equivalent to $\ran(E_\mu)\subset\ran(E_\nu)$). Particularly for each $u\in\H$ the function
\[
\bbR\ni\mu\mapsto\langle u,E_\mu u\rangle
\]
is non-decreasing.
\item Condition $(ii)$ constitutes continuity from the right of the family $\{E_\mu:\ \mu\in\bbR\}$. In addition we will use the following notation
\[
E_{\mu-}:=\slim{k\to\infty}E_{\mu-1\bsh k}.
\]
\item The family $\{E_\mu:\ \mu\in\bbR{}\}$ is called the \emph{spectral family} (or \emph{resolution of identity}) of an operator $\LL$. It is "concentrated" on the spectrum of $\LL$ in the following sense \cite[Prop. 5.4.1(a), p. 171]{BlaExnHav2008}: if $-\infty<a<b<\infty$ we have
\[
\sp(\LL)\cap(a,b)=\emptyset\ \text{iff}\ E_{a}=E_{b-}.
\]
Hence for all $u,v\in\H$ the function
\[
\mu\mapsto\langle v,E_\mu u\rangle
\]
is constant on components of $\rez(\LL)$. In particular it means that in the spectral theorem \ref{thm:spectral} it is enough to define the function $f$ on $\sp(\LL)$ (Tietze theorem).
\item In theorem \ref{thm:spectral} the \emph{continuous functional calculus} for self-adjoint operator is implicitly defined. More explicitly for every self-adjoint $\LL\in\lin(\H)$ we have an injective map of algebra $C(\bbR,\bbC)$  into $\lin(\H)$
\[
C(\bbR,\bbC)\ni f\mapsto f(\LL),
\]
given by the spectral integral. If $f$ is real-valued, $f(\LL)$ is also self-adjoint and if $f$ is bounded, $f(\LL)\in\bound(\H)$ with $\|f(\LL)\|=\sup_{\bbR}|f|$. For all properties of this calculus see \cite[Ch. VI.5.2]{Kato1995} and in a more general borel-setting we refer to \cite[Ch. 5.2]{BlaExnHav2008}. Let us also stress out that in case $f=id_{\bbR}$ we have the spectral representation of an operator $\LL$ itself in the form
\[
\LL=\int_{-\infty}^{\infty}\mu dE_\mu.
\]
\item\label{thm:spectr_mapping}
\emph{Spectral mapping theorem} characterizes the spectrum of an operator $f(\LL)$, where $f\in C(\sp(\LL))$, by identity \cite[Prop. 5.5.3, p. 178]{BlaExnHav2008}
\[
\sp(f(\LL))=\cl{f(\sp(\LL))}.
\]
When $\LL\in\bound(\H)$ we have
\[
\sp(f(\LL))=f(\sp(\LL)).
\]
\item Finally let us describe the point and the continuous spectrum of a self-adjoint operator $\LL\in\lin(\LL)$ using its spectral family \cite[Rem. 5.4.2(b), p. 172]{BlaExnHav2008}.
\[
\psp(\LL)=\{\lam\in\bbR{}:\ E_{\lam-}\neq E_{\lam+}\}
\]
This means that eigenvalues of $\LL$ are precisely the points of discontinuity of its resolution of identity. Moreover we have $\ran(E_{\lam+}-E_{\lam-})=\ker(\LL-\lam)$.
\[
\csp(\LL)=\{\lam\in\bbR{}:\ E_{\lam-}=E_\lam,\ (\forall\eps>0)\ E_{\lam-\eps}\neq E_{\lam+\eps}\}
\]
Hence the elements of the continuous spectrum are points of continuity of the spectral family such, that the map $\mu\mapsto E_\mu$ is non-constant in their vicinity.
\end{enumerate}
\end{remark}

Now we will introduce another decomposition of the spectrum of a self-adjoint operator which will be particularly important in the formulation of the problem investigated in this work.

\begin{definition}\label{def:ess_spec}
Let $\LL\in\lin(\H)$ be self-adjoint and let $\{E_\mu:\ \mu\in\bbR{}\}$ be its spectral family. The \emph{essential spectrum} of $\LL$ is defined as
\[
\essp(\LL):=\{\lam\in\bbR{}:\ (\forall\eps>0)\ \dim\ran(E_{\lam+\eps}-E_{\lam-\eps})=\infty\}.
\]
The \emph{discrete spectrum} of $\LL$ is the completion in spectrum of $\essp(\LL)$, i.e.,
\[
\dsp(\LL):=\sp(\LL)\bsh\essp(\LL).
\]
\end{definition}
We already know from remark \ref{rem:spectral}$(c)$ that a number $\lam\in\bbR{}$ belongs to the spectrum of a self-adjoint operator $\LL$ iff $\dim\ran(E_{\lam+\eps}-E_{\lam-\eps})>0$ for each $\eps>0$. In consequence we can say that essential spectrum collects this elements $\lam\in\sp(\LL)$ where the "essential" jump in dimension of $\ran(E_\mu)$ occurs at the threshold $\mu=\lam$. In particular if $\lam$ is an eigenvalue with infinite multiplicity then it is always a member of $\essp(\LL)$. To be more specific we give below full (topological) characterisation \cite[Thm. 5.4.4., p. 173]{BlaExnHav2008}.

\begin{lemma}\label{lem:ess_char}
Let $\LL\in\lin(\LL)$ be self-adjoint. Then
\begin{align*}
\essp(\LL)&=\{\lam\in\bbR:\ \lam\ \text{is an accumulation point of}\ \sp(\LL)\\
          &\hspace{2.1cm}\text{or an eigenvalue of infinite multiplicity}\}\\
\dsp(\LL) &=\{\lam\in\sp(\LL):\ \lam\text{ is an isolated point of}\ \sp(\LL)\ \text{and}\  \dim\ran(E_{\lam}-E_{\lam-})<\infty\}.
\end{align*}
Particularly $\essp(\LL)$ is closed and $\dsp(\LL)\subset\psp(\LL)$.
\end{lemma}
According to this lemma we can describe the discrete spectrum of $\LL$ as a collection of its isolated eigenvalues with finite multiplicity.\medskip

\begin{remark}\label{rem:fredholm}
If $\X$ is a complex Banach space then, in an analogous fashion, we can define the families of linear, closed and bounded operators in $\X$. Definitions of the spectrum of a closed operator and its pointwise, continuous and residual parts are also similarly given. For a comprehensive treatment see \cite[Ch. 1.7]{BlaExnHav2008} or \cite[Ch. III.5]{Kato1995}. For $\LL\in\lin(\X)$ we can also use the natural pairing between $\X$ and $\X^*$ to define the adjoint operator $\LL^*\in\lin(\X^*)$. However, in this framework the notion of self-adjointness is of course meaningless.

Nevertheless, the notion of essential spectrum can be carried to some extent onto family $\clos(\X)$. Its definition is based on certain subfamily of closed operators. We say that an operator $\LL\in\clos(\X)$ is \emph{semi-Fredholm} if
\begin{enumerate}[(a)]
\item $\ran(\LL)$ is closed,
\item at least one of the following holds 
\[
\dim\ker\LL<\infty\ \text{and}\ \dim(\X/\ran(\LL))<\infty.
\]
\end{enumerate}
If two conditions from point $(b)$ are satisfied simultaneously we say that $\LL$ is \emph{Fredholm}. For semi-Fredholm operator we define its index as
\[
\ind\LL:=\dim\ker\LL-\dim(\X/\ran(\LL))\in[-\infty,\infty].
\]
The definition of essential spectrum takes into account the fact that semi-Fredholm (Fredholm) operators can be thought of as close to being invertible. It is not, however, unique in the literature so let us briefly give the most prominent definitions. If $\LL\in\clos(\X)$ then its semi-Fredholm, Fredholm, Weyl, Browder domain is the following set
\begin{align*}
&\Delta_{sF}(\LL)=\{\lam\in\bbC{}:\ \LL-\lam\ \text{is semi-Fredholm}\},\\
&\Delta_{F}(\LL)=\{\lam\in\bbC{}:\ \LL-\lam\ \text{is Fredholm}\},\\
&\Delta_{W}(\LL)=\{\lam\in\Delta_{F}:\ \ind(\LL-\lam)=0\},\\
&\Delta_{B}(\LL)=\{\lam\in\Delta_{W}:\ \lam\ \text{is not an accumulation point of}\ \sp(\LL)\}.
\end{align*}
Then we define the semi-Fredholm, etc., essential spectrum of $\LL$ as
\[
\essp^{sF}(\LL)=\bbC{}\bsh\Delta_{sF},\ \essp^{F}(\LL)=\bbC{}\bsh\Delta_{F},\ \essp^{W}(\LL)=\bbC{}\bsh\Delta_{W},\ \essp^{B}(\LL)=\bbC{}\bsh\Delta_{B}.
\]
Of course we have inclusions $\Delta_{B}\subset\Delta_{W}\subset\Delta_{F}\subset\Delta_{sF}$ from which we obtain
\[
\essp^{sF}\subset\essp^{F}\subset\essp^{W}\subset\essp^{B}.
\]
When $\LL$ is a self-adjoint operator on Hilbert space all these definitions coincide with the one given in definition \ref{def:ess_spec} \cite[Rem. 1.11, p. 520]{Kato1995}.
\end{remark}\medskip

\subsection{Formulation of the problem and motivation}\label{ssec:formulation}

Let $\H$ be a complex and separable Hilbert space. Consider the semilinear equation in $\H$
\begin{equation}\label{eq:main}
\LL u + \N(u) = h,
\end{equation}
where $h\in\H$ is a fixed non-zero vector. We will assume that operator $\LL\in\lin(\H)$ is self-adjoint and satisfies
\begin{enumerate}
\renewcommand{\theenumi}{$(\LL_\arabic{enumi})$}
\renewcommand{\labelenumi}{\theenumi}
\item\label{h:spectr} $0\in\essp(L)$,
\item\label{h:gap} $(-\del,0)\subset\rez(L)$\ \text{for some}\ $\del>0$,
\item\label{h:bound}$\inf\sp(\LL)>-\infty$.
\end{enumerate}
In general operator $\N\colon\H\to\H$ is non-linear and such that $\N(0)=0$.\medskip

The basic goal of this work is to give additional conditions for $\LL$ and $\N$ which will guarantee the existence of a solution of equation (\ref{eq:main}). Conditions \ref{h:spectr} - \ref{h:bound} arise mainly when $\LL$ comes from elliptic differential operators considered on an unbounded domain. Before describing this situation more precisely let us recall the case of bounded domains.\medskip

\begin{remark}\label{rem:non-triviality}
Note that when $h\in\H\bsh\{0\}$ the assumptions on operator $\N$ yield that equation (\ref{eq:main}) cannot have trivial (zero) solution. This is important since our methods, based on topological degree, do not guarantee non-triviality of solution.
\end{remark}\medskip

In most applications operators of mathematical physics arising from boundary value problems on bounded subsets of $\bbR[n]$ are self-adjoint and have compact resolvent. Firstly, it means that their spectrum is discrete, i.e. $\sp(\LL)=\dsp(\LL)$ (see definition \ref{def:ess_spec}). Secondly, for every $\lam\notin\sp(\LL)$ the resolvent $\R{\LL}(\lam)$ is compact linear operator, i.e., it takes bounded sets into relatively compact ones. Hence in this case we have precisely two possibilities. If $0\notin\sp(\LL)$ then equation (\ref{eq:main}) can be transformed into fixed point problem
\[
u=\R{\LL}(0)(h-\N(u)).
\]
Due to compactness of the resolvent we can use Schauder Fixed Point Theorem to study the solvability of equation \ref{eq:main} (cf. \cite[Ex. 8.2, p. 60]{Deimling2010} or \cite[Sect. 2.6 \& 6.5]{Zeidler1986}). In the second case we have $0<\dim\ker\LL<\infty$ and we say that equation (\ref{eq:main}) is \emph{at resonance}. Moreover $0$ is always an isolated point of the spectrum (cf. lemma \ref{lem:ess_char}). The conditions which guarantee the solvability are connected with the behaviour of non-linear part $\N$ on $\ker\LL$.

The first result concerning this case for semilinear PDEs was given by Landesman \& Lazer in \cite{LandesmanLazer1970}. The authors considered equation (\ref{eq:main}) in the real Hilbert space $\H=L^2(\Om,\bbR)$ where $\Om\subset\bbR[n]$ is bounded, $\LL$ comes from uniformly elliptic, formally self-adjoint linear differential expression
\[
\LL = \sum_{|\al|,|\be|\leqslant 1}(-1)^{|\al|}D^{\al}\left(a_{\al,\be}D^{\be}\right)
\]
and is defined on $H^{1}_{0}(\Om,\bbR)$. This corresponds to the Dirichlet boundary value problem on $\Om$. The non-linear part $\N$ was the superposition operator induced by a bounded Caratheodory function $f\colon\Om\times\bbR\to\bbR$, i.e., $N(u)(x)=f(x,u(x))$ for every $u\in L^2(\Om,\bbR)$, such that the limits
\[
f_{\pm}(x):=\lim_{s\to\pm\infty}f(x,s)
\]
exist for almost all $x\in\Om$ and belong to $L^2(\Om,\bbR)$. If we assume that $0$ is a simple eigenvalue, i.e. $\ker\LL=1$, and we fix $w\in L^2(\Om,\bbR)$ with norm one such that $\ker\LL=\linspan\{w\}$ then the sufficient condition for solvability of (\ref{eq:main}) is
\begin{equation}\tag{$\star$}\label{eq:lan-laz}
\int_{[w>0]}(f_{+}(x)-h(x))w(x)\,dx+\int_{[w<0]}(f_{-}(x)-h(x))w(x)\,dx>0
\end{equation}
where $[w\gtrless0]=\{x\in\Om:\ w(x)\gtrless 0\}$. Assumption (\ref{eq:lan-laz}) is nowadays called the \emph{Landesman-Lazer condition} and it can be read as stating that the function
\[
\bbR\ni c\mapsto\int_{\Om}(f(x,cw(x))-h(x))w(x)\,dx
\]
changes sign on $\bbR$. The original proof of this result is quite complicated but a simpler argument, which makes use of the perturbation technique, can be found in \cite{Hess1974}.

The more abstract approach, which is close to our philosophy, was developed by Br\'{e}zis \& Nirenberg in \cite{BrezisNirenberg1978}. The authors provided some general methods for the study of solvability of equation (\ref{eq:main}) in a real Hilbert space $\H$ where $\LL\in\clos(\H)$ satisfies
\begin{equation}\tag{I}\label{eq:prop_I}
\ran(\LL)=(\ker\LL)^{\perp}\quad \text{and}\quad (\LL_{|\ran(\LL)})^{-1}\colon\ran(\LL)\to\ran(\LL)\cap\dom(\LL)\ \text{is compact}.
\end{equation}
It is in particular true when $\LL$ is self-adjoint with compact resolvent. The main goal was the geometrical characterization of the range of the left hand side of equation (\ref{eq:main}) in the case when operator $\N$ is monotone and demicontinuous (see section \ref{sec:monot_oper}). The main result, under appropriate assumptions, reads as follows (\cite[Thm. I.1]{BrezisNirenberg1978})
\[
\ran(\LL+\N)\approx\ran(\LL)+\conv\ran(\N)\ (\text{or}\ \ran(\LL)+\ran(\N)),
\]
where $"\approx\,"$ means that the interiors and closures of the sets on both sides are equal. In order to formulate the analytical condition the authors introduced the following functional
\[
J_{\N}(u)=\liminf_{\substack{t\to\infty\\ v\to u}}\langle\N(tv),v\rangle,\ u\in\H,
\]
which they call the \emph{recession function} \cite[Sect. II.1, p. 259]{BrezisNirenberg1978}. In case when $\N$ is demicontinuous and satisfies certain growth assumption the condition
\begin{equation}\tag{$\star\star$}\label{eq:bre-nir}
J_{\N}(u)>\langle h,u\rangle,\quad \text{for all}\ u\in\ker\LL,
\end{equation}
is sufficient for $h\in\inter\ran(\LL+\N)$ when in addition $\dim\ker\LL<\infty$ (\cite[Thm. III.1]{BrezisNirenberg1978}). In particular if $\N$ is a superposition operator in $L^2(\Om,\bbR)$, $\Om\subset\bbR[n]$, generated by function $f$ we have (\cite[Prop. II.4]{BrezisNirenberg1978})
\[
J_{\N}(u)\geqslant\int_{[u>0]}f_{+}(x)u(x)\,dx+\int_{[u<0]}f_{-}(x)u(x)\,dx,
\]
where $f_{+}(x)=\liminf_{s\to\infty}f(x,s)$, $f_{-}(x)=\limsup_{s\to-\infty}f(x,s)$. Hence (\ref{eq:bre-nir}) is a generalization of the Landesman-Lazer condition (see also \cite[Thm. IV.5]{BrezisNirenberg1978} for a generalisation of the result from \cite{LandesmanLazer1970}).

As an example of a more recent development let us mention the paper of Berkovits \& Fabry \cite{BerkovitsFabry2005}. The authors studied semilinear equations in Hilbert spaces when linear part satisfies property (\ref{eq:prop_I}) but has an infinite-dimensional kernel. The main tool in their investigations was certain family of mappings involving operators of class $(S_+)$ (see definition \ref{def:mon_type}(\ref{def_mon_type_S})) and the appropriate topological degree for this family. We take the similar approach in chapter \ref{ch:inf_essential} where we also consider some extension of operators of class $(S_+)$ (see definition \ref{def:(S+)_0L}), but different than in \cite{BerkovitsFabry2005}. Then we apply topological degree for this family to prove the existence theorems.

There are also many papers concerned with equation (\ref{eq:main}) in Banach spaces. Here $\LL$ has the non-trivial kernel but is a Fredholm operator. This condition replaces the assumption that $0$ belongs to discrete spectrum and it generalizes lack of essential spectrum which we have in Hilbert space setting (see remark \ref{rem:fredholm}). Some results in this framework are summarized in \cite{Przeradzki1993}.
\bigskip

When we consider semilinear differential equations on an unbounded domains in $\bbR[n]$ the resolvents of associated linear operators are no longer compact. In consequence the essential spectrum of these operators can be nonempty. Hence in this situation operator $\LL$ may not be (continuously) invertible because $0$ belongs to $\essp(\LL)$. This is indicated by condition \ref{h:spectr}. Lemma \ref{lem:ess_char} implies that $0$ is an eigenvalue of infinite multiplicity or an accumulation point of $\sp(\LL)$. In particular, $\ker\LL$ may be trivial and so the Landesman-Lazer type conditions like (\ref{eq:lan-laz}) or (\ref{eq:bre-nir}) are not sufficient to guarantee the solvability of equation (\ref{eq:main}).

As an example we will consider the Schr\"{o}dinger operators on $\LC[n]$ which are the most important operator class in non-relativistic quantum mechanics. They are of the form
\begin{equation}\label{eq:schrodinger_op}
S=-\Del+V
\end{equation}
where $\Del$ is the Laplace operator on whole $\bbR[n]$ and $V\colon\bbR[n]\to\bbR$ is a measurable function. Since the domains of the Laplace operator on $\LC[n]$ and the multiplication operator induced by $V$ can have in general trivial intersection, the proper definition of $S$ and the question of its self-adjointness are by no means trivial. We shall say that a function $V$ belongs to $(L^2+L^{\infty})(\bbR[n],\bbR)$ if it can be expressed in the form of a sum $V=V_{2}+V_{\infty}$, with $V_{2}\in\LR[n]$ and $V_{\infty}\in L^{\infty}(\bbR[n],\bbR)$. %We also define
%\[
%V_{\pm}:\ V_{\pm}(x)=\pm\max\{\pm V(x),0\}.
%\]
For the matter of convenience we will restrict ourselves to the case $n=3$.
\begin{theorem}\label{thm:schrodinger_op}
%\begin{enumerate}[(i)]
%\item
\cite[Thm. 8.2.2, p. 159]{Davies1995} If $V\in (L^{2}+L^{\infty})(\bbR[3],\bbR)$ then operator $S$ defined as an operator sum is self-adjoint and bounded from below with $\dom(S)= H^{2}(\bbR[3],\bbC)$.
%\item\cite[Thm. 14.1.6, p. 446]{BlaExnHav2008}
%If the potential $V$ is such that $V_{+}\in L^{2}_{loc}(\bbR[3],\bbR)$ and $V_{-}\in(L^2+L^{\infty})(\bbR[n],\bbR)$ then operator $S$, defined as a form sum, is selfadjoint and semibounded on $L^{2}(\bbR[3],\bbR)$.
%\end{enumerate}
\end{theorem}

The spectra of Schr\"{o}dinger operators can have a very diverse shape. We will only concentrate on points of essential spectrum lying on a boundary of a spectral gap. This is exactly the situation described abstractly by conditions \ref{h:spectr} and \ref{h:gap} at the beginning of this section. Firstly let us concentrate on points lying at the infimum of essential spectrum. These always exist since the whole spectrum is bounded from below (condition \ref{h:bound}). If $0$ is such point it clearly belongs to essential spectrum and below it we have a finite or infinite set of eigenvectors. Since the infinite set converges always to $0$, the gap in the spectrum arises only when there are a finite number of eigenvalues on the negative real axis. This is precisely the situation which we investigate in chapter \ref{ch:inf_essential}. In case of Schr\"{o}dinger operators the situation can depend on the behaviour of potential $V$ at infinity as the following theorem exemplifies.
\begin{theorem}
\cite[Thm. 8.5.1, p. 168]{Davies1995} \& \cite[Thm. 14.3.9, p. 459]{BlaExnHav2008} Let $S$ be a Schr\"{o}dinger operator (\ref{eq:schrodinger_op}) on $\LC[3]$ with a potential $V\in\LR[3]$. Then $S$ has an essential spectrum equal to $[0,+\infty)$ and
\begin{itemize}[(i)]
\item if there are $c\in[0,1/4)$ and $r_0>0$ such that $V(x)\geqslant-cr^{-2}$ holds for any $r:=|x|\geqslant r_0$, then $\dsp(S)$ is finite;
\item if there are positive $r_0,d,\eps$ such that $V(x)\leqslant-dr^{-2+\eps}$ for all $r\geqslant r_0$, then $\dsp(S)$ is infinite and converges to $0$.
\end{itemize}
\end{theorem}

Another situation when we can deal with point on the boundary of spectral gaps is when the spectrum of Schr\"{o}dinger operator has so called band structure. It means that the essential spectrum consists of a number of disjoint intervals and the discrete spectrum contains a finite or infinite sequences of eigenvalues located between the intervals. Condition \ref{h:gap} is satisfied whenever there is no sequence of eigenvalues converging to a boundary point of one of these intervals. When $0$ is such point then the essential spectrum can clearly lie on the negative real axis and this general situation is treated in chapter \ref{ch:trivial_ker}.

As an example let us mention a particular case for Schr\"{o}dinger operator $S$ when the potential $V$ is a periodic continuous function. Then the spectrum of $S$ has a band structure and is purely continuous (i.e., there are no eigenvalues at all) \cite[Thm. XIII.100, p. 309]{ReedSimon1978}. Under this assumption the existence of solutions to non-linear stationary Schr\"{o}dinger equation
\[
-\Del u(x)+Vu(x)=f(x,u),\quad x\in\bbR[n]
\]
was studied in a number of papers during last twenty years. The case when $0$ is a boundary point of a spectral gap was investigated for example by Bartsch \& Ding in \cite{BartschDing1999}, Willem \& Zou in \cite{WillemZou2003} and Yang, Chen \& Ding in \cite{YanCheDin2010} with the use of variational methods. At the end of chapter \ref{ch:trivial_ker} we show an application of our abstract results to this equation which does not use the periodicity and allows $0$ to be even an eigenvalue of operator $S$.

%In the second chapter we additionally assume that $\LL$ is bounded from below and $0$ is the infinum of its essential spectrum. In this case the spectrum below $-\del$ is discrete and finite and the main assumption about operator $\N$ is its sublinearity and quasi-monotonicity.

%In the third chapter we will allow the spectrum of $\LL$ below $-\del$ to be non-discrete but operator $\N$ will have to satisfy some monotonicity conditions.\medskip

\subsection{The decomposition of a space}\label{ssec:decomposition}

Let $\{E_\mu:\ \mu\in\bbR{}\}$ be the spectral family of an operator $\LL$ which satisfies the assumptions \ref{h:spectr} and \ref{h:gap}. Put
\begin{equation}\label{eq:decomp}
\lH = \ran(E_{-\del}),\quad \rH = \ran(1-E_{0-})=\ker E_{0-},
\end{equation}
and let $\lP,\ \rP\in\proj(\H)$ be associated orthoprojections. 
Then, according to \ref{h:gap} and the properties of $E_\mu$ (see remark \ref{rem:spectral}$(c)$), we have the following decomposition into orthogonal sum
\[
\H=\lH\oplus\rH.
\] 
For $u\in\H$ we will write $u=u_1+u_2$ where $u_1\in\lH$, $u_2\in\rH$. By $\lN$ and $\rN$ we will denote the components of values of operator $\N$ respectively in $\lH$ and $\rH$. Having introduced this decomposition of a Hilbert space $\H$ let us also write down the equation, which we will call the \emph{perturbed equation}. It will be used to solve $(\ref{eq:main})$. Using notations from this and the previous subsection it has the following form:
\begin{equation}\label{eq:perturbed}
\eps u_2+\LL u+\N(u)=h.
\end{equation}
Now let us go back to the introduced decomposition of a space and show its basic properties.
\begin{lemma}\label{lem:decomp}
Assume that $\LL\in\lin(\H)$ is self-adjoint and satisfies \ref{h:spectr} and \ref{h:gap}.
\begin{enumerate}[(i)]
\item Subspaces $\lH$ and $\rH$ are invariant with respect to $\LL$.\medskip
\[
 \hspace{-.7cm}\boxed{\text{Henceforward we will denote by}\ \lL\ \text{and}\ \rL\ \text{the parts of}\ \LL\ \text{respectively in}\ \lin(\lH)\ \text{and}\ \lin(\rH).}
\]
\item $\sp(\lL)=\sp(\LL)\cap(-\infty,-\del]$, i.e., operator $\lL$ is invertible and for $K:=\lL^{-1}$ we have $\|K\|\leqslant\del^{-1}$. If in addition $\LL$ satisfies \ref{h:bound} then $\sp(\lL)\subset[-\ga,-\del]$ hence  $\lL\in\bound(\lH)$ and $\|\lL\|\leqslant\ga$.
\item $\sp(\rL)=\sp(\LL)\cap[0,\infty)$, i.e., operator $\rL$ is non-negative.
\end{enumerate}
\end{lemma}

\begin{myproof}
Ad$(i)$. Firstly we will show that for each $\lam\in\bbR$ the subspace $\H_\lam=\ran(E_\lam)$ \emph{reduces} operator $\LL$, i.e., if  $u\in \dom(\LL)$ then 
\[
E_\lam u\in \dom(\LL)\ \text{and}\ E_\lam\LL u=\LL E_\lam u.
\]
In fact, from theorem \ref{thm:spectral} we know that
\[
\dom(\LL)=\left\{u\in\H:\ \int_{-\infty}^{\infty}|\mu|^2d\langle u,E_\mu u\rangle<\infty\right\}
\]
and $E_\mu E_\lam= E_{\min\{\mu,\lam\}}$ (see also remark \ref{rem:spectral}$(a)$). Accordingly we have
\[
\int_{-\infty}^{\infty}|\mu|^2d\langle u,E_\mu(E_\lam u)\rangle= \int_{-\infty}^{\lam}|\mu|^2d\langle u,E_\mu u\rangle+ \int_{\lam}^{\infty}|\mu|^2d\langle u, E_\lam u\rangle<\infty,
\]
since the first integral on the right hand side is finite due to $u\in \dom(\LL)$ and the second one is zero because we integrate with respect to constant function. Hence $E_\lam u\in\dom(\LL)$ and eventually
\[
\langle v,E_\lam\LL u\rangle=\langle E_\lam v,\LL u\rangle = \int_{-\infty}^{\infty}\mu\, d\langle E_\lam v,E_\mu u\rangle= \int_{-\infty}^{\infty}\mu\, d\langle v,E_\mu E_\lam u\rangle=\langle v,\LL E_\lam u\rangle
\]
for all $v\in\H$.

Returning to decomposition $(\ref{eq:decomp})$ let $u_1\in\lH\cap\dom(\LL)$ and $u_2\in\rH\cap\dom(\LL)$. Then
\[
\LL u_1=\LL E_{-\del} u_1=E_{-\del}\LL u_1\in\lH,
\]
and
\[
E_0\LL u_2= \LL E_0 u_2=0\ \text{so}\ \LL u_2\in\rH.
\]
This means that operators, which are clearly symmetric, given by formulas $\lL u_1=\LL u_1,\ u_1\in\lH\cap\dom(\LL)$, and $\rL u_2=\LL u_2,\ u_2\in\rH\cap\dom(\LL)$, belong to $\lin(\lH)$ and $\lin(\rH)$ respectively.

Ad $(ii)$ i $(iii)$. Put $f_1:=id_{\bbR}\chi_{(-\infty,-\del]},\ f_2:=id_{\bbR}\chi_{(0,\infty)}$. Then, since $f_1+f_2=id_{\sp(\LL)}$, making use of functional calculus given by theorem \ref{thm:spectral} (see also remark \ref{rem:spectral}$(d)$) we have
\[
\LL=f_1(\LL)+f_2(\LL).
\]
According to spectral mapping theorem we have
\begin{align*}
\sp(f_1(\LL))= &\ \cl{f_1(\sp(\LL))}=\sp(\LL)\cap(-\infty,-\del]\cup\{0\},\\
\sp(f_2(\LL))= &\ \cl{f_2(\sp(\LL))}=\sp(\LL)\cap[0,\infty).
\end{align*}
On the other hand, using additivity of Riemann-Stieltjes integral, for all $u\in\dom(\LL)$ and $v\in\H$ we obtain
\[
\langle v,\LL u\rangle=\int_{-\infty}^{\infty}\mu\, d\langle v,E_\mu u\rangle = \int_{-\infty}^{-\del}\mu\, d\langle v,E_\mu u\rangle + \int_{0}^{\infty}\mu\, d\langle v,E_\mu u\rangle.
\]
Note that if $u\in\lH$ then $E_\mu u=E_{-\del}u$ for $\mu\geqslant 0$ (i.e. $\mu\mapsto E_\mu u$ is constant) so the second integral vanishes and if $u\in\rH$ then $E_\mu u=0$ for $\mu\leqslant-\del$, hence the first integral vanishes. This implies that $\lL=f_1(\LL)_{|\lH}$ and $\rL=f_2(\LL)_{|\rH}$ so in particular we have
\begin{align*}
&\sp(\lL)\subset\sp(f_1(\LL))=\sp(\LL)\cap(-\infty,-\del]\cup\{0\},\\
&\sp(\rL)\subset\sp(f_2(\LL))=\sp(\LL)\cap[0,\infty).
\end{align*}
From the above it also follows that for $u_1\in\lH\cap\dom(\LL)$ and any $v_1\in\lH$ we get
\[
\langle v_1,(\lL-\lam)u_1\rangle=\int_{-\infty}^{-\del}(\mu-\lam)\, d\langle v_1,E_\mu u_1\rangle,
\]
for all $\lam\in\bbR$. If $\lam\notin\sp(\LL)\cap(-\infty,-\del]$ then the function $g_\lam\colon \sp(\LL)\cap(-\infty,0)\to\bbR$ given by formula
\[
g_\lam(\mu)=\frac{1}{\mu-\lam},\quad \mu\in\sp(\LL)\cap(-\infty,0),
\]
is continuous and bounded. So we can define an operator $R_1(\lam)\in\bound(\lH)$ satisfying for all $u_1,v_1\in\lH$ the identity
\[
\langle v_1,R_1(\lam)u_1\rangle = \int_{-\infty}^{-\del}g_\lam(\mu)\, d\langle v_1,E_\mu u_1\rangle
\]
Consequently
\[
\langle v_1,(\lL-\lam)R_1(\lam)u_1\rangle = \int_{-\infty}^{-\del}(\mu-\lam)g_\lam(\mu)\, d\langle v_1,E_\mu u_1\rangle = \int_{-\infty}^{-\del}1\, d\langle v_1,E_\mu u_1\rangle=\langle v_1,u_1\rangle.
\]
So $(\lL-\lam)R_1(\lam)=I_{|\lH}$ and similarly $R_1(\lL-\lam)=I_{|\lH\cap \dom(\lL)}$. Hence we see that  $\lam\in\rez(\lL)$ and eventually $\sp(\lL)=\sp(\LL)\cap(-\infty,-\del]$. The equality $\sp(\rL)=\sp(\LL)\cap[0,\infty)$ can be proven in an analogous way. In particular since $\lL$ and $\rL$ are symmetric and have real spectrum it follows from remark \ref{rem:spectrum}$(b)$ that they are self-adjoint.

Finally, let us note that if $-\ga=\inf\sp(\LL)>-\infty$  and $u=E_{-\del}u\in\lH$ then as in the first part of the proof we have
\[
\int_{-\infty}^{\infty}|\mu|^2d\langle u,E_\mu(E_{-\del} u)\rangle= \int_{-\ga}^{-\del}|\mu|^2d\langle u,E_\mu u\rangle\leqslant\max\{|\ga|^2,|\del|^2\}\cdot\|u\|<\infty.
\]
So $u\in\dom(\LL)$ and consequently $\lH\subset\dom(\LL)$. Hence $\dom(\lL)=\lH$ and its spectrum is bounded because $\sp(\lL)=\sp(\LL)\cap(-\infty,-\del]\subset[-\ga,-\del]$. This means that $\lL\in\bound(\lH)$ as remarked in \ref{rem:spectrum}$(c)$.
\end{myproof}
\medskip

\section{Monotone operators}\label{sec:monot_oper}

In this section we will define some classes of non-linear mappings related to the concept of monotonicity in Hilbert space and we will present their basic properties. At the beginning let us recall different kinds of continuity which are used when dealing with non-linear maps.
\begin{definition}\label{def:var_cont}
Let $\H$ be a real or complex Hilbert space and let $G\subset\H$. Let $A\colon G\subset\H\to\H$ and $\{u_k\}_{k\in\bbNo}\subset G$. We say that operator $A$ is
\begin{enumerate}[(1)]
\item \emph{bounded} if the image of every bounded subset of $G$ is bounded,
%\item hemi-ciągłe jeżeli $A(u+t_nv)\weak A(u)$ dla dowolnych $u,v\in G$ oraz $t_n\to 0^+$,
\item \emph{demicontinuous} if $u_k\to u_0$ implies $A(u_k)\weak A(u_0)$,
\item \emph{continuous} if $u_k\to u_0$ implies $A(u_k)\to A(u_0)$,
\item \emph{weakly sequentially continuous} if $u_k\weak u_0$ implies $A(u_k)\weak A(u_0)$,
\item \emph{strongly sequentially continuous} if $u_k\weak u_0$ implies $A(u_k)\to A(u_0)$,
\item \emph{completely continuous} if it is continuous and for every bounded subset $D\subset\H$ the set $\cl{A(G\cap D)}$ is compact in $\H$.
\end{enumerate}
\end{definition}

\begin{remark}
The demicontinuity of the mapping $A$, although it is defined with the use of sequences, is equivalent to the continuity of $A$ from $\H$ endowed with norm topology to $\H$ with its weak topology $\si(\H,\H^\star)$. Weak sequential continuity and strong sequential continuity are essentially weaker notions than continuity from $\si(\H,\H^\star)$ to $\si(\H,\H^\star)$ and from $\si(\H,\H^\star)$ to norm topology respectively.
\end{remark}
\medskip

\subsection{Monotone and maximal monotone operators on Hilbert spaces}\label{sssec:max_mon}

Now we will proceed to monotonicity conditions for non-linear mappings. Throughout this and next subsection we always assume that $\H$ is a real Hilbert space. We will show how to generalize all introduced notions to mappings in complex Hilbert spaces in subsection \ref{sssec:complex_sp}.\medskip

Let $A\colon\H\to\power{\H}$ be a multi-valued operator (multifunction) defined on a real Hilbert space $\H$, where $\power{\H}$ denotes the family of all subsets of $\H$. The \emph{domain} of operator $A$ is the set
\[
\dom(A)=\{u\in\H:\ A(u)\neq\emptyset\},
\]
and its \emph{image} is defined as follows 
\[
\ran(A)=\bigcup_{u\in\H}A(u).
\]
If $A,B\colon\H\to\power{\H}$ and $\al,\be\in\bbR$ for all $u\in\H$ we set
\[
(\al A+\be B)(u)=\{\al v+\be v':\ v\in A(u),\ v'\in B(u)\}
\]
and then $\dom(\al A+\be B)=\dom(A)\cap\dom(B)$. The \emph{inverse operator} $A^{-1}\colon\H\to\power{\H}$ is defined for all $v\in\H$ by formula
\[
A^{-1}(v)=\{u\in\H:\ v\in A(u)\}.
\]
If we identify $A$ with its graph $\graph(A)=\{(u,v)\in\H\oplus\H:\ v\in A(u)\}$ then $A^{-1}$ is the operator whose graph is symmetric with respect to the graph of $A$, i.e., $(u,v)\in\graph(A)$ iff $(v,u)\in\graph(A^{-1})$. Of course we have $\dom(A^{-1})=\ran(A)$.\medskip

\begin{remark}
If for all $u\in\dom(A)$ the set $A(u)$ consists precisely of one element (so it is a singleton) we call $A$ a single-valued operator. In this case we can attribute to $u$ the unique element of $A(u)$, call it $\tilde{A}(u)$, so that we have a map $\tilde{A}\colon\dom(A)\supset\H\to\H$. Henceforward we will not distinguish between $A$ and $\tilde{A}$. 
\end{remark}\medskip

\begin{definition}\label{def:monotony}
We say that operator $A\colon\H\to\power{\H}$ is
\begin{enumerate}[(1)]
\item \emph{monotone} if for all $u,u'\in\dom(A)$
\[
\langle A(u)-A(u'),u-u'\rangle\geqslant 0,
\]
by which we mean that for all $v\in A(u),v'\in A(u')$ we have
\[
\langle v-v',u-u'\rangle\geqslant 0,
\]
\item \emph{strictly monotone} if it satisfies $(1)$ but the equality can only happen if $u=u'$,
\item \emph{strongly monotone} if there is a constant $C>0$ such that for all $u,u'\in\dom(A)$ we have
\[
\langle A(u)-A(u'),u-u'\rangle\geqslant C\|u-u'\|^2.
\]
\end{enumerate}
\end{definition}

The definition of monotone operator is a natural generalisation of a notion of non-decreasing function $f\colon\bbR\to\bbR$. In fact, we can write the condition $f(x)\leqslant f(y)$ if $x\leqslant y$ equivalently as
\[
(f(x)-f(y))\cdot\!(x-y)\geqslant 0,
\]
for all $x,y\in\bbR$, and the replacement of multiplication with a scalar product gives the definition from \ref{def:monotony}$(1)$. Examples of multi-valued functions from $\bbR$ to $\bbR$ are also easy to exhibit. For instance, let $f(x)=0$ if $x<0$, $f(x)=1$ if $x>0$ and let $f(0)$ be a subset of $[0,1]$.

In infinite-dimensional spaces the simplest examples of monotone operators are among linear and single-valued. Indeed, if $\LL\in\lin(H)$ operator $\LL$ is monotone iff it is non-negative, i.e., $\LL\geqslant 0$ (see subsection \ref{sssec:lin_oper}). Other example is concerned with metric projection. Let $C$ be a nonempty and closed subset of $\H$ and for $u\in\H$ let $P(u)$ be the unique element of $C$ that satisfies
\[
\|P(u)-u\|=\inf\{\|u'-u\|:\ u'\in C\}.
\]
Then it can be shown \cite[Ex. 1.2(f), p. 195]{Phelps1997} that $P$  has the following property
\[
\langle P(u)-P(v),u-v\rangle\geqslant\|P(u)-P(v)\|^2,
\]
which of course implies \ref{def:monotony}$(1)$. For more examples see \cite{Phelps1997} or \cite{Brezis1973}.

%\begin{remark}\label{rem:monot_sup}
%One useful property of single-valued monotone operators on Hilbert space, which they share with their one-dimensonal counterparts, is the fact that the supremum of their values over any compact subset is attained at the boundary of this set. More precisely we have \cite[Prop. A.2, p. 313]{BrezisNirenberg1978}: if $A\colon\H\to\H$ is monotone then for each $R>0$
%\[
%\sup_{|u|\leqslant R}\|A(u)\|=\sup_{\|u\|=R}\|A(u)\|.
%\]
%\end{remark}
\medskip

Since the family of monotone operators is inductive with respect to inclusion of graphs the following definition seems natural.

\begin{definition}\label{def:max_mon}
Assume that an operator $A\colon\H\to\power{\H}$ is monotone. Then we say that $A$ is \emph{maximal monotone} if for all $u\in\dom(A)$ and $u_0,v_0\in \H$ the condition
\[
\langle A(u)-v_0,u-u_0\rangle\geqslant 0,
\]
which we understand as
\[
\langle v-v_0,u-u_0\rangle\geqslant 0,\quad \text{for all}\ v\in A(u),
\]
implies that
\[
u_0\in\dom(A)\ \text{and}\ v_0\in A(u_0).
\]
\end{definition}\medskip

Maximal monotonicity means precisely that a given monotone operator does not have non-trivial monotone extension. Indeed, if we had $\langle A(u)-v_0,u-u_0\rangle\geqslant 0$ on $\dom(A)$ but $u_0\notin\dom(A)$ then defining $\tilde{A}\colon\dom(A)\cup\{u_0\}\to\power{\H}$ as $\tilde{A}(u)=A(u)$ for $u\in\dom(A)$, $\tilde{A}(u_0)=\{v_0\}$, we would get the non-trivial monotone extension. If $u_0\in\dom(A)$ but $v_0\notin A(u_0)$, we can define the monotone extension $\tilde{A}\colon\dom(A)\to 2^{\H}$ by formula $\tilde{A}(u)=A(u)$ for $u\in\dom(A)\bsh\{u_0\}$ and $\tilde{A}(u_0)=A(u_0)\cup\{v_0\}$.

Note that a non-decreasing function $f\colon\bbR\to\bbR$ has always the maximal extension $\tilde{f}\colon\bbR\to\power{\bbR}$ which for all $x\in\bbR$ satisfies
\[
\tilde{f}(x)=[f(x^-),f(x^+)],
\]
where
\[
f(x^\pm):=\lim_{t\to x^\pm}f(t).
\]

The characterisation given in next theorem is fundamental in the study and applications of maximal monotone operators.
\begin{theorem}\label{thm:char_max_mon}
Let us assume that operator $A\colon\H\to\power{\H}$ is monotone. The following conditions are equivalent
\begin{enumerate}[$(i)$]
\item $A$ is maximal monotone,
\item $\ran(A+\lam)=\H$ for some $\lam>0$,
\item $\ran(A+\lam)=\H$ for all $\lam>0$.
\end{enumerate}
\end{theorem}
For the proof see for example \cite[Proposition 2.2, p. 23]{Brezis1973}. 
\medskip

We will finish this subsection with some basic facts concerning maximal monotone operators. Firstly, we recall how this notion behaves under various operations on operators.
\begin{lemma}\label{lem:pers_max_mon}
Assume that $A,B\colon\H\to\power{\H}$ are maximal monotone. Then
\begin{enumerate}[$(i)$]
\item for every $\lam>0$ operator $\lam A$ is maximal monotone,
\item the operator $A^{-1}$ is maximal monotone,
\item if \ $\textnormal{int}\dom(A)\cap\dom(B)\neq\emptyset$ operator $A+B$ is maximal monotone and
\[
\cl{\dom(A)\cap\dom(B)}=\cl{\dom(A)}\cap\cl{\dom(B)}.
\]
\end{enumerate}
\end{lemma}
Items $(i)$ and $(ii)$ follow easily from the definition of maximal monotonicity. For the proof of $(iii)$ see \cite[Cor. 2.7, p. 36]{Brezis1973}.
\medskip

Secondly, let us mention an easy result concerning surjectivity of maximal monotone operators (for more see \cite[p. 30-34]{Brezis1973}).
\begin{lemma}\label{lem:sur_max_mon}
If $A\colon\H\to\power{\H}$ is maximal and strongly monotone operator, then $\ran(A)=\H$.
\end{lemma}
\begin{myproof}
Strong monotonicity (see definition \ref{def:monotony}$(iii)$) implies that $A-C$ is also maximal monotone for some constant $C>0$. So conclusion follows from theorem \ref{thm:char_max_mon}.
\end{myproof}
\bigskip

\subsection{Monotone type operators}

Now we will proceed to definition and description of two classes of non-linear maps which belong to so called monotone type operators.

\begin{definition}\label{def:mon_type}
Let $\H$ be a real Hilbert space and assume that $G\subset\H$. We say that operator $A\colon G\subset\H\to\H$ is
\begin{enumerate}[(1)]
\item\label{def_mon_type_S} of class $(S_+)$ on $G$ if for every sequence $\{u_k\}_{k\in\bbNo}\subset G$ conditions
\[
u_k\weak u_0,\quad \limsup_{k\to\infty}\langle A(u_k),u_k-u_0\rangle\leqslant 0
\]
imply that $u_k\to u_0$,
\item quasi-monotone on $G$ if $\{u_k\}_{k\in\bbNo}\subset G$ and $u_k\weak u_0$ implies
\[
\limsup_{k\to\infty}\langle A(u_k),u_k-u_0\rangle\geqslant 0.
\]
\end{enumerate}
\end{definition}
The condition $(S_+)$ was used for the first time by Browder in \cite[Definition 1, p. 651]{Browder1967} and \cite[Definition 1(b), p. 1002]{Browder1970} in the study of  quasi-linear elliptic equations in generalized divergence form (see also the book by Skrypnik \cite{Skrypnik1994}). In Hilbert space framework it generalises the notion of Leray-Schauder map, i.e., completely continuous (see definition \ref{def:var_cont}) perturbation of identity and the appropriate construction of degree function for this class can be performed, see \cite{Oinas2007}. However, in this dissertation we will use a certain modification of condition $(S_+)$ appropriate for the study of densely defined maps. Nevertheless we introduce this concept here to explain how it works in practice and later compare the two notions. The family of quasi-monotone maps is interesting since the condition $(S_+)$ is preserved under such perturbations. For more information on monotone type mappings and their applications to differential equations we refer to \cite{Francu1990}.\medskip 

\begin{remark}\label{rem:S_+}
\begin{enumerate}[(a)]
\item In definition \ref{def:mon_type}$(\ref{def_mon_type_S})$ it is enough to demand that in consequent there is only a subsequence of $\{u_k\}_{k\in\bbN}$ which converges in norm to $u_0$. This is because the subsequence argument implies then the convergence of the whole sequence.
\item Every strongly monotone operator $A$ (see definition \ref{def:monotony}) belongs to class $(S_+)$.
\item Let us also note that if $G$ is weakly closed (e.g., it is a ball closed in norm or a whole space $\H$) then in the definition of quasi-monotonicity we can postulate in consequent
\[
\limsup_{k\to\infty}\langle A(u_k)-A(u_0),u_k-u_0\rangle\geqslant 0.
\]
This follows from the fact that
\[
\lim_{k\to\infty}\langle A(u),u_k-u_0\rangle=0
\]
according to weak convergence of $\{u_k\}_{k\in\bbN}$. In this case the sufficient conditions for operator $A$ to be quasi-monotone are, for example, monotonicity or strong continuity.
\item\label{rem:Galerkin_approx} Condition $(S_+)$ allows to improve the convergence of solutions to Galerkin approximations of non-linear equations. To explain it let us consider the following example.

Let $A\colon G\subset\H\to\H$ be a bounded and demicontinuous operator of class $(S_+)$ defined on a closed subset $G$ of separable (real) Hilbert space $\H$. Consider the problem of the existence of solutions of the equation
\[
A(u)=h,
\]
where $h\in\H$ is a fixed vector. Taking advantage of the separability of Hilbert space we can choose a sequence of subspaces $\{F_k\}_{k\in\bbN}$ of $\H$ with properties
\[
\dim F_k=k,\quad F_k\subset F_{k+1},\quad \cl{\bigcup_{k}F_k}=\H.
\]
For example, they can be generated by finite systems of an orthonormal basis in $\H$. Now we can study, with fixed $k\in\bbN$, the solvability in $F_k$ of equations
\[
\langle A(u),v\rangle=\langle h,v\rangle\quad (\forall v\in F_k),
\]
which are called \emph{Galerkin approximations} of initial equation. So in fact, having chosen a basis, we have to solve systems of non-linear equations in $\bbR[k]$.

Let us assume further that $\{u_k\}_{k\in\bbN}\subset G$, $u_k\in F_k$ is a \emph{bounded} sequence of solutions of Galerkin approximations. Thus we can pass to a subsequence (we do not change the enumeration) such that
\[
u_k\weak u\in\H\ \text{and}\ A(u_k)\weak w\in\H.
\]
Then $w=h$. Indeed, for every $v\in\H$ we can choose a sequence $\{v_k\}_{k\in\bbN}$, $v_k\in F_k$ such that $v_k\to v$. Accordingly
\[
\lim_{k\to\infty}\langle A(u_k),v_k\rangle= \lim_{k\to\infty}\langle h,v_k\rangle=\langle h,v\rangle,
\]
since $u_k$ are the solutions to Galerkin approximations, but on the other hand
\[
\lim_{k\to\infty}\langle A(u_k),v_k\rangle= \lim_{k\to\infty}(\langle A(u_k),v_k-v\rangle +\langle A(u_k),v\rangle)= \langle w,v\rangle.
\]
Next, making use of the fact that $u_k\in F_k$, we have
\[
\limsup_{k\to\infty}\langle A(u_k),u_k-u\rangle= \limsup_{k\to\infty}(\langle h,u_k\rangle-\langle A(u_k),u\rangle)=0.
\]
Because $A$ belongs to class $(S_+)$ we deduce that $u_k\to u$ and in consequence that $u\in G$. We have, therefore, improved the convergence of Galerkin approximations and now the demicontinuity of $A$ implies $A(u_k)\weak A(u)$ and hence $A(u)=h$.
\end{enumerate}
\end{remark}
\medskip

\subsection{Complex spaces}\label{sssec:complex_sp}

If $\H$ is a complex Hilbert space with an inner (hermitean) product $\langle\cdot\,,\cdot\rangle$ we can define the \emph{real Hilbert space obtained from $\H$} denoted as $\reH$. It is a vector space $\H$ over the field $\bbR$ of real numbers, we "forget" about multiplication by complex numbers, endowed with the inner product $\langle\cdot\,,\cdot\rangle_r$ given with formula
\[
\langle u,v\rangle_r=\re\langle u,v\rangle,\quad \big(\forall u,v\in\reH\big).
\]
This product induces, however, the same norm. In general $\reH$ has less linear dependences between vectors, for example if $\H$ is $n$-dimensional then $\dim\reH=2n$. We also have more mutually orthogonal vectors in $\reH$ than in $\H$, in particular if $\{e_k\}_{k\in\bbN}\subset\H$ is an orthonormal (o.n. for short) basis, then the sequence $\{f_k\}_{k\in\bbN}$ defined as
\[
f_{2k-1}=e_k,\ f_{2k}=ie_k,\quad k=1,2,\ldots,
\]
is the o.n. basis of $\reH$.

If $A\colon G\subset\H\to\power{\H}$, where $\H$ is complex, we can define monotonicity, maximal monotonicity and monotone type conditions using the same definitions but applied to operator $A\colon G\subset\reH\to\reH$, which basically amounts to taking the inner product $\langle\cdot\,,\cdot\rangle_r$ in appropriate definitions. \textbf{Hence when we talk about monotonicity, etc. of operator $A$ defined on complex space $\H$, we will always mean the respective property when treating $A$ as an operator in real Hilbert space $\reH$ obtained form $\H$}.

\section{Topological degree for densely defined operators}\label{sec:top_deg}

In \cite{KartsatosSkrypnik1999} Kartsatos \& Skrypnik introduced a modification of condition $(S_+)$ for densely defined operators in reflexive Banach spaces in order to study some quasilinear elliptic problems with very mild growth conditions on coefficients in the highest order term. Moreover, the existence of an appropriate topological degree for this class of operators was proved in this article (see also the paper of Berkovits \cite{Berkovits1999}). In this section we will recall the basic definitions and theorem about existence of degree in a framework of Hilbert spaces. 

\begin{definition}\label{def:(S+)_0L}
Let $L\subset\H$ be a dense linear submanifold of a real and separable Hilbert space $\H$ and let $\{e_i\}_{i\in\bbN}\subset L$ be an orthonormal (o.n. for short) basis of $\H$. We say that operator $A:\dom(A)\subset\H\to\H$ is of class $(S_+)_{0,L}$ if the following conditions are satisfied
\begin{enumerate}[(i)]
\item $\displaystyle L\subset \dom(A)$,
\item for every sequence $\{u_k\}_{k\in\bbN}\subset L$ conditions
\[
u_k\weak u_0\in\H,\quad \limsup_{k\to\infty}\langle A(u_k),u_k\rangle\leqslant 0
\]
\[
\lim_{k\to\infty}\langle A(u_k),e_i\rangle=0\quad (\forall i\in\bbN)
\]
imply that
\[
u_k\to u_0,\quad u_0\in \dom(A),\quad A(u_0)=0.
\]
\end{enumerate}

%W przypadku gdy $(ii)$ zachodzi dla $w=0$ to mówimy że $A$ spełnia warunek $(S_+)_{0,L}$. Jeżeli warunek "$(u_k)\subset L$" zastąpimy przez "$(u_k)\subset \dom(A)$" to mówimy że odwzorowanie $A$ spełnia warunek $(S_+)_{\dom(A)}$.
\end{definition}

\begin{remark}\label{rem:(S+)_0L}
\begin{enumerate}[(a)]
\item Note that due to assumption that $L$ is dense in $\H$, we can always find an o.n. basis of $\H$ which is contained in $L$.
\item In the paper of Kartsatos \& Skrypnik \cite{KartsatosSkrypnik1999} the framework was set up in a separable and reflexive Banach space and instead of o.n. bases the authors used a sequence of finite dimensional spaces $\{F_i\}_{i\in\bbN}$ of $L$ whose union is dense in $\H$. We can recover the original definition by putting $F_i:=\linspan\{e_1,\ldots,e_i\}$. It is also clear that if we have a family of subspaces $\{F_i\}_{i\in\bbN}$ of Hilbert space $\H$ satisfying condition $(2.2)$ from \cite{KartsatosSkrypnik1999} we can choose an o.n. basis, as defined above, which generates this family.
\item There is, however, an essential difference between our definition and the original one. Namely, we fix the o.n. basis in advance, so in fact the property of being an operator of class $(S_+)$ depends on particular basis (which our notation does not take into account). It is not so in definition $(2.1)$ of \cite{KartsatosSkrypnik1999} and to be in accordance with it we should have assumed that $(ii)$ from definition \ref{def:(S+)_0L} is valid \emph{for every o.n. basis $\{e_i\}_{i\in\bbN}\subset L$} of space $\H$. This assumption seems crucial to guarantee the independence of topological degree for this class from the particular choice of basis. However, the proof of this independence given in theorem 2.2 of \cite{KartsatosSkrypnik1999} seems to be incorrect, see also discussion in \cite[section 3.4.1]{Oinas2007}. Hence throughout this paper we will comfort ourselves with the degree which is dependent, at least in principle, on a particular choice of basis, as is in a case of A-proper mappings, and we will not further emphasize this feature.
\item Let us finally note that condition
\[
\lim_{k\to\infty}\langle A(u_k),e_i\rangle=0\quad (\forall i\in\bbN)
\]
from definition \ref{def:(S+)_0L}$(ii)$ does not imply that the sequence $\{A(u_k)\}_{k\in\bbN}$ is weakly convergent, since we do not know if it is bounded. It is for example true when operator $A$ is locally bounded, meaning that every $u\in\dom(A)$ has a neighbourhood $U$ such that $A(U\cap\dom(A))$ remains bounded.  
\end{enumerate}
\end{remark}\medskip

Now we will proceed to a definition of admissible maps, which have some additional continuity properties crucial in the construction of the topological degree. We also give a notion of appropriate deformation (homotopy) of mappings that will preserve the values of degree.

\begin{definition}\label{def:admissibility}
Let $L\subset\H$ be a dense linear submanifold of a real and separable Hilbert space $\H$ and let $\{e_i\}_{i\in\bbN}\subset L$ be an o.n. basis of $\H$.
\begin{enumerate}[(1)]
\item\label{def:adm_oper} We say that $A:\dom(A)\subset\H\to\H$ is \emph{an admissible operator of class $(S_+)_{0,L}$} if
	\begin{enumerate}
	\item $A$ is of class $(S_+)_{0,L}$,
	\item for every finite dimensional subspace $F\subset L$ and each $v\in L$ functional $a(F,v)\colon F\to\bbR{}$ given by the formula
	\[
	a(F,v)(u)=\langle A(u),v\rangle
	\]
	is continuous.
	\end{enumerate}
\item\label{def:adm_homot} We say that one-parameter family $\{A^t:\ t\in[0,1]\}$ of operators $A^t\colon\dom(A^t)\subset\H\to\H$ is \emph{an admissible homotopy of class $(S_+)_{0,L}$} if
	\begin{enumerate}
	\item $\displaystyle L\subset \dom(A^t)$ for all $t\in[0,1]$,
	\item for all sequences $\{u_k\}_{k\in\bbN}\subset L$, $\{t_k\}_{k\in\bbN}\subset[0,1]$ conditions
	\[
	u_k\weak u_0\in\H,\quad t_k\to t_0\in[0,1],\quad \limsup_{k\to\infty}\langle A^{t_k}(u_k),u_k\rangle\leqslant 0
	\]
	\[
	\lim_{k\to\infty}\langle A^{t_k}(u_k),e_i\rangle=0\quad (\forall i\in\bbN)
	\]
	imply that
	\[
	u_k\to u_0,\quad u_0\in\dom(A^{t_0}),\quad A^{t_0}(u_0)=0,
	\]
	%Ponadto własność tę posiada każda rodzina $A^t-w_t$, $t\in[0,1]$, gdzie $w_t$ jest ciągłą krzywą w $\H$,
	\item for every finite dimensional subspace $F\subset L$ and each $v\in L$ functional $a(F,v)\colon F\times[0,1]\to\bbR{}$ given by formula 
	\[
	a(F,v)(t,u)=\langle A^t(u),v\rangle
	\]
	is continuous.
	\end{enumerate}
\end{enumerate}
\end{definition}

\begin{remark}\label{rem:admissibility}
\begin{enumerate}[(a)]
\item If $A\in(S_+)_{0,L}$ is demicontinuous, or more generally demicontinuous on the finite dimensional subspaces, i.e., for every subspace $F\subset\H$ such, that $\dim F<+\infty$ operator $A_{|F}$ is demicontinuous, then $A$ is admissible. The same remark refers to homotopies in definition \ref{def:admissibility}(\ref{def:adm_homot}), that is, a sufficient condition for item $(c)$ is the demi-continuity of the map $(u,t)\mapsto A^t(u)$ on $F\times[0,1]$ for every finite dimensional subspace $F\subset\H$.
\item Note that if $\{A^t:\ t\in[0,1]\}$ is an admissible homotopy of class $(S_+)_{0,L}$ then for each fixed $t_0\in[0,1]$ operator $A^{t_0}$ is admissible of class $(S_+)_{0,L}$. It is not the case with the original definition of Kartsatos \& Skrypnik (see \cite[p. 148]{Berkovits1999}), and this difference is caused by a fact of \textit{a priori} fixing a basis in our definitions, as pointed out in remark \ref{rem:(S+)_0L}.
\end{enumerate}
\end{remark}\medskip

Before stating the main result in this section, we consider the question of relation of class $(S_+)_{0,L}$ with previously introduced class $(S_+)$. To this end we first recall a definition introduced in \cite{Berkovits1999} which is useful for a structural study of class $(S_+)_{0,L}$.

\begin{definition}
Let $L\subset\H$ be a dense linear submanifold of a real and separable Hilbert space $\H$. We say that operator $A:\dom(A)\subset\H\to\H$ is of class $(S_+)_{L}$ if for all $h\in\H$ operator $A-h$ satisfies condition $(S_+)_{0,L}$.
\end{definition}

Now using this subclass we have the following facts contained in \cite[Thm. 3.3 \& Thm. 3.4]{Berkovits1999}.

\begin{lemma}\label{lem:(S+)_L-prop}
Let $L\subset\H$ be a dense linear submanifold of a real and separable Hilbert space $\H$ and let $A:\dom(A)\subset\H\to\H$. % be such that $L\subset\dom(A)$. 
\begin{enumerate}[(i)]
\item\label{lem:(S+)_L-bound} If $A$ is bounded then $A$ is of class $(S_+)_{L}$ iff $\dom(A)=\H$ and  $A$ is demicontinuous of class $(S_+)$.
\item\label{lem:(S+)_L-pert} Let $N\colon\H\to\H$ be bounded, demicontinuous and quasi-monotone and let $A\in(S_+)_{0,L}$ ($A\in(S_+)_{L}$). Then $A+N\in(S_+)_{0,L}$ ($A+N\in(S_+)_{L}$). If in addition $A$ is admissible $A+N$ is also admissible.
\end{enumerate}
\end{lemma}

We will not exploit item $(\ref{lem:(S+)_L-bound})$ of the foregoing lemma which we only recall to give connection between the two classes. It shows that the theory presented in this section is not interesting for bounded maps, since it basically reduces to the better known theory of operators of class $(S_+)$. In our setting we will use it to investigate operators having form $A=\LL+\N$ where $\LL$ is an unbounded linear operator in Hilbert space, so the theory of $(S_+)$ mapping does not apply to $\LL$ and we have to use the more general one presented here. Note also that we can perceive our operator $A$ as a non-linear perturbation of $\LL$ and the assumption about quasi-monotonicity of $\N$ seems natural since such perturbations preserve condition $(S_+)_{0,L}$ (lemma \ref{lem:(S+)_L-prop}$(\ref{lem:(S+)_L-pert})$).

Now let us state the theorem about existence of topological degree for admissible maps of class $(S_+)_{0,L}$.  

\begin{theorem}\label{thm:degree_existence}
Let $L\subset\H$ be a dense linear submanifold of a real and separable Hilbert space $\H$ with o.n. basis $\{e_i\}_{i\in\bbN}\subset L$. Put
\begin{align*}
\Xi=\big\{&(A,G,h):\ G\subset\H\ \text{is open and bounded},\\
      & A\ \text{is an admissible operator of class}\ (S_+)_{0,L},\ h\in\H\bsh A(\partial G\cap\dom(A))\big\}
\end{align*}
Then there exists a function $\deg\colon\Xi\to\bbZ$ with the following properties:
\begin{enumerate}[(i)]
\item\label{h:deg_solv} if $\deg(A,G,h)\neq 0$ for some $(A,G,h)\in\Xi$ the equation $A(u)=h$ has at least one solution in $G$;
\item\label{h:deg_add} if $G_1,G_2\subset\H$ are open and disjoint subsets of a set $G\subset\H$, $(A,G,h)\in\Xi$ and $h\notin A[(\cl{G}\bsh(G_1\cup G_2))\cap\dom(A)]$
\[
\deg(A,G,h)=\deg(A,G_1,h)+\deg(A,G_2,h);
\]
\item\label{h:deg_hom} if $\{A^t:\ t\in[0,1]\}$ is an admissible homotopy of class $(S_+)_{0,L}$, $(A^0,G,h^0),(A^1,G,h^1)\in\Xi$ and $h^t\in \H\bsh A^t[\partial G\cap\dom(A^t)]$ for all $t\in(0,1)$
\[
\deg(A^0,G,h^0)=\deg(A^1,G,h^1);
\]
\item\label{h:deg_norm} if $I\colon\H\to\H$ denotes the identity in $\H$ for every $u\in G$
\[
\deg(I-u,G,h)=\left\{\begin{array}{ll}
					1 & h\in G-u\\
					0 & h\notin \cl{G}-u.
				\end{array}\right.
\]
\end{enumerate}
\end{theorem}
\begin{myproof}
It is based on a property of improving Galerkin approximations for admissible operators of class $(S_+)_{0,L}$. We discussed this feature for mappings of class $(S_+)$ in remark \ref{rem:S_+}(\ref{rem:Galerkin_approx}). The degree is constructed with the use of finite dimensional operators having form
\[
A^j(u)=\sum_{i=1}^{j}\langle A(u),e_i\rangle e_i,\quad u\in\linspan\{e_1,\ldots,e_j\},\ j\in\bbN.
\]
In \cite{KartsatosSkrypnik1999} the authors exploit the Brouwer degree for this operators. Berkovits in \cite{Berkovits1999} have noticed that we can equivalently use the Leray-Schauder degree for completely continuous mappings (see definition \ref{def:var_cont})
\[
T^j=I-F^j+A^j\circ F^j\colon\H\to\H,
\]
where for $j\in\bbN$ operator $F^j$ is the linear projection onto $\linspan\{e_1,\ldots,e_j\}$. If $G\subset\H$ is open and $h\in\H\bsh A(\partial G\cap\dom(A))$, for $j$ big enough the Leray-Schauder degree $\deg_{LS}(T^j,G,h)$ is well defined, the sequence $\{\deg_{LS}(T^j,G,h)\}$ stabilises and we have
\begin{equation}\label{eq:degree_def}
\deg(A,G,h)=\lim_{j\to+\infty}\deg_{LS}(T^j,G,h).
\end{equation}
See \cite{KartsatosSkrypnik1999} and \cite{Berkovits1999} for details and further comments.
\end{myproof}

\begin{remark}
\begin{enumerate}[(a)]
\item Note that if assumptions of theorem \ref{thm:degree_existence}$(\ref{h:deg_hom})$ are satisfied then, as we already remarked in \ref{rem:admissibility}, for every $t\in[0,1]$ the operator $A^t$ is admissible of class $(S_+)_{0,L}$. Thus we can strengthen the assertion of this item to
\[
\deg(A^t,G,h^t)\equiv\text{const},\quad \text{for all}\ t\in[0,1].
\]
\item In case when assumptions of \ref{thm:degree_existence}$(\ref{h:deg_hom})$ hold we say simply that operators $A_0$ and $A_1$ are homotopic.
\item Let us point out the fact that affine homotopies of the form $A^t=(1-t)A^0+tA^1,\ t\in[0,1]$, when $(A^0,G,h),\ (A^1,G,h)\in\Xi$ and even when the condition $h\notin A^t[\partial G\cap \dom(A^t)]$ holds for all $t\in[0,1]$, cannot be used in general. In fact, we cannot guarantee, for example, under foregoing assumptions that condition $(b)$ from definition \ref{def:admissibility}$(\ref{def:adm_homot})$ is true for affine homotopy without any further assumptions. This phenomenon is also related to the lack of \textit{conical structure} in class $(S_+)_{0,L}$ which is possessed by class $(S_+)$ and similar ones (see \cite[p. 143]{Berkovits1999}).
\item If we restrict ourselves to operators of class $(S_+)_{L}$ which are locally bounded the degree is unique \cite[section 5]{Berkovits1999} so, in particular, it does not depend on the choice of basis.
\end{enumerate}
\end{remark}\medskip

Let us now cover the case of complex Hilbert spaces. To this end we shall modify the definition of condition $(S_+)_{0,L}$ in the following way.

\begin{definition}[complex version of $(S_+)_{0,L}$]\label{def:(S+)_0L-complex}
Let $L\subset\H$ be a dense linear submanifold of a complex and separable Hilbert space $\H$ and let $\{e_i\}_{i\in\bbN}\subset L$ be an o.n. basis of $\H$. We say that operator $A:\dom(A)\subset\H\to\H$ is of class $(S_+)_{0,L}$ if the following conditions are satisfied
\begin{enumerate}[(i)]
\item $\displaystyle L\subset \dom(A)$,
\item for every sequence $\{u_k\}_{k\in\bbN}\subset L$ conditions
\[
u_k\weak u_0\in\H,\quad \limsup_{k\to\infty}\langle A(u_k),u_k\rangle_r\leqslant 0
\]
\[
\lim_{k\to\infty}\langle A(u_k),e_i\rangle=0\quad (\forall i\in\bbN)
\]
imply that
\[
u_k\to u_0,\quad u_0\in \dom(A),\quad A(u_0)=0.
\]
\end{enumerate}
\end{definition}\medskip

Note that we use the scalar product of the real version $\reH$ of $\H$ only in one place. In this way if $A$ satisfies the foregoing definition it is of course of class $(S_+)_{0,L}$ on the real space $\reH$ with underlying basis $\{f_k\}_{k\in\bbN}$, as defined in subsection \ref{sssec:complex_sp}. We can also modify the definition \ref{def:admissibility} in a similar way, considering the complex valued functionals $a(F,v)$. Since the topology stays unchanged when we pass from a complex space to its real version, theorem \ref{thm:degree_existence} gives us the existence of degree in a complex case.

To close this section let us say something about the case of linear self-adjoint operators in a complex Hilbert space $\H$. If $\LL\in\lin(\H)$ is such an operator it is natural to ask when it belongs to the class $(S_+)_{0,\,\dom(\LL)}$. Suppose that there is $\lam\leqslant 0$ such that $\lam\in\essp(\LL)$. This is equivalent to the fact that there exists a sequence $\{u_k\}_{k\in\bbN}\subset\dom(\LL)$, $\|u_k\|=1$ for each $k\in\bbN$, $u_k\weak 0$ and $(\LL-\lam)u_k\to 0$ \cite[Thm. 5.4.4, p. 173]{BlaExnHav2008}. Then we have
\[
\langle\LL u_k,u_k \rangle=\langle(\LL-\lam)u_k,u_k\rangle+\lam \to\lam\leqslant 0
\]
and \emph{for every o.n. basis} $\{e_i\}_{i\in\bbN}$
\[
\langle\LL u_k,e_i \rangle=\langle(\LL-\lam)u_k,e_i\rangle+\lam\langle u_k,e_i\rangle\to 0,
\]
for each fixed $i\in\bbN$. Thus we see that the antecedent of conditional in definition \ref{def:(S+)_0L-complex}$(ii)$ is satisfied for this particular sequence, however, it cannot of course fulfil the consequent because the convergence in norm is impossible in this case. To sum up, we have shown that when $\essp(\LL)\cap(-\infty,0]\neq\emptyset$ for some self-adjoint operator $\LL$ then it does not belong to class $(S_+)_{0,\,\dom(\LL)}$ (no matter the choice of o.n. basis). This demonstrates that, due to condition \ref{h:spectr}, we cannot apply this theory directly to the linear part in our main equation from subsection \ref{ssec:formulation}.

To give a positive result assume that there is a $C>0$ such that $\sp(\LL)\subset[C,\infty)$. If we take any sequence $\{u_k\}_{k\in\bbN}\subset\dom(\LL)$ such that $u_k\weak u_0\in\H$ and $\limsup\langle\LL u_k,u_k\rangle\leqslant 0$ we have according to remark \ref{rem:spectrum}
\[
0\geqslant\limsup_{k\to\infty}\langle\LL u_k,u_k\rangle\geqslant \limsup_{k\to\infty}C\|u_k\|^2.
\]
Hence, $u_k\to 0=u_0\in\dom(\LL)$, so (no matter the choice of o.n. basis) $\LL\in(S_+)_{0,\,\dom(\LL)}$.

Therefore we have proved the following fact.

\begin{lemma}\label{lem:suff_for_(S+)L}
Assume that $\LL\colon\dom(\LL)\subset\H\to\H$ is a linear self-adjoint operator in a complex and separable Hilbert space $\H$. If $\sp(\LL)\subset[C,\infty)$ for some $C>0$ then $\LL$ belongs to class $(S_+)_{0,\,\dom(\LL)}$. In the case when $\dsp(\LL)\cap(-\infty,0]=\emptyset$ this condition is also necessary.
\end{lemma}

\chapter[Equation $Lu+N(u)=h$ when $\sp(\LL)$ is discrete below zero][Equation $Lu+N(u)=h$ when $\sp(\LL)$ is discrete]{Equation $Lu+N(u)=h$ when $\sp(\LL)$ is discrete below zero}\label{ch:inf_essential}
In this chapter we strengthen condition \ref{h:spectr} assuming that a self-adjoint operator $L\in\lin(\H)$ satisfies
\begin{enumerate}
\renewcommand{\theenumi}{$(\LL'_\arabic{enumi})$}
\renewcommand{\labelenumi}{\theenumi}
\item\label{h:spectr'} $0=\inf\essp(\LL)$%$0\in\essp(\LL)\subseteq[0,+\infty]$,
\renewcommand{\theenumi}{$(\LL_\arabic{enumi})$}
\renewcommand{\labelenumi}{\theenumi}
\item $(-\del,0)\subset\rez(\LL)$,
\item $\inf\sp(\LL)=-\ga$, $\ga\geqslant 0$.
\end{enumerate}
From conditions \ref{h:spectr'}, \ref{h:gap} and \ref{h:bound} it follows that spectrum of $\LL$ below $-\del$ is discrete and finite \cite[Thm. 4.2.5, p. 90]{Davies1995}, thus it is composed of finite number of eigenvalues each with finite multiplicity. Firstly, let us note two facts.

\begin{lemma}\label{lem:decomp_real}
Assume that $\LL\in\lin(\H)$ satisfies \ref{h:spectr'}, \ref{h:gap} and \ref{h:bound}, and let $\H=\lH\oplus\rH$ as defined in subsection \ref{ssec:decomposition}.
\begin{enumerate}[(i)]
\item $\dim\lH<+\infty$
\item If $N:\H\to\H$ is bounded and quasi-monotone then $\rN\colon\H\to\H$ given by $\rN(u)=\N(u)_{2}\in\rH$ for all $u\in\H$ is also quasi-monotone.
\end{enumerate}
\end{lemma}

\begin{myproof}
$(i)$ Since the spectral family $\{E_{\mu}:\ \mu\in\bbR\}$ of $\LL$ is concentrated on its spectrum (see remark \ref{rem:spectral}) we get 
\[
\lH=\ran(E_0)=\bigcup_{1=1}^{d}\ran(E_{\lam_i}-E_{\lam_i-})
\]
where $\sp(\LL)\cap(-\infty,0)=\{\lam_1,\ldots,\lam_d\}$. So the conclusion follows from the fact that
\[
\dim\ran(E_{\lam_i}-E_{\lam_i-})=\dim\ker(\LL-\lam_i)<+\infty
\]
for every $i=1,\ldots,d$.

\noindent$(ii)$ Let $u_{1k}+u_{2k}\weak u_1+u_2$. Then from the assumption we have
\[
\limsup_{k\to\infty}\langle\N(u_k),u_{k}-u\rangle_r\geqslant 0
\]
and $u_{1k}\to u_1$ from $(i)$. Hence
\begin{multline*}
\limsup_{k\to\infty}\langle\rN(u_{1k}+u_{2k}),u_{2k}-u_2\rangle_r=
\limsup_{k\to\infty}\langle\lN(u_{1k}+u_{2k}),u_{1k}-u_1\rangle_r\ +\\ +\limsup_{k\to\infty}\langle\rN(u_{1k}+u_{2k}),u_{2k}-u_2\rangle_r =
\limsup_{k\to\infty}\langle\N(u_k),u_{k}-u\rangle_r\geqslant 0.\qedhere
\end{multline*}
\end{myproof}

\section{Perturbed equation}

Firstly, making use of a topological degree for densely defined operators presented in section \ref{sec:top_deg}, we will prove the continuation theorem for perturbed equation $(\ref{eq:perturbed})$. Then with help of an additional assumption on the growth of non-linear part $\N$ we will show the solvability of this equation.

Throughout this section we use the degree for admissible mappings of class $(S_+)_{0,\dom(\LL)}$ with respect to an o.n. basis $\{e_i\}\subset\dom(\LL)$ such that
\begin{equation}\label{eq:basis}
\lH=\linspan\{e_1,\ldots,e_s\},\quad\rH=\cl{\linspan}\{e_{s+1},\ldots\},
\end{equation}
where $s=\dim\lH$.
\begin{theorem}\label{thm:continuation}
Suppose that $\LL\in\lin(\H)$ is a self-adjoint operator satisfying \ref{h:spectr'}, \ref{h:gap} and \ref{h:bound} and $\N\colon\H\to\H$ is bounded, demicontinuous and quasi-monotone. Let $G\subset\H$ be open and bounded and let $h\in\H$. If for fixed $\eps>0$ we have
\begin{enumerate}[(i)]
\item\label{h:cont_i} $\eps u_2+\LL u+t\lN(u)+\rN(u)\neq th_1+h_2$ for all $u\in\partial G\cap\dom(\LL)$ and $t\in(0,1]$,
\item\label{h:cont_ii} $(\rL+\eps)u_2+\rN(u_2)\neq h_2$ for all $u_2\in\partial G\cap\dom(\rL)$,
\item\label{h:cont_iii} $\deg(\rL+\eps+\rN,G\cap\rH,h_2)\neq 0$,
\end{enumerate}
then the perturbed equation (\ref{eq:perturbed})
\[
\eps u_{2}+\LL u+\N(u)=h
\]
has a solution in $G\cap\dom(\LL)$.
\end{theorem}
From lemma \ref{lem:suff_for_(S+)L} and quasi-monotonicity of $\rN|_{\rH}$ it follows that the degree in foregoing theorem is well defined.\medskip

\begin{myproof}
Step 1. Equation (\ref{eq:perturbed}) is equivalent to equation $A(u_1+u_2)=Kh_1+h_2$ where operator
\[
A\colon\lH\oplus\dom(\rL)=\dom(\LL)\subset\lH\oplus\rH\to\lH\oplus\rH, 
\]
is given by the formula
\[
A(u_1+u_2)=\underbrace{u_1+K\lN(u_1+u_2)}_{A_1(u_1+u_2)}\ +\ \underbrace{(\rL+\eps)u_2+\rN(u_1+u_2)}_{A_2(u_1+u_2)}.
\]
We will show that $A$ is an admissible operator of class $(S_+)_{0,\dom(\LL)}$ with respect to an o.n basis $\{e_i\}_{i\in\bbN}\subset\dom(\LL)$ satisfying (\ref{eq:basis}). To this end assume that for $\{u_{1k}+u_{2k}\}_{k\in\bbNo}$ it holds
\begin{subequations}
\begin{gather}
\dom(\LL)\ni u_{1k}+u_{2k}\weak u_{10}+u_{20}\label{eq:cont1}\\
\limsup_{k\to\infty}\langle A(u_{1k}+u_{2k}),u_{1k}+u_{2k}\rangle_r\leqslant 0\label{eq:cont2}\\
\lim_{k\to\infty}\langle A(u_{1k}+u_{2k}),e_i\rangle= 0,\quad \forall i\in\bbN.\label{eq:cont3}
\end{gather}
\end{subequations}
From (\ref{eq:cont3}) we get
\begin{subequations}
\begin{gather}
\lim_{k\to\infty} \langle A_1(u_{1k}+u_{2k}),e_i\rangle= 0,\quad \forall i\in\{1,\ldots,s\}\label{eq:cont4}\\
\lim_{k\to\infty} \langle A_2(u_{1k}+u_{2k}),e_i\rangle= 0,\quad \forall i\in\{s+1,\ldots\}\label{eq:cont5}.
\end{gather}
\end{subequations}

Therefore (\ref{eq:cont4}) and the finiteness of $\dim\lH$ imply that $A_1(u_{1k}+u_{2k})\to 0$, and in view of (\ref{eq:cont2})
\begin{equation}\label{eq:cont6}
\limsup_{k\to\infty}\langle A_2(u_{1k}+u_{2k}),u_{2k}\rangle_r=\limsup_{k\to\infty}\langle A(u_{1k}+u_{2k}),u_{1k}+u_{2k}\rangle_r\leqslant 0.
\end{equation}
Let us note that according to lemmas \ref{lem:decomp_real}, \ref{lem:suff_for_(S+)L} and \ref{lem:(S+)_L-prop}$(\ref{lem:(S+)_L-pert})$ operator
\[
\dom(A)\ni u_1+u_2\to(\rL+\eps)u_2+\rN(u_1+u_2)\in\rH
\]
belongs to class $(S_+)_{0,\dom(\LL)}$, so from (\ref{eq:cont1}), (\ref{eq:cont5}) and (\ref{eq:cont6}) we get $u_{1k}+u_{2k}\to u_{10}+u_{20}$, $u_{10}+u_{20}\in\dom(\LL)$ and $A_2(u_{10}+u_{20})=0$. Moreover from (\ref{eq:cont4}) it follows that $A_1(u_{10}+u_{20})=0$.

Operator $A_1$ is admissible because $K\lN$ is demicontinuous, and admissibility of $A_2$ is guaranteed by lemma \ref{lem:(S+)_L-prop}.
\medskip

\noindent Step 2. Let us consider the family of operators $A^t\colon\lH\oplus\dom(\rL)\subset\lH\oplus\rH\to\lH\oplus\rH$, $t\in[0,1]$, defined as
\[
A^t(u_1+u_2)=u_1+tKN_1(u_1+u_2)+A_2(u_1+u_2).
\]
Proceeding analogously as in step 1 we can demonstrate that this is an admissible homotopy of class $(S_+)_{0,\dom(\LL)}$ and that $A^0$ is an admissible operator of class $(S_+)_{0,\dom(\LL)}$.
\medskip

\noindent Step 3. According to steps 1 \& 2 and assumptions $(\ref{h:cont_i})$ and $(\ref{h:cont_ii})$ we get from theorem \ref{thm:degree_existence} that
\begin{equation}\label{eq:cont_deg}
\deg(A,G,h)=\deg(A^0,G,h_2).
\end{equation}
Operator $A^0$ is defined as
\[
A^0(u_1+u_2)=u_1+(\rL+\eps)u_2+\rN(u_1+u_2).
\]
We will reduce the degree $\deg(A^0,G,h_2)$ to the subspace $\rH$. To this end recall that in the proof of theorem \ref{thm:degree_existence} we mentioned the formula (\ref{eq:degree_def}) which connects this degree with the Leray-Schauder degree of appropriate completely continuous mappings. In the present case we have
\[
\deg(A^0,G,h_2)=\lim_{j\to+\infty}\deg_{LS}(T^{0j},G,h_2)
\]
where $T^{0j}=I-F^j+A^{0j}\circ F^j$, $F^j$ is the projection onto $\linspan\{e_1,\ldots,e_j\}$ and for each $u\in H$ we have
\[
(A^{0j}\circ F^j)(u)=\sum_{i=1}^{j}\langle A^0(F^ju),e_i\rangle e_i= \sum_{i=1}^{j}\langle (F^ju)_1+(\rL+\eps)(F^ju)_2+N_2(F^ju),e_i\rangle e_i,
\]
where $(F^ju)_1,(F^ju)_2$ are the parts of $F^ju$ in $\lH,\rH$ respectively. Hence for every $j>s$ we have using (\ref{eq:basis}) and the mutual orthogonality of $\lH$ and $\rH$
\begin{multline}\label{eq:reduction}
(F^j-A^{0j}\circ F^j)(u)=F^ju-\sum_{i=1}^{s}\langle (F^ju)_1,e_i\rangle e_i- \sum_{i=s+1}^{j}\langle(\rL+\eps)(F^ju)_2+N_2(F^ju),e_i\rangle e_i=\\ =F^ju-(F^ju)_1-\sum_{i=s+1}^{j}\langle(\rL+\eps)(F^ju)_2+N_2(F^ju),e_i\rangle e_i.
\end{multline}
The right hand side of (\ref{eq:reduction}) belongs to $\rH$ so we can use the reduction property of Leray-Schauder degree (see \cite[Thm. 8.7, p. 59]{Deimling2010}) and denoting $G_2=G\cap\rH$ we get 
\[
\deg(A^0,G,h_2)=\lim_{j\to+\infty}\deg_{LS}(T^{0j},G,h_2)= \lim_{j\to+\infty}\deg_{LS}(T^{0j}_{|G_2},G_2,h_2).
\]
The sequence on the right hand side approximates the degree of $A^0$ relative to the subspace $\rH$ so finally we have
\[
\deg(A^0,G,h_2)=\deg(A^0_{|G_2},G_2,h_2)\neq 0,
\]
according to $(\ref{h:cont_iii})$. The conclusion follows now from equality (\ref{eq:cont_deg}) and the properties of the degree.
\end{myproof}

\begin{theorem}\label{thm:pert_ex}
Assume that $\LL\in\lin(\H)$ is self-adjoint, satisfies \ref{h:spectr'}, \ref{h:gap} and \ref{h:bound}, and $\N\colon\H\to\H$ is bounded, demicontinuous, quasi-monotone with
\[
\lim_{k\to+\infty}\frac{\|N(u_k)\|}{\|u_k\|}=0,
\]
for each $\{u_k\}_{k\in\bbN}\subset\H$ such that $\|u_k\|\to+\infty$. Then for all $0<\eps<1$ and $h\in\H$ equation (\ref{eq:perturbed}) has a solution.
\end{theorem}

\begin{myproof}
Step 1. Let $0<\eps<1$. Firstly, we will show that condition $(\ref{h:cont_i})$ from theorem \ref{thm:continuation} holds for $G=B(0,r)$ when $r>0$ is sufficiently large. Suppose the contrary, that there is a sequence $r_k\to+\infty$ and $u_k\in S(0,r_k)\cap\dom(\LL)$, $t_k\in(0,1]$ such that
\[
\left\{\begin{array}{l}
\lL u_{1k} + t_k\lN(u_k)=t_kh_1\\
(\rL+\eps)u_{2k}+\rN(u_k)=h_2,
\end{array}\right.
\]
for every $k\in\bbN$. We can find $M>0$ so big that for $k\geqslant M$ we have
\[
\frac{\|\lN(u_k)\|}{\|u_k\|},\frac{\|h_1\|}{\|u_k\|}\leqslant\frac{1}{2\sqrt{2}\|K\|}
\]
and therefore 
\[
\|u_{1k}\|^2=t_k^2\|K(\lN(u_k)-h_1)\|^2\leqslant 2\|K\|^2\big(\|\lN(u_k)\|^2+\|h_1\|^2\big)\leqslant\frac{1}{2}\|u_k\|^2
\]
which implies $\|u_{1k}\|^2\leqslant\|u_{2k}\|^2$. But then
\[
\frac{\eps^2}{2}\leqslant\frac{\|(\rL+\eps)u_{2k}\|^2}{2\|u_{2k}\|^2}=\frac{\|h_2-\rN(u_k)\|^2}{\|u_{2k}\|^2+\|u_{2k}\|^2}\leqslant 2\frac{\|h_2\|^2+\|\rN(u_k)\|^2}{\|u_k\|^2}
\]
which violates the sublinearity of $\N$ for $k$ large enough. Thus we see that there is $r_1>0$, $r_1>\frac{2\|h_2\|}{\eps}$, such that for all $r\geqslant r_1$ condition $(\ref{h:cont_i})$ in theorem \ref{thm:continuation} holds with $G=B(0,r)$.\medskip

\noindent Step 2. We will now prove item $(\ref{h:cont_iii})$ from theorem \ref{thm:continuation}. To this end we will consider two homotopies acting accordingly with the following scheme
\[
(\rL+\eps+\rN-h_2)\xrightarrow{\text{Step 3}}(\rL+\eps) \xrightarrow{\text{Step 4}} I_{|\rH}.
\]\medskip

\noindent Step 3. Let $A^t\colon\rH\to\rH$, $t\in(0,1]$, be given with the formula
\[
A^t(u_2)=(\rL+\eps)u_2+t\tilde{\rN}(u_2),
\]
where $\tilde{\rN}=\rN-h_2$. Let us assume that there is a sequence $r_k\to+\infty$ and $u_{2k}\in S(0,r_k)\cap\rH$, $t_k\in[0,1]$ such that
\begin{equation*}
A^{t_k}(u_{2k})=0.
\end{equation*}
Then we get
\[
\eps^2\leqslant\frac{\|(\rL+\eps)u_{2k}\|^2}{\|u_{2k}\|^2}=\frac{t_k^2\|\tilde{\rN}(u_{2k})\|^2}{\|u_{2k}\|^2}\leqslant2\frac{\|\rN(u_{2k})\|^2+\|h_2\|^2}{\|u_{2k}\|^2}
\]
which contradicts the sublinearity of $\N$ for large enough $k$. Hence there is $r_2>0$ such that for all $r\geqslant r_2$ equation $A^t(u)=0$ does not have a solution on $S(0,r)\cap\rH$ for every $t\in(0,1]$. Note that in particular for $t=1$ we get condition $(\ref{h:cont_ii})$ from theorem \ref{thm:continuation}.

Now we will demonstrate that $\{A^t:\ t\in[0,1]\}$ is an admissible homotopy of class $(S_+)_{0,\dom(\rL)}$ on $\rH$. To this end fix an o.n. basis $\{e_i\}_{i\in\bbN}$ of $\H$ and take the sequences $\{u_{2k}\}_{k\in\bbNo}\subset\rH$ and $\{t_k\}_{k\in\bbNo}\subset[0,1]$ such that
\begin{subequations}
\begin{gather}
\dom(\rL)\ni u_{2k}\weak u_{20},\ t_k\to t_0\label{eq:pert_ex1}\\
\limsup_{k\to\infty}\langle A^{t_k}(u_{2k}),u_{2k}\rangle_r\leqslant 0\label{eq:pert_ex2}\\
\lim_{k\to\infty}\langle A^{t_k}(u_{2k}),e_{2i}\rangle= 0,\quad \forall i\in\bbN.\label{eq:pert_ex3}
\end{gather}
\end{subequations}
We can assume, perhaps passing to a subsequence, that $\tilde{\rN}(u_{2k})\weak w\in\rH$. Making use of (\ref{eq:pert_ex2}) we have
\begin{multline}\label{eq:pert_ex4}
\limsup_{k\to\infty}\langle(\rL+\eps)u_{2k}+t_0w,u_{2k}\rangle_r=\limsup_{k\to\infty}\big(\langle A^{t_k}(u_{2k}),u_{2k}\rangle_r +\\
+\langle t_0w-t_k\tilde{\rN}(u_{2k}),u_{2k}\rangle_r\big)
\leqslant\limsup_{k\to\infty}\langle t_0w-t_k\tilde{\rN}(u_{2k}),u_{2k}\rangle_r
\end{multline}
Note that $t_k\tilde{\rN}(u_{2k})\weak t_0w$ so $\langle t_0w-t_k\tilde{\rN}(u_{2k}),u_{20}\rangle\to0$ and because of $(\ref{eq:pert_ex4})$ we have
\begin{multline}\label{eq:pert_ex5}
\limsup_{k\to\infty}\langle(\rL+\eps)u_{2k}+t_0w,u_{2k}\rangle_r\leqslant\limsup_{k\to\infty}\langle t_0w-t_k\tilde{\rN}(u_{2k}),u_{2k}-u_{20}\rangle_r=\\
=-\liminf_{k\to\infty}t_k\langle\tilde{\rN}(u_{2k}),u_{2k}-u_{20}\rangle_r\leqslant 0.
\end{multline}
The last inequality holds because $\N$ is quasi-monotone. Subsequently using (\ref{eq:pert_ex3}) we have for all $i\in\bbN$
\begin{equation}\label{eq:pert_ex6}
\lim_{k\to\infty}\langle(\rL+\eps)u_{2k}+t_0w,e_i\rangle=\lim_{k\to\infty}\langle A^{t_k}(u_{2k}),e_i\rangle + \lim_{k\to\infty}\langle t_0w-t_k\tilde{\rN}(u_{2k}),e_i\rangle=0.
\end{equation}
Thus from $(\ref{eq:pert_ex5})$ and $(\ref{eq:pert_ex6})$, on account of lemma \ref{lem:suff_for_(S+)L}, we get: $u_{2k}\to u_{20}\in \dom(\rL)$ and $(\rL+\eps)u_{20}+t_0w=0$. Finally, taking advantage of demi-continuity of $\rN$ we have $\tilde{\rN}(u_{2k})\weak~\rN(u_{20})-h_2$ and since the weak limit is unique it follows that $\rN(u_{20})-h_2=w$, from which $A^{t_0}(u_{20})=0$.

The continuity condition from definition \ref{def:admissibility} is valid on account of lemmas \ref{lem:(S+)_L-prop} and \ref{lem:suff_for_(S+)L}. Therefore, using theorem \ref{thm:degree_existence}, we have
\begin{equation}\label{eq:pert_ex7}
\deg(\rL+\eps+\rN,B(0,r)\cap\rH,h_2)=\deg(\rL+\eps,B(0,r)\cap\rH,0),
\end{equation}
for every $r\geqslant r_2$.\medskip

\noindent Step 4. Consider now homotopy given by formula
\[
A^t=(1-t)(\rL+\eps)+tI_{|\rH}.
\]
Then for every $u_2\in\rH$
\[
\langle A^t(u_2),u_2\rangle\geqslant(1-t)\eps\|u_2\|^2+t\|u_2\|^2\geqslant\eps\|u_2\|^2
\]
since $\eps\leqslant 1$. The foregoing estimate implies that $A^t(u_2)\neq 0$ when $u_2\neq 0$. We can also check, analogously as in previous step, that $\{A^t:\ t\in[0,1]\}$ is an admissible homotopy of class $(S_+)_{0,\dom(\rL)}$. Hence we have
\begin{equation}\label{eq:pert_ex8}
\deg(\rL+\eps,B(0,r)\cap\rH,0)=\deg(I,B(0,r)\cap\rH,0)=1
\end{equation} 
for every $r>0$.\medskip

\noindent Step 5. To sum up put $r_0=\max\{r_1,r_2\}$. Then for all $r>r_0$ the items $(\ref{h:cont_i})$ and $(\ref{h:cont_ii})$ from theorem \ref{thm:continuation} holds for each $G=B(0,r)$ and because of $(\ref{eq:pert_ex7})$ and $(\ref{eq:pert_ex8})$ we have
\[
\deg(\rL+\eps+\rN,B(0,r)\cap\rH,h_2)=1.
\]
Therefore the conclusion of theorem \ref{thm:continuation} gives the existence of a solution to equation $(\ref{eq:perturbed})$ in $B(0,r)$.
\end{myproof}

\section{The limiting step}

\begin{lemma}\label{lem:limit_sol}
Let $\LL\in\lin(\H)$ and $\N\colon\H\to\H$ satisfy the assumptions of theorem \ref{thm:pert_ex}. Let $u_\eps\in\H$ be a solution of perturbed equation $(\ref{eq:perturbed})$ with $\eps>0$, the existence of which is guaranteed by theorem \ref{thm:pert_ex}. If there is a constant $C>0$ such that $\|u_\eps\|\leqslant C$ for all $\eps>0$ small enough then equation $(\ref{eq:main})$ has a solution.
\end{lemma}

\begin{myproof}
Take any decreasing sequence $\eps_k\to 0$ and let $u_k=u_{\eps_k}$ be a solution of equation (\ref{eq:perturbed}) with $\eps=\eps_k$. We can assume, passing to a subsequence if needed, that $u_k\weak u\in\H$. Then $\eps_ku_k\to 0$. Since $u_k$ is a solution of perturbed equation $(\ref{eq:perturbed})$ we have
\[
\LL u_k=h-\N(u_k)-\eps_ku_{2k},
\]
from which it follows that the sequence $\{\LL u_k\}_{k\in\bbN}$ is bounded, so $\LL u_k\weak w\in\H$. Using symmetry of $\LL$ we get for all $v\in\dom(\LL)$
\[
\langle u,\LL v\rangle=\lim_{k\to\infty}\langle u_k,\LL v\rangle= \lim_{k\to\infty}\langle\LL u_k, v\rangle=\langle w,v\rangle.
\]
So $u\in\dom(\LL^*)$ and $\LL^*u=w$. However $\LL$ is self-adjoint so we finally get that $u\in\dom(\LL)$ and $\LL u=w$. Next we compute
\begin{multline*}
\limsup_{k\to\infty}\langle\N(u_k),u_k-u\rangle=\limsup_{k\to\infty}\langle h-(\LL u_k+\eps_ku_{2k}),u_k-u\rangle=\limsup_{k\to\infty}\langle-\lL u_{1k},u_{1k}-u_1\rangle+\\
+\limsup_{k\to\infty}\langle-(\rL+\eps_k)u_{2k},u_{2k}-u_2\rangle=\limsup_{k\to\infty}\langle-(\rL+\eps_k)u_{2k},u_{2k}-u_2\rangle.
\end{multline*}
The last equality follows from the fact that $u_{1k}\to u_1$ because $\dim\lH$ is finite. Subsequently, taking advantage of the convergence $\eps_ku_{2k}\to 0$ and non-negativity of $\rL$ we get
\begin{equation*}
\limsup_{k\to\infty}\langle\N(u_k),u_k-u\rangle=\limsup_{k\to\infty}\langle-\rL u_{2k},u_{2k}-u_2\rangle= -\liminf_{k\to\infty}\langle\rL(u_{2k}-u_2),u_{2k}-u_2\rangle\leqslant 0.
\end{equation*}
Hence from quasi-monotonicity of operator $\N$ the foregoing estimate implies that $\N(u_k)\weak\N(u)$. Finally, making $k$ in the identity
\[
\eps_ku_{2k}+\LL u_k+\N(u_k)=h
\]
go to infinity we arrive at
\[
\LL u+\N(u)=h.\qedhere
\]
\end{myproof}
\bigskip

Now let us introduce the following notation. %If $\LL\in\clos(\H)$ where $\H$ is a Hilbert space by $\H_\LL$ we will denote the Hilbert space $\big(\dom(L),\langle\cdot\,,\cdot\rangle_\LL\big)$ where
%\[
%\langle u,v\rangle_\LL=\langle u,v\rangle+\langle\LL u,\LL v\rangle,\quad u,v\in\dom(\LL).
%\]
%The convergence of sequences in norm of $\H_\LL$ will be indicated as $\overset{\LL}{\to}$ and in weak topology as $\overset{\LL}{\weak}$.
\begin{definition}\label{def:recession} %\cite{BreNir}
Assume that $\H$ is a Hilbert space and $\N\colon\H\to\H$. Define the functional $J_{\N}\colon\H\to[-\infty,+\infty]$ with formula
\[
J_{\N}(u)=\inf\left\{\liminf_{k\to+\infty}\langle N(t_kv_k),v_k\rangle_r:\ t_k\to+\infty,\ \{v_k\}_{k\in\bbN}\subset\H,\ v_k\weak u\right\}.
\]
\end{definition}

\begin{remark}\label{rem:recession}
\begin{enumerate}[(a)]
\item The functional $J_{\N}$ is an example of so called sequential recession function introduced in \cite[Rem. 2.17, p. 157]{Baiocchi1988}, along with more general topological recession function \cite[Def. 2.2, p. 152]{Baiocchi1988}, to study the abstract minimization problems with non-coercive and non-convex energy functional.
\item Let $\N\colon\H\to\H$ %be as in definition \ref{def:recession}.
Define functional $\psi_\N\colon\H\to\bbR$ as
\[
\psi_\N(u)=\langle\N(u),u\rangle_r.
\]
Then we see that
\[
J_{\N}(u)=\liminf_{\substack{t\to\infty\\ v\weak u}} \frac{\psi_{N}(tv)}{t}.
\]
Hence recession function $J_{\N}(u)$ describes the growth of $\psi_{\N}$ as we are heading to infinity and weakly in the direction $u$. 
\end{enumerate}
\end{remark}

\begin{lemma}\label{lem:recession}
Let $\N\colon\H\to\H$. % be as in definition \ref{def:recession}.
\begin{enumerate}[(1)]
\item $J_{\N}(\lam u)=\lam J_{\N}(u)$, for all $\lam>0$.
\item $J_{\N}(0)\in\{-\infty, 0\}$
\end{enumerate}
\end{lemma}

\begin{myproof}
To prove \textit{(1)} note that for each $t_k\in\bbR$, $v_k\in\H$
\[
\langle\N(t_kv_k),v_k\rangle=\lam\left\langle\N\left(\lam t_k\frac{1}{\lam} v_k\right),\frac{1}{\lam}v_k\right\rangle.
\]
Moreover $t_k\to+\infty$ and $v_k\weak\lam u$ if and only if $\lam t_k\to+\infty$ and $\frac{1}{\lam}v_k\weak u$.

Ad. \textit{(2)}. Taking $v_k=0$ for each $k\in\bbR$ we see that $J_{\N}(0)\leqslant 0$. From \textit{(1)} $J_{\N}(0)=\lam J_{\N}(0)$ for every $\lam>0$. Hence $J_{\N}(0)\in\{-\infty,0\}$.
\end{myproof}

%\noindent{\bf Uwaga.}
%Jeżeli $\LL\in\lin(\H)$ spełnia \ref{h:inv} oraz $\N\colon\H\to\H$ to zawsze zakładamy że $J_{\LL\N}$ jest określony na $\H_\LL$.

\begin{theorem}\label{thm:main_ex}
Let $\LL\in\lin(\H)$ be a self-adjoint operator that  satisfies \ref{h:spectr'}, \ref{h:gap} and \ref{h:bound}, and let $h\in~\H$. Moreover assume that $\N\colon\H\to\H$ is bounded, demicontinuous, quasi-monotone and
\begin{enumerate}[(i)]
\item $\displaystyle\lim_{k\to+\infty}\frac{\|\N(u_k)\|^2}{\|u_k\|}=0$ for each sequence $\{u_k\}_{k\in\bbN}\subset\H$ such that $\|u_k\|\to+\infty$,
\item $\displaystyle \limsup_{k\to\infty}\frac{\langle\N(u_k),u_k\rangle_r}{\|u_k\|}>0$ for all $\{u_k\}_{k\in\bbN}\subset\H$ such that $\|u_k\|\to+\infty$.
%\[
%\|u_k\|\to\infty,\ \frac{u_k}{\|u_k\|}\overset{\LL}{\weak}0
%\]
%as $k\to\infty$.
\item $J_{\N}(u)>\langle h,u\rangle_r$ for all $u\in\ker\LL\cap S(0,1)$,
\end{enumerate}
Then equation (\ref{eq:main}) has a solution.
\end{theorem}

\begin{remark}
Let us note that item \textit{(ii)} of the foregoing theorem already implies that $J_{\N}(u)\geqslant0$ for each $u\neq 0$ such that $J_{\N}(u)>-\infty$. Indeed, let $\eps>0$ and take $u\in\H\bsh\{0\}$. We can chose $v_k\weak u$ and $0<t_k\to+\infty$ such that
\[
\lim_{k\to+\infty}\langle\N(t_kv_k),v_k\rangle<J_{\N}(u)+\eps.
\]
Since $\liminf\|v_k\|\geqslant\|u\|>0$ vectors $u_k:=t_kv_k$ satisfy $\|u_k\|\to+\infty$. Hence from \textit{(ii)} we get
\[
0<\limsup_{k\to+\infty}\frac{\langle\N(u_k),u_k\rangle_r}{\|u_k\|}= \limsup_{k\to+\infty}\frac{\langle\N(t_kv_k),v_k\rangle_r}{\|v_k\|}\leqslant \frac{1}{\|u\|}\lim_{k\to+\infty}\langle\N(t_kv_k),v_k\rangle.
\]
Thus we have
\[
0<J_{\N}(u)+\eps
\]
and the assertion follows because $\eps>0$ was chosen arbitrary.
\end{remark}

\begin{myproof}
Let $\eps_k\to 0$ as $k\to\infty$ and for every $k\in\bbN$ fix a solution $u_k\in\H$ of perturbed equation $(\ref{eq:perturbed})$ with $\eps=\eps_k$ whose existence is guaranteed by theorem \ref{thm:pert_ex}. Due to lemma \ref{lem:limit_sol} it is enough to demonstrate that a sequence $\{u_k\}_{k\in\bbN}$ is bounded in $\H$.\medskip

\noindent Step 1. Firstly, we will show that from boundedness of sequence $\{u_{2k}\}_{k\in\bbN}$ follows the boundedness of $\{u_{1k}\}_{k\in\bbN}$. Indeed, assume that there is $C>0$ such that $\|u_{2k}\|\leqslant C$ for all $k\in\bbN$ and that $\|u_{1k}\|\to\infty$ as $k\to\infty$. Then for $k\in\bbN$ so large that $\|u_{1k}\|\geqslant C$ we have
\[
\frac{\|K\|^{-2}}{2}=\frac{\|K\|^{-2}\|u_{1k}\|^2}{2\|u_{1k}\|^2}\leqslant\frac{\|\lL u_{1k}\|^2}{\|u_{1k}\|^2+C^2}=\frac{\|h_1-\lN(u_k)\|^2}{\|u_{1k}\|^2+C^2} \leqslant 2\frac{\|h_1\|^2+\|\lN(u_k)\|^2}{\|u_k\|^2},
\]
where $K=\lL^{-1}$, which leads to contradiction with $(i)$ for sufficiently large $k$.\medskip

\noindent Step 2. Therefore we assume now that $\|u_{2k}\|\to\infty$ as $k\to\infty$. We will prove that in this case
\begin{equation}\label{eq:main_ex1}
\lim_{k\to\infty}\frac{\|u_{1k}\|}{\|u_{2k}\|}=0.
\end{equation}
To this end let $\eps>0$ and choose $k_0\in\bbN$ so large that for all $k\geqslant k_0$ we have
\[
\frac{\|h_1-\lN(u_k)\|}{\|u_k\|}\leqslant\frac{\eps}{\|K\|}.
\]
Making use of the foregoing estimate we get for $k\geqslant k_0$
\[
\|u_{1k}\|^2=\|K(h_1-\lN(u_k))\|^2\leqslant\|K\|^2\|h_1-\lN(u_k)\|^2\leqslant\eps^2\|u_k\|^2=\eps^2\big(\|u_{1k}\|^2+\|u_{2k}\|^2\big)
\]
hence
\[
\frac{\|u_{1k}\|^2}{\|u_{2k}\|^2}\leqslant\frac{\eps^2}{1-\eps^2}.
\]
Since $\eps>0$ was chosen arbitrary we arrive at (\ref{eq:main_ex1}). In particular, due to sublinearity of $\N$, from (\ref{eq:main_ex1}) the following equality follows
\begin{equation}\label{eq:main_ex1a}
\lim_{k\to\infty}\frac{\|\N(u_k)\|}{\|u_{2k}\|}=0.
\end{equation}
In fact, let $\eps>0$. Then for all $k$ sufficiently large
\begin{equation*}
\frac{\|u_{1k}\|}{\|u_{2k}\|}<1\ \text{and}\ \frac{\|N(u_k)\|}{\|u_k\|}<\eps.
\end{equation*}
From the foregoing estimates we get
\[
\frac{\|N(u_k)\|^2}{2\|u_{2k}\|^2}<\frac{\|N(u_k)\|^2}{\|u_{1k}\|^2+\|u_{2k}\|^2}<\eps^2,
\]
hence
\[
\frac{\|\N(u_k)\|}{\|u_{2k}\|}<\sqrt{2}\eps
\]
which gives (\ref{eq:main_ex1a}). Finally let us prove that
\begin{equation}\label{eq:main_ex1b}
\lim_{k\to\infty}\frac{\|u_{1k}\|^2}{\|u_{2k}\|}=0.
\end{equation}
Let $\eps>0$ and choose $k_0\in\bbN$ so large that for all $k\geqslant k_0$ we have
\[
\frac{\|h_1-\lN(u_k)\|^2}{\|u_k\|}\leqslant\frac{\eps}{\|K\|^2}.
\]
So we get for $k\geqslant k_0$
\[
\|u_{1k}\|^2=\|K(h_1-\lN(u_k))\|^2\leqslant\|K\|^2\|h_1-\lN(u_k)\|^2\leqslant\eps\|u_k\|\leqslant C\eps(\|u_{1k}\|+\|u_{2k}\|)
\]
hence
\[
\frac{\|u_{1k}\|^2}{\|u_{2k}\|}\leqslant C\eps\left(1+\frac{\|u_{1k}\|}{\|u_{2k}\|}\right).
\]
Making $k$ converge to $+\infty$ we get from (\ref{eq:main_ex1})
\[
\limsup_{k\to\infty}\frac{\|u_{1k}\|^2}{\|u_{2k}\|}\leqslant C\eps
\]
and since $\eps>0$ was chosen arbitrary we arrive at (\ref{eq:main_ex1b}) 
\medskip

\noindent Step 3. We can assume that $u_{2k}/\|u_{2k}\|\weak u\in\rH$. Then using step 1 we have
\[
v_k:=\frac{u_k}{\|u_{2k}\|}\weak u.
\]
%Accordingly, since for each $k\in\bbN$ vector $u_k$ is a solution of $(\ref{eq:perturbed})$ with $\eps=\eps_k$,
%\[
%\LL v_k=\frac{\LL u_k}{\|u_{2k}\|}=\frac{h-\N(u_k)-\eps_ku_{2k}}{\|u_{2k}\|}
%\]
%and taking into account (\ref{eq:main_ex1a})
%\[
%\|\LL v_k\|\leqslant \frac{\|h\|}{\|u_{2k}\|}+\frac{\|\N(u_k)\|}{\|u_{2k}\|}+\eps_k\to 0.
%\]
%To sum up $\{v_k\}_{k\in\bbN}\subset\dom(\LL)$ satisfies: $v_k\weak u$ and $\LL v_k\to 0$ as $k\to\infty$. Repeating the argument used in a proof of lemma \ref{lem:limit_sol} we can deduce from these two facts that $u\in\ker\LL$. Moreover, due to the definition of scalar product in $\H_\LL$ it follows that $v_k\overset{\LL}{\weak}u$.\medskip
%\noindent Step 4. 
Note that from \ref{h:bound} and (\ref{eq:main_ex1b}) we get
\begin{equation*}
\limsup_{k\to\infty}\frac{\langle\lN(u_k),u_{1k}\rangle_r}{\|u_{2k}\|}=\limsup_{k\to\infty}\frac{\langle h_1-\lL u_{1k},u_{1k}\rangle_r}{\|u_{2k}\|}\leqslant \limsup_{k\to\infty}\frac{\langle h_1,u_{1k}\rangle_r +\ga\|u_{1k}\|^2}{\|u_{2k}\|}=0.
\end{equation*}
Taking advantage of non-negativity of operator $\LL+\eps_k$ on $\rH$ for each $k\in\bbN$ we also have
\begin{equation*}
\limsup_{k\to\infty}\frac{\langle\rN(u_k),u_{2k}\rangle_r}{\|u_{2k}\|}\leqslant\limsup_{k\to\infty}\frac{\langle(\rL+\eps_k)u_{2k}+\rN(u_k),u_{2k}\rangle_r}{\|u_{2k}\|}=\limsup_{k\to\infty}\frac{\langle h_2,u_{2k}\rangle_r}{\|u_{2k}\|}=
\langle h,u\rangle_r.
\end{equation*}
Therefore from two foregoing estimates it follows that
\begin{equation}\label{eq:main_ex2}
\limsup_{k\to\infty}\frac{\langle\N(u_k),u_{k}\rangle_r}{\|u_{2k}\|}\leqslant \langle h,u\rangle_r.
\end{equation}
Let us consider two cases. If $u=0$ then from (\ref{eq:main_ex2})
\[
\limsup_{k\to\infty}\frac{\langle\N(u_k),u_k\rangle_r}{\|u_{k}\|}\leqslant 0
\]
which contradicts $(ii)$. If $u\neq 0$ then using (\ref{eq:main_ex2}) once more we get
\[
\langle h,u\rangle_r\geqslant \limsup_{k\to\infty}\frac{\langle\N(u_k),u_k\rangle_r}{\|u_{2k}\|}\geqslant \liminf_{k\to\infty}\left\langle\N\left(\|u_{2k}\|\frac{u_k}{\|u_{2k}\|}\right),\frac{u_k}{\|u_{2k}\|}\right\rangle_r\geqslant J_{\N}(u),
\]
which, making use of property \textit{(1)} from lemma \ref{lem:recession}, contradicts $(iii)$ this time.
\end{myproof}

\section{Example}

We will illustrate the results of this chapter with a semi-abstract example. Let $P_B\colon\H\to\H$ be the metric projection on a closed ball $B(0,1)$ defined for each $u\in\H$ as
\[
P_B(u)=\left\{\begin{array}{cr}
u, & \|u\|\leqslant 1,\\
\frac{1}{\|u\|}u, & \|u\|>1.
\end{array}\right.
\]
Define operator $\N\colon\H\to\H$ as
\begin{equation}\label{eq:nonlin_ex}
\N(u)=\phi(\|u\|)P_B(u),\quad u\in\H,
\end{equation}
where $\phi\colon\bbR_+\cup\{0\}\to\bbC$ and $\bbR_+$ denotes the set of positive real numbers.

\begin{lemma}\label{lem:main_example}
Assume that $\phi\colon\bbR_+\cup\{0\}\to\bbR_+\cup\{0\}$ is continuous. Then operator $\N$ is bounded, continuous, quasi-monotone and for any $\{u_k\}_{k\in\bbN}\subset\H$ such that $\|u_k\|\to+\infty$ we have
\begin{equation}\label{eq:sign_cond_ex}
\limsup_{k\to+\infty}\frac{\langle\N(u_k),u_k\rangle_r}{\|u_k\|}=\limsup_{k\to+\infty}\phi(\|u_k\|).
\end{equation}
\end{lemma}

\begin{myproof}
Let us take a vector $u\in\H$ with $\|u\|\leqslant r$ where $r>0$ is fixed. Then
\[
\|N(u)\|=\phi(\|u\|)\|P_B(u)\|\leqslant\phi(\|u\|),
\]
since $P_B$ is a projection on the unit closed ball. The continuity of $\phi$ allows us to make the following estimate
\[
\|N(u)\|\leqslant\sup_{t\in[0,r]}\phi(t).
\]
To show the continuity of $\N$ take a sequence $\{u_k\}_{k\in\bbNo}\subset\H$ such that $u_k\to u_0$. Then $\|u_k\|\to\|u_0\|$ and since $\phi$ and $P_B$ are continuous we get
\[
\N(u_k)\to \N(u_0)
\]
due to continuity of scalar multiplication.

For the proof of quasi-monotonicity let us take any $\{u_k\}_{k\in\bbNo}\subset\H$ satisfying $u_k\weak u_0$. Then
\begin{multline*}
\limsup_{k\to+\infty}\langle\N(u_k),u_k-u_0\rangle_r= \limsup_{k\to+\infty}\phi(\|u_k\|)\big(\langle P_B(u_k),u_k\rangle-\langle P_B(u_k),u_0\rangle_r\big)\geqslant\\
\geqslant\limsup_{k\to+\infty}\phi(\|u_k\|)\|P_B(u_k)\|\big(\|u_k\|-\|u_0\|\big)\geqslant 0
\end{multline*}
since $\phi$ is non-negative and due to weak convergence we have $\liminf\|u_k\|\geqslant\|u_0\|$.

Finally to prove (\ref{eq:sign_cond_ex}) consider $\{u_k\}_{k\in\bbN}\subset\H$ with $\|u_k\|\to+\infty$. Then we have
\[
\limsup_{k\to+\infty}\frac{\langle\N(u_k),u_k\rangle_r}{\|u_k\|}= \limsup_{k\to+\infty}\frac{\phi(\|u_k\|)\langle P_B(u_k),u_k\rangle}{\|u_k\|}= \limsup_{k\to+\infty}\phi(\|u_k\|)\|P_B(u_k)\|
\]
and since $\|P_B(u_k)\|=1$ for all $k$ large enough we get (\ref{eq:sign_cond_ex}).
\end{myproof}

\begin{theorem}
Let $\LL\in\lin(\H)$ be a self-adjoint operator satisfying \ref{h:spectr'}, \ref{h:gap}, \ref{h:bound} and let $h\in\H$. Assume that $\phi\colon\bbR_+\cup\{0\}\to\bbR_+\cup\{0\}$ is continuous such that
\begin{enumerate}[(i)]
\item $\displaystyle \lim_{t\to+\infty}\frac{\phi(t)^2}{t}=0$,
\item $\displaystyle \liminf_{t\to+\infty}\phi(t)>\|h\|$.
\end{enumerate}
Then equation
\[
\LL u +\N(u)=h,
\]
where $\N$ is given by formula (\ref{eq:nonlin_ex}), has a solution.
\end{theorem}

\begin{myproof}
From the assumptions and lemma \ref{lem:main_example} it follows that operator $\N$ is bounded, continuous and satisfies condition \textit{(ii)} of theorem \ref{thm:main_ex}. In order to use this theorem we still have to show that conditions \textit{(i)} and \textit{(iii)} hold.

To prove that \textit{(i)} from theorem \ref{thm:main_ex} is true let us take a sequence $\{u_k\}_{k\in\bbN}\subset\H$ such that $\|u_k\|\to+\infty$. Making use of formula (\ref{eq:nonlin_ex}) we get
\[
\limsup_{k\to+\infty}\frac{\|\N(u_k)\|^2}{\|u_k\|}= \limsup_{k\to+\infty}\frac{\phi(\|u_k\|)^2\|P_B(u_k)\|^2}{\|u_k\|}\leqslant \limsup_{k\to+\infty}\frac{\phi(\|u_k\|)^2}{\|u_k\|}
\]
since $\|P_B(u_k)\|\leqslant 1$ for all $k\in\bbN$. Hence, in view of assumption \textit{(i)} from this lemma, the limit exists and we have
\[
\lim_{k\to+\infty}\frac{\|\N(u_k)\|^2}{\|u_k\|}\leqslant \lim_{k\to+\infty}\frac{\phi(\|u_k\|)^2}{\|u_k\|}=0.
\]
Now we are going to check condition \textit{(iii)} from theorem \ref{thm:main_ex}. To this end let us take $u\in\ker\LL\cap S(0,1)$ and sequences $\{t_k\}\subset\bbR_+$ and $\{v_k\}\subset\H$ such that
$t_k\to+\infty$ and $v_k\weak u$. By definition in (\ref{eq:nonlin_ex}) we have
\begin{equation}\label{eq:main_example1}
\liminf_{k\to+\infty}\langle\N(t_kv_k),v_k\rangle_r= \liminf_{k\to+\infty}\phi(t_k\|v_k\|)\|P_B(t_kv_k)\|\|v_k\|.
\end{equation}
Since $\liminf\|v_k\|\geqslant\|u\|=1$ because of the weak convergence it follows that $t_k\|v_k\|\to+\infty$. Hence $\|P_B(t_kv_k)\|=1$ for all $k$ large enough. Passing to a limit in (\ref{eq:main_example1}) we get
\[
\liminf_{k\to+\infty}\langle\N(t_kv_k),v_k\rangle_r\geqslant \big(\liminf_{t\to+\infty}\phi(t)\big)\|u\| = \liminf_{t\to+\infty}\phi(t).
\]
Finally we have
\[
J_{\N}(u)\geqslant\liminf_{t\to+\infty}\phi(t)>\|h\|\geqslant\langle u,h\rangle
\]
for every $u\in\ker\LL\cap S(0,1)$.
%Firstly we will show that the sequence $\{\|v_k\|\}_{k\in\bbN}$ is asymptotically separated from zero. If this was not the case we would have $\|v_k\|\to 0$, perhaps passing to a subsequence. In particular that would imply
%\begin{equation}\label{eq:main_example2}
%\langle v_k,v\rangle\to 0,\quad \text{for each}\ v\in\dom(\LL).
%\end{equation}
%Because $\{v_k\}$ converges weakly in $\H_{\LL}$ to $u$ and $u\in\ker\LL$ we also have
%\begin{equation}\label{eq:main_example3}
%\langle v_k-u,v\rangle+\langle\LL v_k, Lv\rangle\to 0,\quad \text{for each}\ v\in\dom(\LL).
%\end{equation}
%Hence (\ref{eq:main_example2}) and (\ref{eq:main_example3}) would imply that

\end{myproof}

\chapter[Equation $\LL u + \N(u)=h$ when $\sp(\LL)$ is non-discrete below zero][Equation $\LL u + \N(u)=h$ when $\sp(\LL)$ is non-discrete]{Equation $\LL u + \N(u)=h$ when $\sp(\LL)$ is non-discrete below zero}\label{ch:trivial_ker}

In this chapter we assume that a self-adjoint operator $\LL\in\lin(\H)$ satisfies
\begin{enumerate}
%\renewcommand{\theenumi}{$(\LL''_\arabic{enumi})$}
%\renewcommand{\labelenumi}{\theenumi}
%\item\label{h:spectr"} $0\in\csp(L)$,
\renewcommand{\theenumi}{$(\LL_\arabic{enumi})$}
\renewcommand{\labelenumi}{\theenumi}
\item $0\in\sp(\LL)$
\item $(-\del,0)\subset\rez(L)$\ \text{for some}\ $\del>0$,
\item $\inf\sp(\LL)=-\ga$, $\ga\geqslant 0$.
\end{enumerate}

This is a more general situation than in chapter \ref{ch:inf_essential}. Condition \ref{h:spectr} allows the spectrum below $-\del$ to be non-discrete, so the space $\lH$ may have infinite dimension. The techniques used in this chapter are based on the work of Brezis \& Nirenberg \cite{BrezisNirenberg1978}. We start with a simple consequence of \ref{h:bound}.

\begin{lemma}\label{lem:bound}
For every $u_1\in\lH$ we have
\[
\langle\lL u_1,u_1\rangle\geqslant -\ga/\del^2\|\lL u_1\|^2.
\]
\end{lemma}

\begin{myproof}
Let $u_1\in\lH$ i $v_1=\lL u_1$. From lemma \ref{lem:decomp} we get
\[
\|\lL u_1\|=\|v_1\|\geqslant\del\|Kv_1\|=\del\|u_1\|.
\]
Hence using \ref{h:bound}
\[
\langle\lL u_1,u_1\rangle\geqslant-\ga\|u_1\|^2\geqslant-\ga/\del^2\|\lL u_1\|^2.\qedhere
\]
\end{myproof}

\section{Perturbed equation}

Firstly, we will show the solvability of perturbed equations defined in (\ref{eq:perturbed}). Recall that we use decomposition $\N(u)=\lN(u)+\rN(u)$ of the values of non-linear part $\N$ according to the decomposition of Hilbert space $\H$ given in (\ref{eq:decomp}).

%Let us underline that the following theorem, and lemma \ref{lem:limit_sol_monot} from the next section, are proven under more general assumptions on linear part from subsection \ref{ssec:formulation}. Condition \ref{h:spectr"} is nevertheless essential in the proof of main result of this chapter - theorem \ref{thm:main_ex_monot}.  

\begin{theorem}\label{thm:pert}
Assume that a self-adjoint operator $\LL\in\lin(\H)$ satisfies \ref{h:spectr}, \ref{h:gap} and \ref{h:bound}. Operator $\N\colon\H\to\H$ is bounded, demicontinuous, $\N(0)=0$ and
\begin{enumerate}[(i)]
\item there exists $\al>\ga/\del^2$such that for all $u,u'\in\H$ 
\[
\langle\N(u)-\N(u'),u-u'\rangle_r\geqslant\al\|\lN(u)-\lN(u')\|^2.
\]
\end{enumerate}
Then for each $\eps>0$ and $h\in\H$ the perturbed equation (\ref{eq:perturbed})
\[
\eps u_2+\LL u+\N(u)=h
\]
admits precisely one solution.
\end{theorem}

\begin{myproof}
Recall from subsection \ref{ssec:decomposition} that $P_1$ and $P_2=I-P_1$ are orthoprojections on $\lH$ and $\rH$ respectively. Put
\[
\LL(\eps)=\lL P_1+\eps P_2\in\bound(\H).
\]
The boundedness of $\lL$ follows from \ref{h:bound}. Then, what is easy to check, $\LL(\eps)$ is a symmetric bijection and hence $0\notin\sp(\LL(\eps))$ (see remark \ref{rem:resolvent}). What is more, from the fact that $\sp(\lL)\subset[-\ga,-\del]$ (lemma \ref{lem:decomp}) we get
\[
\sp(\LL(\eps))\subset[-\ga,-\del]\cup\{\eps\}
\]
so $\LL(\eps)^{-1}\in\bound(\H)$ is self-adjoint and due to spectral mapping theorem (see remark \ref{rem:spectral}$(\ref{thm:spectr_mapping})$)
\[
\sp(\LL(\eps)^{-1})\subset[-1/\del,-1/\ga]\cup\{1/\eps\}.
\]
Moreover if $u=\LL(\eps)^{-1} v$ then $v_1=\lL u_1 $ and $v_2=\eps u_2$ so we have
\begin{multline}\label{eq:thm_pert1}
\langle\LL(\eps)^{-1} v, v\rangle=\langle u,\LL(\eps) u\rangle=\langle u_1,\lL u_1\rangle+\eps\|u_2\|^2\geqslant\\
\geqslant1/\eps\|\eps u_2\|^2-\ga/\del^2\|\lL u_1\|^2=1/\eps\|v_2\|^2-\ga/\del^2\|v_1\|^2,
\end{multline}
where in the estimate from below we used lemma \ref{lem:bound}. Note that with introduced notation we can convert equation (\ref{eq:perturbed}) as follows:
\begin{align*}
\eps u_2 +\lL u_1 = h-\N(u)-\rL u_2, && \text{by definition of}\ \LL(\eps)\ \text{we get},\\
\LL(\eps) u=h-\N(u)-\rL u_2, && \text{we invert,}\\
u=\LL(\eps)^{-1}(h-\N(u)-\rL u_2), && \text{and finally,}\\
u=\LL(\eps)^{-1} h-\LL(\eps)^{-1}(\rL P_2+\N)(u).
\end{align*}
Let $v=(\rL P_2+\N)(u)$ and put
\[
A=(\rL P_2+\N)^{-1}\colon\H\supset\dom(A)\to 2^\H.
\]
According to the transformations just performed equation (\ref{eq:perturbed}) is equivalent to
\begin{equation}\label{eq:thm_pert2}
\LL(\eps)^{-1} h\in A(v)+\LL(\eps)^{-1} v,\quad v\in\dom(A).
\end{equation}

Since $\sp(\rL P_2)\subset[0,\infty)$, by lemma \ref{lem:decomp}, for all $\lam>0$ we have $-\lam\in\rez(\rL P_2)$. Theorem \ref{thm:char_max_mon} shows that $\rL P_2$ is maximal monotone. Because of $(i)$ operator $\N$ is in particular monotone. It is also maximal monotone. Indeed, assume that $u_0,v_0\in\H$ satisfy
\[
\langle N(u)-v_0,u-u_0\rangle_r\geqslant 0
\]
for every $u\in\H$. Then for all $u\in\H$ and $t>0$ we have in particular $\langle N(u_0+tu)-v_0,u\rangle_r\geqslant 0$ and making $t$ converge to $0$ we arrive at
\[
\langle N(u_0)-v_0,u\rangle_r\geqslant 0,
\]
since $N$ is demicontinuous. This gives $N(u_0)=v_0$ because $u$ is arbitrary. Hence operator $\N+\rL P_2$ is also maximal monotone (lemma \ref{lem:pers_max_mon}$(iii)$) and in consequence $A$ is maximal monotone (lemma \ref{lem:pers_max_mon}$(ii)$).

Let $v,v'\in\dom(A)$. Then for all $u\in A(v),\ u'\in A(v')$
\begin{multline*}
\langle u-u',v-v'\rangle_r= \langle u-u',\rL u_2+\N(u)-\rL u'_2-\N(u')\rangle_r= \langle u_2-u_2',\rL(u_2-u'_2)\rangle+\\
+\langle u-u',\N(u)-\N(u')\rangle_r\geqslant\al\|\lN(u)-\lN(u')\|^2=\al\|v_1-v'_1\|^2.
\end{multline*}
In the estimate from below we used the non-negativity of operator $\rL$ (see lemma \ref{lem:decomp}) and assumption $(i)$. Hence for all $v,v'\in\dom(A)$
\begin{equation}\label{eq:thm_pert3}
 \langle A(v)-A(v'),v-v'\rangle_r\geqslant\al\|v_1-v'_1\|^2.
\end{equation}
Put $A(\eps)=A+\LL(\eps)^{-1}$. From (\ref{eq:thm_pert1}) and (\ref{eq:thm_pert3}) we get for each $v,v'\in\dom(A)$ that
\begin{multline*}
 \langle A(\eps)(v)-A(\eps)(v'),v-v'\rangle_r= \langle A(v)-A(v'),v-v'\rangle_r+\langle \LL(\eps)^{-1}(v-v'),v-v'\rangle\geqslant\\
\geqslant\al\|v_1-v'_1\|^2+1/\eps\|v_2-v'_2\|^2-\ga/\del^2\|v_1-v'_1\|^2= 1/\eps\|v_2-v'_2\|^2+(\al-\ga/\del^2)\|v_1-v'_1\|^2.
\end{multline*}
So operator $A(\eps)$ is strongly monotone (definition \ref{def:monotony}$(iii)$) with constant $C=\min\{1/\eps,\al-\ga/\del^2\}>0$. In particular it means that it is one-to-one in a sense that if $v\neq v'$ then $A(\eps)(v)\cap A(\eps)(v')=\emptyset$. In fact, if there was a vector $w\in A(\eps)(v)\cap A(\eps)(v')$ then from the last estimate we would get
\[
0=\langle w-w,v-v'\rangle\geqslant C\|v-v'\|^2.
\]

Now we will show that operator $A(\eps)$ is maximal monotone. To this end note that
\[
A(\eps)+\lam=A+(\LL(\eps)^{-1}+\lam)
\]
is surjective for all $\lam>1/\del$. Indeed, operator $A$ is, as we already proved, maximal monotone and $\LL(\eps)^{-1}+1/\del$ also belongs to this class according to theorem \ref{thm:char_max_mon} and the fact that $(-\infty,-1/\del)\subset\rez(\LL(\eps)^{-1})$. From lemma \ref{lem:pers_max_mon}$(iii)$ we get that $A+\LL(\eps)^{-1}+1/\delta$ is maximal monotone so, making use of theorem \ref{thm:char_max_mon} once more, for every $\lambda>1/\delta$ the operator $A+\LL(\eps)^{-1}+\lam$ is surjective. 

To sum up, we showed that operator $A(\eps)$ is strongly and maximal monotone so due to lemma \ref{lem:sur_max_mon} it is surjective. This ensures the existence of solution to (\ref{eq:thm_pert2}) which is unique since set-values of $A(\eps)$ are disjoint, this ends the proof of the theorem.
\end{myproof}

\section{The limiting step}

Now we will show that the boundedness of solution of perturbed equations, when the perturbation parameter $\eps$ is also bounded, is sufficient for the solvability of the main equation.

\begin{lemma}\label{lem:limit_sol_monot}
Assume that $\LL$ and $\N$ satisfy the assumptions of theorem \ref{thm:pert} and let $u_\eps\in\H$ be a solution of perturbed equation (\ref{eq:perturbed}) with $\eps>0$. If there is a constant $C>0$ such that $\|u_\eps\|\leqslant C$ for each sufficiently small $\eps>0$ then equation (\ref{eq:main}) has a solution.
\end{lemma}

\begin{myproof}
Fix a sequence $\eps_k\to 0$ and let $u_k:=u_{\eps_k}\in\H$ be a solution of equation \ref{eq:perturbed} with $\eps=\eps_k$. We can assume, choosing a subsequence if needed, that there is $u\in\H$ such that $u_k\weak u$. Since $\N$ is monotone we have for all $v\in\H$
\[
\langle \N(u_k)-\N(v),u_k-v\rangle_r\geqslant 0,
\]
and further
\[
 \langle h-\eps_ku_{2k}-\LL u_k-\N(v),u_k-v\rangle_r\geqslant 0.
\]
Hence, taking advantage of non-negativity of $\rL$, we have for all $v\in\dom(\LL)$
\begin{multline}\label{eq:lem_pert}
\langle h-\N(v),u_k-v\rangle_r\geqslant
 \langle\LL u_k+\eps_ku_{2k},u_k-v\rangle_r= \langle\lL u_{1k},u_{1k}-v_1\rangle_r+\\
+\langle\rL(u_{2k}-v_2),u_{2k}-v_2\rangle + 
\langle\rL v_2,u_{2k}-v_2\rangle_r+ \langle\eps_ku_{2k},u_{2k}-v_2\rangle_r\geqslant\\
\geqslant\langle\lL u_{1k},u_{1k}-v_1\rangle_r+ 
\langle\rL v_2,u_{2k}-v_2\rangle_r + \langle\eps_ku_{2k},u_{2k}-v_2\rangle_r.
\end{multline}
Since $\eps_ku_{2k}\to 0$ and $u_{2k}\weak u_2$ it follows that
\[
\lim_{k\to\infty}\langle\eps_ku_{2k},u_{2k}-v_2\rangle\to 0
\]
and
\[
\lim_{k\to\infty}\langle\rL v_2,u_{2k}-v_2\rangle= 
\langle\rL v_2,u_{2}-v_2\rangle.
\]
We have to study the convergence of the first term on the right hand side of (\ref{eq:lem_pert}). Firstly, we will show that the sequence $\{\lL u_{1k}\}_{k\in\bbN}$ is convergent in norm. To this end fix $k,l\in\bbN$ and note that making use of assumption $(i)$ from theorem \ref{thm:pert} we have
\begin{multline*}
\|\lL u_{1k}-\lL u_{1l}\|^2=\|\lN(u_k)-\lN(u_l)\|^2\leqslant 1/\al\langle \N(u_k)-\N(u_l),u_k-u_l\rangle_r=\\
=-1/\al \langle\LL u_k-\LL u_l,u_k-u_l\rangle+ 1/\al \langle\eps_lu_{2l}-\eps_ku_{2k},u_{2k}-u_{2l}\rangle_r.
\end{multline*}
The non-negativity of $\rL$, lemma \ref{lem:bound} and condition \ref{h:bound} implies that
\begin{multline*}
 \langle\LL u_k-\LL u_l,u_k-u_l\rangle= \langle\lL(u_{1k}-u_{1l}),u_{1k}-u_{1l}\rangle+ \langle\rL(u_{2k}-u_{2l}),u_{2k}-u_{2l}\rangle\geqslant\\
\geqslant-\ga/\del^2\|\lL(u_{1k}-u_{1l})\|^2.
\end{multline*}
From this two estimates we eventually get
\[
\|\lL u_{1k}-\lL u_{1l}\|^2\leqslant\ga/\al\del^2\|\lL(u_{1k}-u_{1l})\|^2 + C(\eps_k+\eps_l),
\]
and since $\al>\ga/\del^2$ we compute
\[
\lim_{k,l\to\infty}(1-\ga/\al\del^2)\|\lL u_{1k}-\lL u_{1l}\|^2\leqslant \lim_{k,l\to\infty}C(\eps_k+\eps_l)=0.
\]
We know from lemma \ref{lem:decomp} that $\lL\in\bound(\lH)$ so $\lL u_{1k}\weak \lL u_1$, according to the weak convergence of $\{u_{1k}\}_{k\in\bbN}$, and since $\{\lL u_{1k}\}_{k\in\bbN}$ is a Cauchy sequence it follows that $\lL u_{1k}\to\lL u_1$. Therefore we have
\[
\lim_{k\to\infty}\langle\lL u_{1k},u_{1k}-v_1\rangle= 
\langle\lL u_{1},u_{1}-v_1\rangle.
\]
Taking the limit on the both sides in (\ref{eq:lem_pert}) as $k\to\infty$ we arrive at
\begin{equation}\label{eq:lem_pert2}
\langle h-\N(v),u-v\rangle_r\geqslant\langle\lL u_{1},u_{1}-v_1\rangle_r + \langle\rL v_2,u_{2}-v_2\rangle_r,
\end{equation} 
for every $v\in\dom(\LL)$.

Putting $v_1=u_1$ in (\ref{eq:lem_pert2}) we have for all $v_2\in\dom(\rL)$
\[
 \langle\rL v_2+\N(u_1+v_2)-h_2,v_2-u_2\rangle_r\geqslant 0
\]
which, from the maximal monotonicity of $\rL(\cdot)+\N(u_1+\cdot\,)$ (derived analogously as in the proof of theorem \ref{thm:pert}), gives $u_2\in\dom(\rL)$ and
\[
\rL u_2+\N(u_1+u_2)=h_2.
\]

Next, fixing $v_2=u_2$ we get from (\ref{eq:lem_pert2})
\[
\langle h_1-\lL u_1-\lN(v_1+u_2),u_1-v_1\rangle_r\geqslant 0.
\]
If we take $v_1=u_1+tw_1\in\lH$, where $w_1\in\lH$ and $t>0$, we then have
\[
\langle h_1-\lL u_1-\lN(u+tw_1),w_1\rangle_r\geqslant 0.
\]
Making $t$ converge to zero and using demi-continuity of $\N$ it follows that
\[
\langle h_1-\lL u_1-\lN(u),w_1\rangle_r\geqslant 0
\]
which, because $w_1\in\lH$ is arbitrary, implies $\lL u_1+\lN(u)=h_1$.
\end{myproof}

We are in position to prove the main theorem of this chapter. Recall that by $J_{\N}$ we denote the recession functional of operator $\N$ with respect to weak convergence introduced in definition \ref{def:recession}.

\begin{theorem}\label{thm:main_ex_monot}
Assume that a self-adjoint operator $\LL\in\lin(\H)$ satisfies \ref{h:spectr}, \ref{h:gap}, \ref{h:bound} and that $h\in\H$. Let $\N\colon\H\to\H$ be bounded, demicontinuous, $\N(0)=0$ and assume that the following conditions hold \begin{enumerate}[(i)]
\item there exists $\al>\ga/\del^2$ such that for all $u,u'\in\H$
\[
 \langle\N(u)-N(u'),u-u'\rangle_r\geqslant\al\|\N(u)-\N(u')\|^2,
\]
\item $\displaystyle\limsup_{k\to+\infty}\frac{\langle\N(u_k),u_k\rangle_r}{\|u_k\|}>\frac{\ga\|h\|}{\del^2\al-\ga}$ for each sequence $\{u_k\}_{k\in\bbN}\subset\H$ such that $\|u_k\|\to+\infty$,
\item $\displaystyle J_{\N}(u)>\frac{\ga\|h\|}{\del^2\al-\ga}+\langle h,u\rangle_r$, for all $u\in\ker\LL\cap S(0,1)$.
\end{enumerate}
Then equation (\ref{eq:main}) has a solution.
\end{theorem}

\begin{myproof}
Let $\eps_k\to 0$ and let $u_k\in\dom(\LL)$ be a solution of perturbed equation (\ref{eq:perturbed}) with $\eps=\eps_k$. According to lemma \ref{lem:limit_sol_monot} it suffices to show that the sequence $\{u_k\}_{k\in\bbN}$ is bounded.\medskip

\noindent Step 1. Firstly, we will show that from boundedness of sequence $\{u_{2k}\}_{k\in\bbN}$ follows the boundedness of $\{u_{1k}\}_{k\in\bbN}$. Therefore assume that the sequence $\{u_{2k}\}_{k\in\bbN}$ is bounded. Since $u_k,\ k\in\bbN$, is a solution of perturbed equation we have in particular
\[
u_{1k}=K(h_1-\lN)(u_k),
\]
where $K=\lL^{-1}$ belongs to $\bound(\lH)$ (see lemma \ref{lem:decomp}). Thus it is sufficient to show that the sequence $\{\lN(u_k)\}_{k\in\bbN}$ is bounded. We will prove more, i.e, that the sequence $\{\N(u_k)\}_{k\in\bbN}$ is bounded. To this end let us note that multiplying both sides of (\ref{eq:perturbed}) with $u_k$ we get
\[
\langle\eps_ku_{2k}+\LL u_k+\N(u_k),u_k\rangle=\langle h,u_k\rangle,
\]
hence
\[
\langle\N(u_k)-h,u_k\rangle=-\langle\LL u_k,u_k\rangle-\eps_k\|u_{2k}\|^2.
\]
On the other hand, making use of assumption $(i)$ with $u=u_k$ and $u'=0$ we have
\[
\langle\N(u_k)-h,u_k\rangle=\langle\N(u_k),u_k\rangle_r-\langle h,u_k\rangle_r \geqslant \al\|\N(u_k)\|^2-\|h\|\|u_k\|.
\]
Taking together the two above formulas we arrive at
\[
\al\|\N(u_k)\|^2-\|h\|\|u_k\|\leqslant-\langle\LL u_k,u_k\rangle-\eps_k\|u_{2k}\|^2 \leqslant\ga/\del^2\|\lL u_{1k}\|^2-\eps_k\|u_{2k}\|^2,
\]
where, in the last estimate, we used lemma \ref{lem:bound}. Next using $\|u_{1k}\|\leqslant1/\del\|\lL u_{1k}\|$ we compute
\begin{multline*}
\al\|\N(u_k)\|^2\leqslant\ga/\del^2\|\lL u_{1k}\|^2+\|h\|\|u_{1k}\|+\|h\|\|u_{2k}\|-\eps_k\|u_{2k}\|^2\leqslant \ga/\del^2\|\lL u_{1k}\|^2+\\
+1/\del\|h\|\|\lL u_{1k}\|+\|h\|\|u_{2k}\|-\eps_k\|u_{2k}\|^2=\ga/\del^2(\|\lL u_{1k}\|+C)^2+\|h\|\|u_{2k}\|-\\
-\eps_k\|u_{2k}\|^2+C\leqslant\ga/\del^2(\|\N(u_k)\|+C)^2+\|h\|\|u_{2k}\|-\eps_k\|u_{2k}\|^2+C.
\end{multline*}
Taking advantage of Cauchy inequality with $\eps'>0$ such that $(1+2\eps')\ga/\del^2<\al$ we get
\[
\al\|\N(u_k)\|\leqslant(1+2\eps')\ga/\del^2\|\N(u_k)\|^2+\|h\|\|u_{2k}\|-\eps_k\|u_{2k}\|^2+C(\eps')
\]
and finally
\begin{equation}\label{eq:thm_main_ex_monot1}
(\al-(1+2\eps')\ga/\del^2)\|\N(u_k)\|^2\leqslant\|h\|\|u_{2k}\|-\eps_k\|u_{2k}\|^2+C(\eps')
\end{equation}
from which follows the boundedness of $\{\N(u_k)\}_{k\in\bbN}$.\medskip

\noindent Step 2. Now assume that $\|u_{2k}\|\to+\infty$. Then
\begin{equation}\label{eq:thm_main_ex_monot2}
\limsup_{k\to\infty}\frac{\|u_{1k}\|^2}{\|u_{2k}\|}\leqslant\frac{\|h\|}{\del^2\al-\ga}.
\end{equation}
Indeed, since $u_{1k}=K(h_1-\lN)(u_k)$, Cauchy inequality with $\eps'$ and estimate (\ref{eq:thm_main_ex_monot1}) imply that
\begin{multline*}
\frac{\|u_{1k}\|^2}{\|u_{2k}\|}\leqslant\frac{(\|\lN(u_k)\|+\|h_1\|)^2}{\del^2\|u_{2k}\|}\leqslant\frac{(1+2\eps')\|\N(u_k)\|^2+C(\epsilon')\|h\|^2}{\del^2\|u_{2k}\|}\leqslant\\ \leqslant\frac{(1+2\eps')\|h\|}{\del^2\al-(1+2\eps')\ga}-\frac{(1+2\eps')\eps_k}{\del^2\al-(1+2\eps')\ga}\|u_{2k}\|+\frac{C(\eps')(1+\|h\|^2)}{\|u_{2k}\|},
\end{multline*}
and we get
\[
\frac{\|u_{1k}\|^2}{\|u_{2k}\|}\leqslant\frac{(1+2\eps')\|h\|}{\del^2\al-(1+2\eps')\ga}+\frac{C(\eps')(1+\|h\|^2)}{\|u_{2k}\|}.
\]
Making $k$ converge to $\infty$ and then $\eps'$ to $0$ we arrive at (\ref{eq:thm_main_ex_monot2}). In particular (\ref{eq:thm_main_ex_monot2}) implies that
\[
\lim_{k\to\infty}\frac{\|u_{1k}\|}{\|u_{2k}\|}=0.
\]
Thus we can assume that there is $u\in\rH$ such that 
\[
v_k:=\frac{u_{k}}{\|u_{2k}\|}\weak u.
\]
%Moreover we have
%\[
%\|\LL v_k\|=\frac{\|\LL u_k\|}{\|u_{2k}\|}\leqslant\frac{\|h\|}{\|u_{2k}\|}+  \frac{\|\N(u_k)\|}{\|u_{2k}\|}+ \eps_k
%\]
%and making use of (\ref{eq:thm_main_ex_monot1}) once more we deduce that $\LL v_k\to 0$.

%To sum up, the sequence $\{v_k\}_{k\in\bbN}\subset\dom(\LL)$ satisfies $v_k\weak v\in\rH$ and $\LL v_k\to 0$. Repeating the argument used in a proof of lemma \ref{lem:limit_sol} we can deduce from these two facts that $u\in\ker\LL$. Moreover, due to the definition of scalar product in $\H_\LL$ it follows that $v_k\overset{\LL}{\weak}u$.\medskip

\noindent Step 3. Note that from \ref{h:bound} and (\ref{eq:thm_main_ex_monot2}) we have
\begin{equation*}
\limsup_{k\to\infty}\frac{\langle\lN(u_k),u_{1k}\rangle_r}{\|u_{2k}\|}=\limsup_{k\to\infty}\frac{\langle h_1-\lL u_{1k},u_{1k}\rangle_r}{\|u_{2k}\|}\leqslant \limsup_{k\to\infty}\frac{\langle h_1,u_{1k}\rangle_r +\ga\|u_{1k}\|^2}{\|u_{2k}\|}\leqslant\frac{\ga\|h\|}{\del^2\al-\ga}.
\end{equation*}
Since operator $\LL+\eps_k$ is non-negative on $\rH$ we also have
\begin{equation*}
\limsup_{k\to\infty}\frac{\langle\rN(u_k),u_{2k}\rangle_r}{\|u_{2k}\|}\leqslant\limsup_{k\to\infty}\frac{\langle(\rL+\eps_k)u_{2k}+\rN(u_k),u_{2k}\rangle_r}{\|u_{2k}\|}=\limsup_{k\to\infty}\frac{\langle h_2,u_{2k}\rangle_r}{\|u_{2k}\|}=\langle h,u\rangle_r.
\end{equation*}
From this two estimates it follows that
\begin{equation*}
\limsup_{k\to\infty}\frac{\langle\N(u_k),u_{k}\rangle_r}{\|u_{k}\|}\leqslant \frac{\ga\|h\|}{\del^2\al-\ga}+\langle h,u\rangle_r.
\end{equation*}
If $u=0$ then
\[
\limsup_{k\to\infty}\frac{\langle\N(u_k),u_{k}\rangle_r}{\|u_{k}\|}\leqslant \frac{\ga\|h\|}{\del^2\al-\ga} 
\]
which is in contradiction with $(ii)$ and if $u\neq 0$ we get
\[
J_{\N}(u)\leqslant\frac{\ga\|h\|}{\del^2\al-\ga}+\langle h,u\rangle_r
\]
which contradicts $(iii)$ this time.
\end{myproof}

\section{Example}

Let $f\colon\bbR[n]\times\bbR\to\bbR$ be a Caratheodory function and for $u\in\LR[n]$ define
\begin{equation}\label{eq:nonlin_monot_ex}
N(u)(x)=f(x,u(x)).
\end{equation}
Let us also introduce the following lower and upper derivatives of function $f$ with respect to the second variable:
\begin{gather*}
\oline{D}_2f(x,t)=\lim_{\del\to 0^+}\sup\left\{\frac{f(x,s)-f(x,t)}{s-t}:\ 0<|s-t|<\del\right\},\\
\uline{D}_2f(x,t)=\lim_{\del\to 0^+}\inf\left\{\frac{f(x,s)-f(x,t)}{s-t}:\ 0<|s-t|<\del\right\}.
\end{gather*}
The limits always exist but can be infinite. In case when $f(x,\cdot)$ is differentiable for some $x\in\bbR[n]$ these are equal and the same as derivative $D_2f(x,\cdot)$. Throughout this section $\langle\cdot\,,\cdot\rangle_{2}$ and $\|\cdot\|_{2}$ will denote respectively the scalar product and norm in $\LC[n]$, i.e.
\begin{align*}
&\langle u,v\rangle_{2}=\int_{\bbR[n]}u(x)\overline{v(x)}\,dx\\
&\|u\|_{2}=\int_{\bbR[n]}|u(x)|^2\,dx
\end{align*}
for all $u,v\in\LC[n]$.
\begin{lemma}\label{lem:superpos_monot}
Let $f\colon\bbR[n]\times\bbR\to\bbR$ be a Caratheodory function and $f(x,0)=0$ for a.e. $x\in\bbR[n]$. Assume further that
\begin{enumerate}[(i)]
\item for a.e. $x\in\bbR[n]$ and each $t\in\bbR$
\[
|f(x,t)|\leqslant a(x)+bt,
\]
where $a\in\LR[n]$ and $b>0$,
\item for a.e. $x\in\bbR[n]$ the function $f(x,\cdot\,)$ is non-decreasing,
\item there are some $\al,\be>0$ such that for a.e. $x\in\bbR[n]$ and each $t\in\bbR$ we have
\[
\be\leqslant\uline{D}_2f(x,t)\leqslant\oline{D}_2f(x,t)\leqslant\frac{1}{\al}.
\]
\end{enumerate}
Then operator $\N\colon\LR[n]\to\LR[n]$ given by formula (\ref{eq:nonlin_monot_ex}) is well-defined, bounded and continuous. Moreover for every $u,u'\in\LR[n]$ it satisfies
\begin{equation*}
\langle\N(u)-N(u'),u-u'\rangle_{2}\geqslant\al\|\N(u)-\N(u')\|_{2}^2
\end{equation*}
and also for each $u\in\LR[n]$ we have
\[
\langle\N(u),u\rangle_{2}\geqslant\be\|u\|_{2}^2.
\]
\end{lemma}

\begin{remark}
\begin{enumerate}[(a)]
\item Condition $(i)$ is a standard polynomial growth assumption which ensures that the superposition operator $\N$ acts in $\LR[n]$ is bounded and continuous. It is even necessary for this to happen in the class of Caratheodory functions (or more generally for sup-measurable functions, see \cite[Theorem 3.1, p. 67]{AppellZabrejko1990}).
\item Note that for almost every $x\in\bbR[n]$ the function $f(x,\cdot\,)$ cannot have jump discontinuities (since it is Caratheodory), so for these $x$ we have
\[
-\infty<\uline{D}_2f(x,t)\leqslant\oline{D}_2f(x,t)<+\infty
\]
for all $t\in\bbR$.
\end{enumerate}
\end{remark}\medskip

\begin{myproof}
As we already remarked, condition $(i)$ ensures that operator $\N$ is well-defined, continuous and bounded. From $(iii)$ we deduce that for a.e. $x\in\bbR[n]$ and any $t,t'\in\bbR$ we have
\begin{equation}\label{eq:ex_monot1}
|f(x,t')-f(x,t)|\leqslant 1/\al|t'-t|.
\end{equation}
Indeed, assume that $t<t'$ and fix $\eps>0$ and $x\in\bbR[n]$. For each $s\in[t,t']$ choose $\del_s>0$ such that
\[
\sup\left\{\frac{f(x,s)-f(x,s')}{s-s'}:\ 0<|s-s'|<\del_s\right\}<\frac{1}{\al}+\eps.
\] 
Taking a finite subcover of $[t,t']$ from $\{B_{\del_s}(s):\ s\in[t,t']\}$, where $B_\del(s)\subset\bbR$ is a ball centred at $s$ with radius $\del>0$, we can choose a sequence $t=t_0<t_1<\ldots<t_m=t'$ such that $f(x,t_k)-f(x,t_{k-1})<1/\al+\eps$, $k=1,\ldots,m$. Hence we have
\[
f(x,t')-f(x,t)=\sum_{k=1}^{m}[f(x,t_k)-f(x,t_{k-1})]<(1/\al+\eps)(t'-t).
\]
Letting $\eps\to 0$ we arrive at (\ref{eq:ex_monot1}). From this and using the monotonicity of $f(x,\cdot\,)$ we get
\begin{equation}\label{eq:ex_monot2}
\al|f(x,s)-f(x,t)|^2\leqslant (f(x,s)-f(x,t))(s-t).
\end{equation}
In an analogous fashion, making use of a lower estimate for $\uline{D}_2f(x,\cdot\,)$ in $(iii)$, we have for a.e. $x\in\bbR[n]$ and all $t,t'\in\bbR$
\[
|f(x,t')-f(x,t)|\geqslant\be|t'-t|.
\]
Hence taking $t'=0$ in the formula above and using the fact that $f(x,t)t\geqslant 0$ we get in particular that
\begin{equation}\label{eq:ex_monot3}
f(x,t)t\geqslant\be t^2.
\end{equation}
Finally, let us take $u,u'\in\LR[n]$. Firstly, applying (\ref{eq:ex_monot2}) we compute 
\begin{multline*}
\langle\N(u)-N(u'),u-u'\rangle_{2}=\int_{\bbR[n]}(f(x,u(x))-f(x,u'(x))(u(x)-u'(x))\,dx\geqslant\\
\geqslant\int_{\bbR[n]}\al|f(x,u(x))-f(x,u'(x))|^2\,dx=\al\|\N(u)-\N(u')\|_{2}^2,
\end{multline*}
which gives the first condition in the assertion of lemma. Secondly, (\ref{eq:ex_monot3}) implies that
\[
\langle\N(u),u\rangle_{2}=\int_{\bbR[n]}f(x,u(x))u(x)\,dx\geqslant\be\|u\|_{2}^2,
\]
so the former one is also true.
\end{myproof}

\begin{theorem}
Let $S=-\Delta+V\in\lin(\LC[3])$ be a Schr\"{o}dinger operator with potential $V\in (L^2+L^{\infty})(\bbR[3],\bbR)$ (see subsection \ref{ssec:formulation}). Assume that $0\in\sp(S)$ and $(-\del,0)\subset\rez(S)$ for some $\del>0$. Moreover assume that $f\colon\bbR[3]\times\bbR\to\bbR$ satisfies the assumptions of lemma \ref{lem:superpos_monot} with $\al>\inf\sp(S)/\del$. Then for each $h\in\LR[3]$ equation
\[
-\Delta u+V(x)u+f(x,u)=h(x)
\]
has a solution $u\in H^2(\bbR[3]\!,\bbR)$.
\end{theorem}

\begin{myproof}
Under assumptions made on potential $V$ the operator $S$ is well defined on $\dom(S)=H^{2}(\bbR[3],\bbC)$ (see theorem \ref{thm:schrodinger_op}). Moreover conditions \ref{h:spectr} - \ref{h:bound} imposed on linear operator in theorem \ref{thm:main_ex_monot} are satisfied. Now let us define operator $\tilde{\N}\colon\LC[3]\to\LC[3]$ as follows
\[
\tilde{\N}(u)(x)=f(x,\xi(x))-if(x,\eta(x))=\N(\xi)(x)-i\N(\eta)(x)
\]
in the notation of lemma \ref{lem:superpos_monot}, where $u=\xi+i\eta\in\LC[3]$. Then we have
\[
\langle\tilde{\N}(u)-\tilde{\N}(u'),u-u'\rangle_{2,r}=\langle\N(\xi)-N(\xi'),\xi-\xi'\rangle_{2}+\langle\N(\eta)-N(\eta'),\eta-\eta'\rangle_{2}
\]
and
\[
\langle\tilde{\N}(u),u\rangle_{2,r}=\langle\tilde{\N}(\xi),\xi\rangle_{2}+ \langle\tilde{\N}(\eta),\eta\rangle_{2}
\]
for all $u,u'\in\LC[3]$, where $\langle\cdot\,,\cdot\rangle_{2,r}=\re\langle\cdot\,,\cdot\rangle_{2}$. Therefore from lemma \ref{lem:superpos_monot} follows that $\tilde{\N}$ is continuous, bounded and
\[
\langle\tilde{\N}(u)-\tilde{\N}(u'),u-u'\rangle_{2,r}\geqslant\al \left(\|\N(\xi)-\N(\xi')\|_{2}^2+\|\N(\eta)-\N(\eta')\|_{2}^2\right)= \al\|\tilde{\N}(u)-\tilde{\N}(u')\|_{2}^2,
\]
and
\[
\limsup_{\|u\|_{2}\to\infty}\frac{\langle\tilde{\N}(u),u\rangle_{2,r}}{\|u\|_{2}}\geqslant\limsup_{\|u\|_{2}\to\infty}\be\|u\|_{2}=+\infty.
\]
Hence assumptions $(i)$ and $(ii)$ of theorem \ref{thm:main_ex_monot} are satisfied. Next take $0<t_k\to\infty$ and $\{v_k\}\subset\LC[3]$ such that $v_k\weak u\neq 0$. %Using Calderon-Zygmund inequality \cite[Thm. 9.11]{GilbargTrudinger2001} we have for every $r>0$ and $v\in\LR[n]_{S}$
%\[
%\|v\|_{H^2(B(0,r),\bbR)}\leqslant C(\|v\|_{2}+\|Sv\|_{2})=C\|v\|_{\LR[n]_{S}},
%\]
%where $B(0,r)\subset\bbR[n]$ is a ball centred in $0$ of radius $r$. Hence the injection $\LR[n]_{S}\hookrightarrow H^{2}_{loc}(\bbR[n],\bbC)$ is continuous and in particular $\LR[n]_{S}$ is compactly embedded into $L^{2}_{loc}(\bbR[n],\bbC)$. Thus we can assume that $v_k$ converges pointwise to $u$ and since $u\neq 0$ we have $\liminf_{k\to\infty}\|v_k\|_{2}>0$. 
Then $\liminf\|v_k\|\geqslant\|u\|=1$ so $t_k\|v_k\|\to+\infty$. Making use of condition (\ref{eq:ex_monot3}) from the proof of lemma \ref{lem:superpos_monot} we get
\[
\langle\tilde{\N}(t_kv_k),v_k\rangle_{2,r}=\int_{\bbR[3]}f(x,t_k\xi_k(x))\xi_k(x)\,dx+\int_{\bbR[3]}f(x,t_k\eta_k(x))\eta_k(x)\,dx\geqslant\be t_k\|v_k\|^{2}_{2}\to\infty,
\]
where for each $k\in\bbN$ we have $v_k=\xi_k+i\eta_k$. Hence
\[
J_{\tilde{\N}}(u)=+\infty
\]
for each $u\neq 0$ so condition $(iii)$ from theorem \ref{thm:main_ex_monot} is also satisfied. This implies that there is $u=\xi+i\eta\in\dom(S)=H^2(\bbR[3]\!,\bbR)$ such that
\[
S u(x)+\tilde{\N}(u)(x)=h(x)
\]
and in particular $\xi\in H^{2}(\bbR[3],\bbR)$ satisfies
\[
S\xi(x)+f(x,\xi(x))=h(x)
\]
in $\LR[3]$.
\end{myproof}
%bounded from below and has the so called band structure, i.e., is purely continuous and consists of closed disjoint intervals \cite[Thm. XIII.100]{ReedSimon1978}

%\include{PhD_general_case}

\backmatter

\bibliographystyle{plain}
\bibliography{C:/Users/pc/Bibliografia/Articles,C:/Users/pc/Bibliografia/Functional_Analysis,C:/Users/pc/Bibliografia/Nonlinear_Analysis,C:/Users/pc/Bibliografia/PDEs,C:/Users/pc/Bibliografia/Mathematical_Physics}

\begin{thebibliography}{10}

\bibitem{AppellZabrejko1990}
J.~Appell and P.P. Zabrejko.
\newblock {\em {Nonlinear superposition operators}}, volume~95 of {\em
  Cambridge Tracts in Mathematics}.
\newblock Cambridge University Press, Cambridge, 1990.

\bibitem{Baiocchi1988}
C.~Baiocchi, G.~Buttazzo, F.~Gastaldi, and F.~Tomarelli.
\newblock {General existence theorems for unilateral problems in continuum
  mechanics}.
\newblock {\em Arch. Rational Mech. Anal.}, 100(2):149--189, 1988.

\bibitem{BartschDing1999}
T.~Bartsch and Y.~Ding.
\newblock {On a nonlinear Schr\"{o}dinger equation with periodic potential}.
\newblock {\em Math. Ann.}, 313:15--37, 1999.

\bibitem{Berkovits1999}
J.~Berkovits.
\newblock {On the degree theory for densely defined mappings of class
  $(S_{+})_{L}$}.
\newblock {\em Abstr. Appl. Anal.}, 4(3):141--152, 1999.

\bibitem{BerkovitsFabry2005}
J.~Berkovits and Ch. Fabry.
\newblock {An extension of the topological degree in Hilbert space}.
\newblock {\em Abstr. Appl. Anal.}, 2005(6):581--597, 2005.

\bibitem{BlaExnHav2008}
J.~Blank, P.~Exner, and M.~Havl\'{i}\v{c}ek.
\newblock {\em {Hilbert space operators in quantum physics}}.
\newblock {T}heoretical and {M}athematical {P}hysics. Springer, 2nd edition,
  2008.

\bibitem{Brezis1973}
H.~Br\'{e}zis.
\newblock {\em {Operateurs maximaux monotones et semi-groupes de contractions
  dans les espaces de Hilbert}}, volume~5 of {\em North-Holland Mathematics
  Studies}.
\newblock North-Holland Publishing Company, Amsterdam, 1973.

\bibitem{BrezisNirenberg1978}
H.~Br\'{e}zis and L.~Nirenberg.
\newblock {Characterizations of the ranges of some nonlinear operators and
  applications to boundary value problems}.
\newblock {\em Ann. Scuola Norm. Sup. Pisa}, 5(2):225--326, 1978.

\bibitem{Browder1967}
F.E. Browder.
\newblock {Nonlinear eigenvalue problems and Galerkin approximations}.
\newblock {\em Bull. Amer. Math. Soc}, 74(4):651--656, 1968.

\bibitem{Browder1970}
F.E. Browder.
\newblock {Nonlinear Elliptic Boundary Value Problems and the Generalized
  Topological Degree}.
\newblock {\em Bull. Amer. Math. Soc}, 76(5):999--1005, 1970.

\bibitem{Davies1995}
E.~B. Davies.
\newblock {\em {Spectral Theory and Differential Operators}}.
\newblock Cambridge University Press, Cambridge, 1995.

\bibitem{Deimling2010}
K.~Deimling.
\newblock {\em {Nonlinear functional analysis}}.
\newblock Dover Publications, Mineola, New York, {D}over edition, 2010.

\bibitem{Francu1990}
J.~Franc\r{u}.
\newblock {Monotone operators. A survey directed to applications to
  differential equations}.
\newblock {\em Aplikace matematiky}, 35(4):257--301, 1990.

\bibitem{Hess1974}
P.~Hess.
\newblock {On a theorem by Landesman and Lazer}.
\newblock {\em Indiana Univ. Math. J.}, 23(9):827--829, 1974.

\bibitem{KarakostasTsamatos2001}
G.L. Karakostas and P.Ch. Tsamatos.
\newblock {On a Nonlocal Boundary Value Problem at Resonance}.
\newblock {\em J. Math. Anal. Appl.}, 259:209--218, 2001.

\bibitem{KartsatosSkrypnik1999}
A.G. Kartsatos and I.V. Skrypnik.
\newblock {Topological degree theories for densely defined mappings involving
  operators of type $(S_+)$}.
\newblock {\em Adv. Differential Equations}, 4(3):413 -- 456, 1999.

\bibitem{Kato1995}
T.~Kato.
\newblock {\em {Perturbation Theory for Linear Operators}}.
\newblock Classics in Mathematics. Springer-Verlag, Berlin Heidelberg, 1995.

\bibitem{LandesmanLazer1970}
E.M. Landesman and A.C. Lazer.
\newblock {Nonlinear perturbations of linear elliptic boundary value problems
  at resonance}.
\newblock {\em J. Math. Mech.}, 19(7):609--623, 1970.

\bibitem{Oinas2007}
J.~Oinas.
\newblock {The degree theory and the index of a critical point for mappings of
  the type $({S}_+)$}.
\newblock {\em Acta Univ. Oul. A}, 488, 2007.

\bibitem{Phelps1997}
R.R. Phelps.
\newblock {Lectures on maximal monotone operators}.
\newblock {\em Extracta Math.}, 12(3):193--230, 1997.

\bibitem{Przeradzki1993}
B.~Przeradzki.
\newblock {Three methods for the study of semilinear equations at resonance}.
\newblock {\em Colloq. Math.}, LXVI:109--129, 1993.

\bibitem{ReedSimon1978}
M.~Reed and B.~Simon.
\newblock {\em {Methods of Mathematical Physics Vol. 4 - Analysis of
  Operators}}.
\newblock Academic Press, San Diego, 1978.

\bibitem{ReedSimon1980}
M.~Reed and B.~Simon.
\newblock {\em {Methods of Mathematical Physics Vol. 1 - Functional Analysis}}.
\newblock Academic Press, San Diego, 1980.

\bibitem{Skrypnik1994}
I.V. Skrypnik.
\newblock {\em {Methods for Analysis of Nonlinear Elliptic Boundary Value
  Problems}}, volume 139 of {\em Translations of Mathematical Monographs}.
\newblock American Mathematical Society, Providence, Rhode Island, 1994.

\bibitem{WillemZou2003}
M.~Willem and W.~Zou.
\newblock {On a Schr\"{o}dinger equation with periodic potential and spectrum
  point zero}.
\newblock {\em Indiana Univ. Math. J.}, 52(1):109--132, 2003.

\bibitem{YanCheDin2010}
M.~Yang, W.~Chen, and Y.~Ding.
\newblock {Solutions for periodic Schr\"{o}dinger equation with spectrum zero
  and general superlinear nonlinearities}.
\newblock {\em J. Math. Anal. Appl.}, 364(2):404--413, 2010.

\bibitem{Zeidler1986}
E.~Zeidler.
\newblock {\em {Nonlinear Functional Analysis and its Applications, Vol. I
  Fixed-Point Theorems}}.
\newblock Springer-Verlag, New York, 1986.

\end{thebibliography}

\end{document}